\documentclass[12pt]{article}
\usepackage{amssymb,latexsym}
\usepackage{amsfonts}
\font\germ=eufm10 \font\eur=eusm10 \font\bcl=cmbsy10
\def\qed{\hfill $\vrule height 2.5mm  width 2.5mm depth 0mm $}

\def\neweq{\setcounter{equation}{0}}
\newtheorem{theorem}{Theorem}[section]
\newtheorem{pr}[theorem]{Proposition}
\newtheorem{cor}[theorem]{Corollary}
\newtheorem{de}[theorem]{Definition}
\newtheorem{con}[theorem]{Conjecture}
\newtheorem{rem}[theorem]{Remark}
\newtheorem{lem}[theorem]{Lemma}

\newtheorem{prb}[theorem]{Problem}
\newtheorem{ex}[theorem]{Example}
\newtheorem{exs}[theorem]{Examples}
\def\g{\hbox{\germ g}}

\def\M{\hbox{\germ M}}

\def\K{\hbox{\germ K}}
\def\L{\hbox{\germ L}}
\def\R{\mathbb{R}}
\def\N{\mathbb{N}}
\def\Z{\mathbb{Z}}
\def\C{\mathbb{C}}
\def\Q{\mathbb{Q}}
\def\F{\mathbb{F}}
\def\Y{\mathbb{Y}}
\def\Ss{\mathbb{S}}

\def\Pe{\hbox{\eur P}}
\def\Pb{\hbox{\bcl P}}
\def\Rb{\hbox{\bcl R}}
\def\ds{\displaystyle}
\def\wt{\widetilde}
\def\wh{\widehat}
\def\ld{\lambda}
\def\ch{\rm ch}
\def\al{\alpha}
\def\goti#1{#1\llap{$#1$\hskip.3pt}\llap{$#1$\hskip.4pt}
   \llap{$#1$\hskip.5pt}}
\def\ggoti#1{\goti{\goti#1}}
\def\bnu{\ggoti\nu}
\def\bdoteq{\buildrel\bullet\over{=\!\!\!=}}
\begin{document}
\begin{center}
{\textbf{UBIQUITY OF KOSTKA POLYNOMIALS}}
\end{center}

\begin{center}
\textsc{anatol n.\ kirillov}
\end{center}
\begin{center}
{\small {\it Graduate School of Mathematics, Nagoya University}} \\
{\small {\it Chikusa--ku, Nagoya 486--8602, Japan}}\\ {\small
{\it and}}
 \\
{\small {\it Steklov Mathematical Institute,}} \\ {\small {\it
Fontanka 27, St.Petersburg, 191011, Russia}}
\end{center}

\vskip -1cm
\begin{abstract}
We report about results revolving around Kostka--Foulkes and
parabolic Kostka polynomials and their connections with
Representation Theory and Combinatorics. It appears that the set
of all parabolic Kostka polynomials forms a semigroup, which we
call {\it Liskova semigroup}. We show that polynomials frequently
appearing in Representation Theory and Combinatorics belong to
the Liskova semigroup. Among such polynomials we study
rectangular $q$--Catalan numbers; generalized exponents
polynomials; principal specializations of the internal product of
Schur functions; generalized $q$--Gaussian polynomials; parabolic
Kostant partition function and its $q$--analog; certain
generating functions on the set of transportation matrices. In
each case we apply rigged configurations technique to obtain some
interesting and new information about Kostka--Foulkes and
parabolic Kostka polynomials, Kostant partition function,
MacMahon, Gelfand--Tsetlin and Chan--Robbins polytopes. We
describe certain connections between generalized saturation and
Fulton's conjectures and parabolic Kostka polynomials; domino
tableaux and rigged configurations. We study also some properties
of $l$--restricted generalized exponents and the stable behaviour
of certain Kostka--Foulkes polynomials.
\end{abstract}
\vskip 0.5cm
\begin{center}
\textsc{Contents}
\end{center}
\nopagebreak
\begin{tabular}{rrl}
~~~~& 1 &  Introduction
\\
& 2 & Kostka--Foulkes polynomials and weight multiplicities
\\
&   & 2.1~~Kostka--Foulkes polynomials
\\
&   & 2.2~~Skew Kostka--Foulkes polynomials
\\
&   & 2.3~~Weight multiplicities
\\
& 3 & Parabolic Kostka polynomials
\\
& 4 & Fermionic formula for parabolic Kostka polynomials
\\
& 5 & Generalized exponents and mixed tensor representations
\\
& 6 & Internal product of Schur functions
\\
& 7 & Liskova semigroup
\\
&    & 7.1~~Realizable polynomials and product formula
\\
&    & 7.2~~Generalized $q$--Gaussian coefficients
\\
&    & 7.3~~Multinomial coefficients
\\
&    & 7.4~~Kostka--Foulkes polynomials and transportation
matrices
\\
&    & 7.5~~Gelfand--Tsetlin polytope and volume of weight
subspaces
\\
& 8 & Stable behavior of Kostka--Foulkes polynomials
\end{tabular}

\newpage

\section{Introduction}
\label{intro}
\neweq

The Kostka--Foulkes polynomials $K_{\ld\mu}(q)$ are defined as the matrix
elements of the transition matrix which expresses the Schur functions
$s_{\ld}(x)$ in terms of the Hall--Littlewood symmetric functions
$P_{\mu}(x;q)$:
\begin{equation}
s_{\ld}(x)=\sum_{\mu}K_{\ld\mu}(q)P_{\mu}(x;q). \label{1.1}
\end{equation}

To our knowledge, at the first time the polynomials
$K_{\ld\mu}(q)$ had appeared implicitly in the middle of fifties
in the familiar series of papers by Green \cite{Gr1,Gr2} about
computation of the characters of the general linear groups over
finite fields and their connections with the Hall--Steinitz
algebra. The famous Green polynomials $Q_{\rho}^{\ld}(q)$ are
related to the polynomials $K_{\ld\mu}(q)$ via the following
relation
$$q^{n(\mu)}Q_{\rho}^{\mu}(q^{-1})=\sum_{\ld}\chi_{\rho}^{\ld}K_{\ld\mu}(q),
$$
where $\chi_{\rho}^{\ld}$ is the value of the irreducible
character $\chi^{\ld}$ of the symmetric group $S_n$ at elements of
cycle--type $\rho$. Later in the paper \cite{Lt2} D.E.~Littlewood
introduced the symmetric functions
\begin{equation}Q_{\ld}'(x;q)=\prod_{1\le i<j\le n}(1-qR_{ij})^{-1}s_{\ld}(x),
\label{1.2}
\end{equation}
where $R_{ij}s_{\ld}(x):=s_{R_{ij}(\ld)}(x)$, and
$R_{ij}(\ld)=(\ld_1,\ldots ,\ld_i+1,\ldots ,\ld_j-1,\ldots
,\ld_n)$ is the so--called {\it raising} operator, see e.g.
\cite{Ma}, Chapter~I, \S 1. Here and further $s_{\ld}(x)$ denotes
the Schur function corresponding to a partition (or composition)
$\ld$, and to the set of variables $x_1,\ldots ,x_n$. In the same
paper D.E.~Littlewood also computed a few examples of the
decomposition of function $Q_{\ld}'(x;q)$ in terms of Schur
functions. The functions $Q_{\ld}'(x;q)$ are called {\it
modified} Hall--Littlewood polynomials and related with
polynomials $K_{\ld\mu}(q)$ via a simple relation
\begin{equation}
Q_{\mu}'(x;q)=\sum_{\ld}K_{\ld\mu}(q)s_{\ld}(x). \label{1.3}
\end{equation}
Polynomials $K_{\ld\mu}(q)$ also appeared implicitly, via the
generalized exponents of certain modules, in a fundamental paper
by Kostant \cite{K}.

Probably, A.O.~Morris \cite{Mo} was the first who introduced a special
notation $f_{\ld\mu}(t)$ for the transition coefficients between the
modified Hall--Littlewood polynomials $Q_{\ld}'(x;t)$ and Schur
functions $s_{\mu}(x)$, i.e. for Kostka--Foulkes polynomials. Using the
Littlewood formula (\ref{1.2}), A.O.~Morris discovered an important
recurrence relation between the polynomials $K_{\ld\mu}(q)$, and
computed all those polynomials 
with $|\ld|=|\mu|\le 4$.

In a survey paper \cite{Fo} H.O.~Foulkes at the first time, as far
as I am aware, gave rise the problem of non--negativity of the
coefficients of polynomials $K_{\ld\mu}(q)$. He stated a problem
"to associate each Young tableau with a precise index of $q$".
This problem was settled by A.~Lascoux and M.-P.~Sch\"utzenberger
\cite{LS,S} who suggested the name "Foulkes polynomials" for
polynomials $K_{\ld\mu}(q)$. Later they renamed these polynomials
by {\it Foulkes--Green} polynomials, see e.g. \cite{LS2}, and
C.R. Acad. Sci., Paris 288 (1979), 95-98.

In the middle of seventies the polynomials $K_{\ld\mu}(q)$ had appeared
in Algebraic Geometry and Representation Theory in a connection with the
study of the so--called "variety of the $N$--stable complete flags
${\mathcal B}_N$"; the crucial papers here are those by R.~Steinberg
(1976), R.~Hotta and T.~Springer (1977), T.~Springer (Inv. Math. 36
(1976), 173-207), H.~Kraft (1981), and C.~De~Concini and C.~Procesi
(1981). Precise references may be found in \cite{GP}. In a paper by
Garsia and Procesi \cite{GP} the polynomials $K_{\ld\mu}(q)$ are named
by $q$--Kostka polynomials.

There exists a vast body of literature about the Schur and
Hall--Littlewood functions, their numerous applications to
Representation theory of the symmetric and general linear groups
\cite{Ma,Sag,Mo2}, to Enumerative combinatorics \cite{St3},
Algebraic geometry \cite{Ful,GP}, to Integrable systems
\cite{Kac,KR,NY}, $\ldots$. But it is surprisingly enough that
till now there are no books or survey articles (except
\cite{DLT}) which would describe a big variety of applications of
Kostka--Foulkes polynomials to different areas of Mathematics and
Mathematical Physics.

The present paper is an attempt to collect together some results
about the Kostka--Foulkes and parabolic Kostka polynomials
related mainly to Representation theory and Combinatorics. Thus,
we almost omitted very interesting connections of polynomials
$K_{\ld\mu}(q)$ with Algebraic geometry, as well as with
Integrable systems. Some examples of such connections one may
find in \cite{GP,HKKOTY,Kir4,Kir5,LLT,KS,NY,ScW,SW,T}.

The main problems we are going to consider and study in the
present paper are {\it Identification} and {\it Characterization}
problems.

i) {\it Characterization problem}: how to characterize polynomials with
non--negative integer coefficients which might appear to be a
Kostka--Foulkes (or parabolic Kostka) polynomials $K_{\ld\mu}(q)$ (or
$K_{\ld R}(q)$) for some partitions $\ld$ and $\mu$ (or for some
partition $\ld$ and dominant sequence of rectangular shape partitions
$R$)?

If such an occasion happened, i.e.
$$P(q)\bdoteq K_{\ld R}(q),$$
one can apply {\it fermionic} formula for polynomials
$K_{\ld\mu}(q)$ (or $K_{\ld R}(q)$), see Theorem~\ref{t4.3} and
Corollary~\ref{c4.4},  to obtain a deep information about
polynomial $P(q)$. Typical example is the generalized binomial
coefficients $\left[\begin{array}{c}N\\ \ld\end{array}\right]_q$.
One can prove \cite{Kir3} that
$$\left[\begin{array}{c}N\\ \ld'\end{array}\right]_q\bdoteq K_{\wt\ld ,\wt\mu}(q),
$$
where $\wt\ld =(N|\ld|, \ld)$, $\wt\mu =(|\ld|^{N+1})$, and the
symbol $\bdoteq$ is explained in the list of special notation
({\bf Notatition}) at the end of Introduction. In this way one
comes to the following identity
\begin{equation}
\left[\begin{array}{c}N\\ \ld'\end{array}\right]_q=\sum_{\{\nu\}}
q^{c(\nu)}\prod_{j,k\ge 1}\left[\begin{array}{c}P_j^{(k)}(\nu)+
m_j(\nu^{(k)})\\ m_j(\nu^{(k)})\end{array}\right]_q, \label{1.4}
\end{equation}
summed over all sequences of partitions $\{\nu\}
=\{\nu^{(1)},\nu^{(2)},\ldots\}$ such that

$\bullet$ $|\nu^{(k)}|=\sum_{j\ge k}\ld_j$, $k\ge 1$,
$\nu^{(0)}:=\emptyset$;

$\bullet$ $P_j^{(k)}(\nu):=j(N+1)\delta_{k,1}+Q_j(\nu^{(k-1)})
-2Q_j(\nu^{(k)})+Q_j(\nu^{(k+1)})\ge 0$, for all $k,j\ge 1$;

$\bullet$ $m_j(\nu^{(k)})$ is the number of parts of partition
$\nu^{(k)}$ of size $j$;

$\bullet$ $c(\nu)=n(\nu^{(1)})-n(\ld)+\ds\sum_{k,j\ge 1}\left(
\begin{array}{c}(\nu^{(k)})_j'-(\nu^{(k+1)})_j'\\ 2\end{array}\right)$;

$\bullet$ $\left[\begin{array}{c}n\\ m\end{array}\right]_q=
\ds\frac{(q;q)_n}{(q;q)_m(q;q)_{n-m}}$, if $0\le m\le n$, and
equals to 0 otherwise.

\vskip 0.5cm
The identity (\ref{1.4}) is a highly non--trivial
identity, which generalizes the so--called KOH identity, see e.g.
\cite{Z,Kir3,Kir6}, or Section~\ref{ffpkp}, (\ref{4.4}). The only
known proof of (\ref{1.4}) is a combinatorial one and based on a
construction of the rigged configurations bijection, see
\cite{Kir3,Kir6}. For further generalization of the identity
(\ref{1.4}) see Section~\ref{ipsf}, (\ref{6.9}), or \cite{Kir6},
Section~10.1. It was shown [{\it ibid}] that "a constructive
proof of the unimodality of the Gaussian coefficients", which is
due to O'Hara \cite{OH}, admits a natural interpretation in terms
of rigged configurations. Note also, that identity (\ref{1.4})
leads to "a one--line high school algebra proof" of the
unimodality of the generalized Gaussian polynomials
$\left[\begin{array}{c}n\\ \ld\end{array}\right]_q$ (cf.
D.~Zeilberger, IMA volumes in Math. and Appl., 18 (1989), 67-75).

{\it Characterization problem} appears to be a rather difficult one, and
we don't know a complete list of conditions for a given polynomial
$P(q)\in\N [q]$ to be a parabolic Kostka polynomial. Examples show that
in general the sequence of coefficients of a  parabolic Kostka
polynomial $K_{\ld R}(q)$ may be not symmetric or unimodal, and may have
gaps, i.e. if
$$K_{\ld R}(q)=q^{\bullet}(a_0+a_1q+\cdots +a_mq^m),~~
a_0a_m\ne 0,$$
then it may be exist an index $j$, $1\le j\le m-1$,
such that $a_j=0$. Nevertheless, it seems plausible the following
\begin{con}\label{c1.1} Let \vskip -0.5cm
$$P(q)=\sum_{k=n}^ma_kq^k
$$
be a polynomial with non--negative integer coefficients such that
$a_na_m\ne 0$, $a_{n-1}=0$. If there exist non--negative integers
$i$ and $j$, $n\le i<j<m$, such that $j-i\ge 2$, and $a_{i-1}\le
a_i>a_{i+1}\ge\cdots\ge a_j<a_{j+1},$ then the polynomial $P(q)$
can't be equal to any parabolic Kostka polynomial.
\end{con}

ii) {\it Identification problem}: even though an answer to the
{\it Characterization problem} is still unknown, it is possible to
identify some polynomials, which frequently  appear in
Combinatorics and Representation theory, with certain parabolic
Kostka polynomials.

One of the main goals of our paper is to show that among polynomials,
which occur as parabolic Kostka polynomials, are the following:

$\bullet$ rectangular $q$--Catalan numbers, Section~\ref{kfpwm},
Exercise~1;

$\bullet$ generalized exponents polynomials, Theorem~\ref{t5.3};

$\bullet$ principal specialization of the internal product of Schur
functions, Theorem~\ref{t6.1};

$\bullet$ generalized $q$--binomial coefficients, Section~\ref{ggc};

$\bullet$ $q$--multinomial and $q^2$--multinomial coefficients,
Section~\ref{mc};

$\bullet$ $q$--analog of Kostant's partition function,
Section~\ref{kfs}, Exercise~2;

$\bullet$ certain generating functions on the set of
transportation matrices, Section~\ref{kfptm};

$\bullet$ polynomials $n(1+q)^m$, where $n\ge 1$ and $m\ge 2$ are
integer numbers, Section~\ref{ffpkp}, Exercise~6.

The main purpose of the present paper is to study combinatorial
properties of these polynomials, which follow from {\it identification}
of the polynomials above with certain parabolic Kostka polynomials, by
means of the rigged configurations bijection.

Let us say a few words about the content of our paper.

In Section~\ref{kfpwm} we recall definitions of Kostka--Foulkes
polynomials, weight multiplicities, and $q$--weight multiplicities.
Exercises to Section~\ref{kfpwm} include, among others, definitions and
basic properties of the rectangular $q$--Catalan and $q$--Narayana
numbers, computation of the volume and $\delta$--vector of MacMahon's
polytope.

In Section~\ref{pkp} we recall the definition of parabolic Kostka
polynomials which are natural generalization of the
Kostka--Foulkes polynomials. More details and proofs may be found
in \cite{KS}. Exercises to Section~\ref{pkp} include, among
others, a $q$--analog of a result by J.~Stembridge \cite{Stm}
about the number of up--down staircase tableaux, see Exercise~7;
connection between Littlewood--Richardson's numbers
$c_{\ld\mu}^{\nu}$ and certain parabolic Kostka polynomials, see
Exercise~3; various combinatorial properties of the convex
polytope
$${\cal P}_n=\{(x_1,\ldots ,x_n)\in\R_+^n|~x_i+x_{i+1}\le
1, 1\le i\le n-1\}, $$ including {\it fermionic} formulae for the
number $i({\cal P}_n;k)$ of integer points inside of the polytope
$k{\cal P}_n$; connections with {\it restricted} Kostka
polynomials and alternating permutations, see Exercise~4.

In Section~\ref{ffpkp} we recall a fermionic formula for parabolic
Kostka polynomials. Proof of the fermionic formula for parabolic
Kostka polynomials is based on the study of combinatorial
properties of the rigged configurations bijection, and may be
found in \cite{KSS}, see also \cite{Kir1}. Exercises to
Section~\ref{ffpkp} include, among others, an explanation of
certain connections between domino and $p$--ribbon tableaux on
one hand and parabolic Kostka polynomials and rigged
configurations on the other, see Exercises~3 and 4; "exotic
examples" of parabolic Kostka polynomials, which have a few
non--zero terms only. It is a very interesting and difficult
problem to find {\it all} parabolic Kostka polynomials with 2,3
or 4 non--zero terms.

In Section~\ref{gemtr} we study the generalized exponents
polynomials for some representations of the Lie algebra $sl(n)$.
The generalized exponents polynomials were introduced and studied
in the familiar paper by Kostant \cite{K} in a connection with the
investigation of the structure of the symmetric algebra of the
adjoint representation. Our main observation is that the
generalized exponents polynomial $G_N(V_{\al}\otimes V_{\beta}^*)$
corresponding to the tensor product $V_{\al}\otimes V_{\beta}^*$
of irreducible $\g l(N)$--representation $V_{\al}$ and the dual
$\g l(N)$--representation $V_{\beta}^*$, after multiplication by
some power of $q$, coincides with a certain parabolic Kostka
polynomial, see Theorem~\ref{t5.3}. This observation leads to a
fermionic formula for generalized exponents polynomial
corresponding to $\g l(N)$--representations of the form
$V_{\al}\otimes V_{\beta}^*$ and gives the first real means for
computing polynomials $G_N(V_{\al}\otimes V_{\beta}^*)$. We give
also a generalization of the Gupta--Hessenlink--Peterson theorem
(Theorem~\ref{t5.0}) to the case of so--called $l$--restricted
generalized exponents polynomial, see Definition~\ref{d5.9} and
Theorem~\ref{t5.12}.

In Section~\ref{ipsf} we study the internal product of Schur
functions and its principal specializations. It looks a
challenging problem to find and prove an analog of the
Lascoux--Sch\"utzenberger theorem about Kostka--Foulkes
polynomials (see e.g. \cite{Ma}, Chapter~III, (6.5), (i)) for the
transition coefficients $L_{\al\beta}^{(\mu)}(q)$ between the
internal product of Schur functions and Hall--Littlewood
polynomials
$$s_{\al}*s_{\beta}(x)=\sum_{\mu}L_{\al\beta}^{(\mu)}(q)
P_{\mu}(x;q),
$$
see Section~\ref{ipsf}, Problem~\ref{prb6.5}, for further details.

Polynomials $L_{\al\beta}^{(\mu)}(q)$ have many interesting
properties (see Proposition~\ref{p6.3} and Exercises {\bf 1} and
{\bf 2} to Section~\ref{ipsf}), and may be considered as a
natural generalization of the Kostka--Foulkes polynomials
$K_{\al\mu}(q)$. Note that numbers
$$L_{\al\beta}^{(\mu)}(0)=g_{\al\beta\mu}$$
are equal to the structural constants for
multiplication of the characters of the symmetric group
$S_{|\al|}$:
$$\chi^{\al}\chi^{\beta}=\sum_{\mu}g_{\al\beta\mu}\chi^{\mu}.
$$
Recall that it is an important open problem to obtain a nice
combinatorial interpretation of the numbers $g_{\al\beta\mu}$.
Based on an observation, see our Theorem~\ref{t6.1}, that the
principal specialization of a Schur function after multiplication
by some power of $q$, coincides with a certain parabolic Kostka
polynomial, we obtain a new proof, see Corollary~\ref{c6.7}, of a
fermionic/combinatorial formula for the principal specialization
of Schur functions first obtained in \cite{Kir2}.

In Section~\ref{kfs} we study {\it Identification} problem. Our
first result states that the set of all polynomials, which may
occur as a {\it parabolic} Kostka polynomial, forms a semigroup,
i.e. closed under multiplication. We call this semigroup by {\it
Liskova semigroup}. We show that Liskova semigroup contains the
generalized $q$--Gaussian polynomials
$\left[\begin{array}{c}n\\\ld'\end{array}\right]_q$,
$q$--multinomial coefficients $\left[\begin{array}{c}N\\
n_1,\ldots ,n_k\end{array}\right]_q$, $q$--analog of Kostant's
partition function, and polynomials related to the
Robinson--Schensted--Knuth correspondence between the set of
transportation matrices and that of pairs of semistandard Young
tableaux of the same shape. More precisely, let $\ld$ and $\mu$
be compositions whose lengths do not exceed $n$. Denote by
${\Pb}_{\ld\mu}$ the set of all $n$ by $n$ matrices of
non--negative integers with row sums $\ld_i$ and column sums
$\mu_j$ (the so--called set of {\it transportation matrices} of
type $(\ld ;\mu)$). The Robinson--Schensted--Knuth correspondence
establishes a bijection
\begin{equation}
{\Pb}_{\ld\mu}\simeq\coprod_{\eta}STY(\eta ,\mu)\times STY(\eta ,\ld)
\label{1.5}
\end{equation}
between the set of transportation matrices ${\Pb}_{\ld\mu}$ and that of
all pairs $(P,Q)$ of semistandard Young tableaux of the same shape and
weights $\mu$ and $\ld$ correspondingly.

Thus, for any transportation matrix $m\in{\Pb}_{\ld\mu}$ one can
associate the following statistics:

$\bullet$ left charge ~~$c_L(m)=c(P)$,

$\bullet$ right charge ~~$c_R(m)=c(Q)$,

$\bullet$ total charge ~~$c(m)=c(P)+c(Q)$.

It is well--known that the generating function
\begin{equation}
\sum_{m\in{\Pb}_{\ld\mu}}q^{c_R(m)}=\sum_{\eta}K_{\eta\mu}
K_{\eta\ld}(q):={\cal P}_{\ld\mu}(q) \label{1.6}
\end{equation}
describes\\
i) dimensions of cohomology groups of the partial unipotent flag
variety ${\cal F}_{\mu}^{\ld}$
$$q^{n(\ld)}{\cal P}_{\ld\mu}(q^{-1})=\sum_{k\ge 0}\dim H^{2k}({\cal
F}_{\mu}^{\ld};\Q )q^k,
$$
see e.g., Hotta~R. and Shimomura~N., Math. Ann. {\bf 241} (1979),
193-208, or \cite{LLT,Kir5};\\
ii) unrestricted one dimensional sum corresponding to the highest
weight $\ld$ and the tensor product of one row representations
$$V_{(\mu_1)}\otimes V_{(\mu_2)}\otimes\cdots\otimes V_{(\mu_l)},
$$
where $\mu =(\mu_1,\ldots ,\mu_l)$ and $l=l(\mu )$. For definition
and basic properties of one dimensional sums, their connections
with Kostka--Foulkes polynomials, with crystal basis, and
integrable models, see e.g. \cite{HKKOTY} and the literature
quoted therein;\\
iii) the number $\al_{\ld}(S;p)$ of chains of subgroups
$$\{ e\}\subseteq H^{(1)}\subseteq\cdots\subseteq H^{(m)}\subseteq G
$$
in a finite abelian $p$--group of type $\ld$ such that each subgroup
$H^{(i)}$ has order $p^{a_i}$, $1\le i\le m$:
$$p^{n(\ld)}\al_{\ld}(S;p^{-1})={\cal P}_{\ld\mu}(p),
$$
where $S=\{ a_1<a_2<\cdots <a_m\}$ is a subset of the interval
$[1,|\ld |-1]$, and
$$\mu :=\mu (S)=(a_1,a_2-a_1,\ldots ,a_m-a_{m-1},|\ld |-a_m),$$
see e.g. \cite{Bu}, or \cite{Kir5};\\
iv) the number of rational points ${\cal F}_{\mu}^{\ld}(\F_q)$
over the finite field $\F_q$ of the partial unipotent flag
variety ${\cal F}_{\mu}^{\ld}$:
$${\cal F}_{\mu}^{\ld}(\F_q)=q^{n(\ld)}{\cal P}_{\ld\mu}(q^{-1}).
$$

In particular, if $\ld =(1^n)$ then
$${\cal P}_{(1^n)\mu}(q)=q^{n(\mu')} \left[\begin{array}{c}n\\
\mu_1,\ldots ,\mu_n\end{array}\right],
$$
see e.g. (\ref{2.19}), and hence ${\cal P}_{(1^n)\mu}(q)$ belongs
to the Liskova semigroup.
 \vskip 0.3cm

\hskip -0.6cm{\bf Question.} Is it true that polynomials ${\cal
P}_{\ld\mu}(q)$ belong to the Liskova semigroup for all
partitions $\ld$ and compositions $\mu$?

\begin{con} Let $P(q)$ be a polynomial with non--negative integer
coefficients and non--zero constant term. There exist two
polynomials $f_1(q)$ and $f_2(q)$ both belonging to the Liskova
semifroup, such that $P(q)=f_1(q)/f_2(q)$.
\end{con}

For a "strong version" of this conjecture, see Section~\ref{rp}.

In Section~\ref{kfptm} we study the generating functions
\begin{eqnarray}
\sum_{m\in{\Pb}_{\ld\mu}}q^{c(m)}&=&\sum_{\eta}K_{\eta\mu}(q)
K_{\eta\ld}(q), \label{1.7}\\
 \sum_{m\in{\Rb}_{\ld\mu}}q^{c(m)}&=&\sum_{\eta}K_{\eta'\mu}(q)
 K_{\eta\ld}(q), \label{1.8}
\end{eqnarray}
where ${\Rb}_{\ld\mu}$ denotes the set of all $(0,1)$--matrices of size
$n$ by $n$ with row sums $\ld_i$ and column sums $\mu_j$.

It happens that the generating function (\ref{1.7}) coincides
with a certain Kostka--Foulkes polynomial, see Theorem~\ref{t7.8},
whereas the generating function (\ref{1.8}) coincides with a
certain "super" Kostka--Foulkes polynomial, see
Theorem~\ref{t7.9}.

In Section~\ref{vws} for any two partitions $\ld$ and $\mu$ of the
same size, we introduce {\it Ehrhart} polynomial ${\cal
E}_{\ld\mu}(t)$ of the weight subspace $V_{\ld}(\mu)$. Namely,
using the specialization $q=1$ of fermionic formula (\ref{4.1a})
for Kostka's polynomials, it is not difficult to see that the
Kostka number $K_{l\ld ,l\mu}$ is a {\it polynomial} in $l$ with
integer coefficients. Based on examples, we state a conjecture
that Ehrhart's polynomial ${\cal E}_{\ld\mu}(t)$ has in fact {\it
non--negative} integer coefficients. This conjecture would follow
from the general theory of Ehrhart's polynomials, see e.g.
\cite{Hi,St3}, if it would be known that the Gelfand--Tsetlin
polytope $G(\ld ,\mu)$, see Section~\ref{vws}, is an integral one
(this is the so--called Berenstein--Kirillov conjecture, see
\cite{KB}, p.101, Conjecture~2.1). More generally, it follows
from the $q=1$ case of the fermionic formula (\ref{4.1}) that for
any partition $\ld$ and any sequence of rectangular shape
partitions $R$ the parabolic Kostka number $K_{l\ld ,lR}$ is a
polynomial ${\cal E}_{\ld ,R}{(l)}$ in $l$ with integer
coefficients. Based on examples, we make a conjecture that the
polynomial ${\cal E}_{\ld ,R}(t)$ has in fact {\it non--negative}
integer coefficients. We expect that the similar statements are
valid for (generalized) parabolic Kostka numbers $K_{l\ld ,l\mu
,\eta}:=K_{l\ld,l\mu,\eta}(q)\vert_{q=1}$.

It looks a challenging problem to find an explicit formula for
the Ehrhart polynomials ${\cal E}_{\ld\mu}(t)$. This problem
should be a very difficult one, however, since, for example, the
polynomial ${\cal E}_{(n^n), ((n-1)^n,1^n)}(t)$ coincides with
the Ehrhart polynomial of the Birkhoff polytope ${\cal B}_n$, see
Example~\ref{e7.16}. The (normalized) leading coefficient of
polynomial ${\cal E}_{\ld\mu}(t)$ is equal to the (normalized)
volume of Gelfand--Tsetlin's polytope $G(\ld ,\mu)$, and is known
in the literature (see e.g., G.J.~Heckman, Inv. Math. {\bf 67}
(1982), 333-356) as a {\it continuous} analog of the weight
multiplicity $\dim V_{\ld}(\mu)$.

Finally, in the Exercises to Section~\ref{kfs} among others, we
study some properties of the parabolic Kostant partition function
$K_{\Phi (\eta)}(\gamma ;1)$, its $q$--analog, as well as its
connections with parabolic Kostka polynomials, and the volume of
so--called Chan--Robbins polytope, see Exercise~2. Kostant's
partition function was introduced and studied by F.A.~Berezin and
I.M.~Gelfand (Proc. Moscow Math. Soc. {\bf 5} (1956), 311-351)
for the case $\g =sl(n)$, and by B.~Kostant (Trans. Amer. Math.
Soc., 93 (1959), 53-73) for arbitrary semi--simple finite
dimensional Lie algebra $\g$. In \cite{Kir8} we had considered in
details algebraic and combinatorial properties of the Kostant
partition function $K_{(1^n)}(\gamma ;1)$, and, in particular,
computed the values of function $K_{(1^n)}(\gamma ;1)$ for
certain vectors $\gamma$, see Exercise~2 {\bf c, e}--{\bf g} to
Section~\ref{kfs}. It is not difficult to see (Exercise {\bf 2d}
to Section~\ref{kfs}) that if
$$\gamma =(\gamma_1,\ldots ,\gamma_n)\in\Z^n,~~ |\gamma|=0,
$$
then the value of Kostant's partition function $K_{(1^n)}(\gamma
;1)$ is equal to the number of $n$ by $n$ {\it skew--symmetric}
integer matrices $m=(m_{ij})_{1\le i,j\le n}$ with row sums
$\gamma_i$ and column sums $-\gamma_j$, such that $m_{ij}\ge 0$,
if $i\le j$. We denote by $SM(\gamma)$ the set of all such
skew--symmetric matrices. It is natural to ask does there exist
an analog of the Robinson--Schensted--Knuth correspondence for
the set of "skew--symmetric" transportation matrices
$SM(\gamma)$? We don't know any answer on this question in
general, but in a particular case when
$$\gamma =(d,d+1,\ldots ,d+n-1,-n(2d+n-1)/2),$$
we suggest such a bijection, see (\ref{3.22}), Exercise~2{\bf c}
to Section~\ref{kfs}. We show also (Exercise~2{\bf h} to
Section~\ref{kfs}) that the Ehrhart polynomial ${\cal E}(CR_n;t)$
of the Chan--Robbins polytope $CR_n$ is equal to the "continuous
analog" ${\cal P}_{\Phi(\eta)}(\beta)$ of the Kostant partition
function $K_{\Phi(\eta)}(\beta ,1)$. Namely,
$${\cal E}(CR_n;t)={\cal P}_{\Phi (1^{n+1})}(t,\underbrace{0,\ldots ,0}_{n-1},-t).
$$
Using a "fermionic formula" for the Kostant partition function
and its continuous analog in the {\it dominant chamber}
$Y_{n+1}^+$ obtained in \cite{Kir8}, we had computed the leading
coefficients of the Ehrhart polynomial ${\cal E}(CR_n;t)$, and
hence, the volume of Chan--Robbins' polytope $CR_n$. Formula for
the volume of Chan--Robbins polytope was conjectured in
\cite{CR}, and has been proved for the first time by
D.~Zeilberger \cite{Z2}. In Exercise~4 we state a generalized
saturation conjecture, and discuss connections between the
saturation conjecture (now a theorem by A.~Knutson and T.~Tao, J.
Amer. Math. Soc. {\bf 12} (1999), 1055--1090) and parabolic
Kostka polynomials. More precisely, for any two partitions $\ld$
and $\mu$, and a composition $\eta$ such that $|\eta |\ge l(\mu)$
we define numbers $a(\ld,\mu;\eta)$ and $b(\ld,\mu;\eta)$ via the
decomposition
$$K_{\ld\mu\eta}(q)=q^{a(\ld,\mu;\eta)}b(\ld,\mu;\eta)+{\rm
higher~degree~terms}.
$$
The numbers $a(\ld,\mu;\eta)$ and $b(\ld,\mu;\eta)$ seem to have
many interesting combinatorial properties. For example, the
Littlewood--Richardson numbers $c_{\ld,\mu}^{\nu}$ are a special
case of the numbers $b(\ld,\mu;\eta)$, see Section~\ref{pkp},
Example~3. In Example~4 (to Section~\ref{kfs}) we state several
conjectures about the numbers $a(\ld,\mu;\eta)$ and
$b(\ld,\mu;\eta)$. The basic two are:

$\bullet$ {\it generalized saturation conjecture:}

~~~$a(n\ld,n\mu;\eta)=na(\ld,\mu;\eta)$ for any positive integer
$n$;

$\bullet$ {\it generalized Fulton's conjecture:}

if $b(n\ld,n\mu;\eta)=1$ for some positive integer $n$, then
$b(N\ld,N\mu;\eta)=1$ for any positive integer $N$.

In Section~\ref{sb} we study the limiting behaviour of certain
Kostka--Foulkes polynomials. We give a generalization, see
Theorem~\ref{t8.1}, of a result obtained by Stanley \cite{St1},
and give answers on some questions had been posed in \cite{Kir6}.

Each Section, except for the first one, ends with a series of
exercises (which total number exceeds one hundred). The exercises
vary in difficulty: some are straightforward applications of the
material presented, while others are more difficult and represent
a brief exposition of original works or taken from various papers
written by the author.

The bibliographic notes are kept very brief. They are not
intended to give a comprehensive historical account or a complete
bibliography of the respective subject areas.

\vskip 0.3cm

\hskip -0.6cm{\bf Notation}

Throughout the paper we follow to Macdonald's book \cite{Ma} as for
notation related to the theory of symmetric functions, and Stanley's
book \cite{St3} as for notation related to Combinatorics. Below we give
a list of some special notation which we will use

1) if $P(q)$ and $Q(q)$ are polynomials in $q$, the symbol $P(q)\bdoteq
Q(q)$ means that the ratio $P(q)/Q(q)$ is a power of $q$;

2) if $a,k_0,\ldots ,k_m$ are (non--negative) integers, the symbol
$q^a(k_0,\ldots ,k_m)$ stands for polynomial $\sum_{j=0}^mk_jq^{a+j}$;

3) (Garsia's symbol $\chi (P)$) if $P$ is any statement, then $\chi
(P)=1$ if it is true, and $\chi (P)=0$ otherwise;

4) if $\ld$ is a partition and $\mu$ is a composition, the symbol
$STY(\ld ,\mu)$ denotes the set of {\it semistandard} Young
tableaux of shape $\ld$ and weight (or content) $\mu$.

\vskip 0.5cm

\section{Kostka--Foulkes polynomials and weight multiplicities}
\label{kfpwm}
\neweq

\subsection{Kostka--Foulkes polynomials}
\label{kfp}

Kostka--Foulkes polynomials are defined as the matrix elements of
the transition matrix
$$K(q)=M(s,P)
$$ from the Schur functions $s_{\lambda}(x)$ to the Hall--Littlewood
functions $P_{\mu}(x;q)$:
\begin{equation}
s_{\lambda}(x)=\sum_{\mu}K_{\lambda\mu}(q)P_{\mu}(x;q).\label{2.1}
\end{equation}
It is well known \cite{Ma}, Chapter~I, that if $\lambda$ and $\mu$
partitions, then

$\bullet$ $K_{\lambda\mu}(q)\ne 0$ if and only if $\lambda\ge\mu$ with
respect to the dominance partial ordering "$\ge$" on the set of partitions
${\cal P}_n$:
\begin{eqnarray*}
&&\lambda\ge\mu ~~{\rm if ~and ~only ~if}\\ &&i)~~|\lambda|=|\mu| ,\\
&&ii)~~\lambda_1+\cdots +\lambda_i\ge\mu_1+\cdots +\mu_i ~~{\rm for
~all}~~ i\ge 1.
\end{eqnarray*}

$\bullet$ If $\lambda\ge\mu$, $K_{\lambda\mu}(q)$ is a monic of degree
$n(\mu)-n(\lambda)$ polynomial with non--negative integer coefficients.

Let us recall that for  any partition $\lambda
=(\lambda_1,\lambda_2,\cdots ,\lambda_p)$,
$$ n(\lambda )=\sum_{i=1}^p(i-1)\lambda_i=\sum_{1\le i<j\le p}\min
(\lambda_i,\lambda_j).
$$

It was H.O.~Foulkes \cite{Fo} who asked for a combinatorial explanation of the
positivity of coefficients of polynomials $K_{\lambda\mu}(q)$. He
conjectured the existence of natural statistics $c(T)$ on the set of
semistandard Young tableaux
$STY(\lambda ,\mu)$ of shape $\lambda$ and weight $\mu$ such that
\begin{equation}
K_{\lambda\mu}(q)=\sum_{T\in STY(\lambda ,\mu)}q^{c(T)}.\label{2.2}
\end{equation}
This conjecture was settled by Lascoux and Sch\"utzenberger \cite{LS,S},
who identified $c(T)$ to be the rank of a poset structure on the set
$STY(\bullet ,\mu)$ of semistandard Young tableaux of fixed weight
$\mu$. Another proof of the positivity of coefficients of polynomials
$K_{\ld\mu}(q)$ has been obtained by Lusztig \cite{Lu2}, who identified
the polynomials
$${\wt K}_{\ld\mu}(q):=q^{n(\mu)-n(\ld)}K_{\ld\mu}(q^{-1})$$
with the certain Kazhdan--Lusztig polynomials
$P_{n_{\ld},n_{\mu}}(q)$ for the affine Weyl group $W(\wh{\g}
l_n)$.

\subsection{Skew Kostka--Foulkes polynomials}
\label{skfp}

Let $\ld ,\mu$ and $\nu$ be partitions, $\ld\supset\mu$, and
$|\ld|=|\mu|+|\nu|$.

\begin{de}\label{d1.1} The skew Kostka--Foulkes polynomials
$K_{\ld\setminus\mu ,\nu}(q)$ are defined as the transition coefficients
from the skew Schur functions $s_{\ld\setminus\mu}(x)$ to
the Hall--Littlewood functions $P_{\nu}(x;q)$:
\begin{equation}
s_{\ld\setminus\mu}(x)=\sum_{\nu}K_{\ld\setminus\mu ,\nu}(q)P_{\nu}(x;q).
\label{2.3}
\end{equation}
\end{de}

It is clear that
$$K_{\ld\setminus\mu ,\nu}(q)=\sum_{\pi}c_{\mu\pi}^{\ld}K_{\pi\nu}(q),$$
where the coefficients $c_{\mu\pi}^{\nu}={\rm
Mult}[V_{\nu}:V_{\mu}\otimes V_{\pi}]$ stand for the
Littlewood--Richardson numbers.

Let us remark that 
\begin{equation}
K_{\ld\setminus\mu ,\nu}(q)=\sum_Tq^{c(T)} \label{2.4}
\end{equation}
summed over all semistandard skew tableaux $T$ of shape $\ld\setminus\mu$
and weight $\nu$, where $c(T)$ denotes the charge of skew tableau $T$.

We define also {\it cocharge} version of the skew Kostka--Foulkes
polynomials:
\begin{equation}
{\overline K}_{\ld\setminus\mu ,\nu}(q)=\sum_{\pi}c_{\mu\pi}^{\ld}
{\overline K}_{\pi\mu}(q), \label{2.5}
\end{equation}
where ${\overline K}_{\ld\mu}(q)=q^{n(\mu)}K_{\ld\mu}(q^{-1})$.

\subsection{Weight multiplicities}
\label{wm}

Let $V$ be a finite dimensional representation of the Lie algebra $\g
l(n)$, and $\mu$ be an integrable weight (i.e. composition), denote by
$V(\mu)$ the weight $\mu$ subspace of the representation $V$. Dimension
${\rm dim}V(\mu)$ of the space $V(\mu)$ is called {\it weight
multiplicity} of $\mu$ in $V$. Let ${\rm ch}V$ denote the character of
representation $V$. It is a symmetric polynomial in variables
$X_n=(x_1,\ldots ,x_n)$. If $P$ is any polynomial in the variables $X_n$, let
$[x^{\mu}]P$ be the coefficient of $x^{\mu}=x_1^{\mu_1}\cdots
x_n^{\mu_n}$ in $P$. It is clear that
$${\rm dim}V(\mu)=[x^{\mu}]{\rm ch}V.
$$

Follow S.~Kato~\cite{Kato}, define a $q$--analog ${\rm dim}_qV(\mu)$ of the
weight multiplicity ${\rm dim}V(\mu)$ as the transition coefficient from
the character ${\rm ch}V$ to the Hall--Littlewood polynomials
$P_{\mu}(X_n;q)$:
\begin{equation}
{\rm ch}V=\sum_{\mu}{\rm dim}_qV(\mu)P_{\mu}(X_n;q). \label{2.6}
\end{equation}
In the present paper we are mainly interested in the following two cases:

i) $V=V_{\ld}$ is an irreducible highest weight $\ld$ representation of
the Lie algebra $\g l(n)$.

In this case the character ${\rm ch}V_{\ld}$
of representation $V_{\ld}$ coincides with the Schur function
$s_{\ld}(X_n)$, and dimension ${\rm dim}V_{\ld}(\mu)$ of the weight $\mu$
subspace of $V_{\ld}$ is equal to the Kostka number
$$K_{\ld\mu}:=K_{\ld\mu}(1),$$
and is equal also to the number of semistandard Young tableaux of
shape $\ld$ and weight $\mu$. It follows from (\ref{2.6}) and
(\ref{2.1}) that the $q$--analog ${\rm dim}_qV_{\ld}(\mu)$ of
weight multiplicity ${\rm dim}V_{\ld}(\mu)$ is equal to the
Kostka--Foulkes polynomial $K_{\ld\mu}(q)$. Polynomials
$K_{\ld\mu}(q)$ and the Lusztig $q$--analog of weight multiplicity
$d_{\mu}(L_{\ld};q)$ in the case of type $A$, see \cite{Ls},
p.215, (6.10), are connected by the following relation
$$d_{\mu}(L_{\ld};q)=q^{n(\mu)-n(\ld)}K_{\ld\mu}(q^{-1}).
$$
In the sequel, we will denote the Lusztig polynomials $d_{\mu}(L_{\ld};q)$ by
$\wt K_{\ld\mu}(q)$, so that
\begin{equation}
\wt K_{\ld\mu}(q)=q^{n(\mu)-n(\ld)}K_{\ld\mu}(q^{-1})=
q^{-n(\ld)}{\overline K}_{\ld\mu}(q). \label{2.7}
\end{equation}

ii) $V=V_{(\mu)}=V_{(\mu_1)}\otimes\cdots\otimes V_{(\mu_r)}$ is the
tensor product of irreducible representations $V_{(\mu_a)}$
corresponding to the one row partitions $(\mu_a)$, $1\le a\le r$, and
$\mu =(\mu_1,\ldots ,\mu_r)$.

In this case the character ${\rm ch}V_{(\mu)}$ is equal to the product
$\ds\prod_{j=1}^rh_{\mu_j}(X_n)$ of complete homogeneous symmetric
polynomials $h_{\mu_j}(X_n)$, $1\le j\le r$. The weight $\ld$
multiplicity ${\rm dim}V_{(\mu)}(\ld)$ is equal to the number ${\cal
P}_{\ld\mu}$ of $n$ by $n$ matrices of non--negative integers with row
sums $\ld_i$ and column sums $\mu_j$:
\begin{equation}
{\rm dim}V_{(\mu)}(\ld)=\sum_{\eta}K_{\eta\mu}K_{\eta\ld}={\cal
P}_{\ld\mu}. \label{2.8}
\end{equation}
We denote by ${\Pb}_{\ld\mu}$ the set of all $n$ by $n$ matrices of
non--negative integers with row sums $\ld_i$ and column sums $\mu_j$
(the set of so--called {\it transportation matrices} of type $(\ld
;\mu)$). There exists a remarkable bijection, known as the
Robinson--Schensted--Knuth correspondence,
\begin{equation}
{\Pb}_{\ld\mu}\simeq\coprod_{\eta}STY(\eta ,\mu)\times STY(\eta ,\ld)
\label{2.9}
\end{equation}
between the set of transportation matrices of type $(\ld ;\mu)$ and that
of all pairs $(P,Q)$ of semistandard Young tableaux of the same shape
and weights $\mu$ and $\ld$, correspondingly. The
Robinson--Schensted--Knuth correspondence (\ref{2.9}) leads to a
combinatorial/bijective proof of the last equality in (\ref{2.8}).

Finally, let us find the $q$--analog, denoted by ${\rm
dim}_qV_{(\mu)}(\ld)$, of multiplicity of the weight $\ld$ in
representation $V_{(\mu )}$. Using (\ref{2.6}), (\ref{2.1}) and
the product formula for the character of representation $V_{(\mu
)}$:
$${\rm ch}V_{(\mu )}=\ds\prod_{j=1}^r h_{\mu_j}(X_n),$$
one can obtain that
\begin{equation}
{\rm dim}_qV_{(\mu)}=\sum_{\eta}K_{\eta\mu}K_{\eta\ld}(q):={\cal
P}_{\ld\mu}(q). \label{2.10}
\end{equation}
Polynomials ${\cal P}_{\ld\mu}(q)$ admit numerous algebraic,
algebro--geometric and combinatorial interpretations. Some of
them were mentioned in Introduction. Here we only mention a
result by R.~Hotta and N.~Shimomura (Math. Ann. {\bf 24} (1979),
193-198) that the product $t^{n(\ld)}{\cal P}_{\ld\mu}(t^{-1})$
is equal to the Poincare polynomial of the partial unipotent flag
variety ${\cal F}_{\mu}^{\ld}$:
$$t^{n(\ld)}{\cal P}_{\ld\mu}(t^{-1})=\sum_{k\ge 0}{\rm dim}H^{2k}({\cal
F}_{\mu}^{\ld},\Q )t^k,
$$
see \cite{Kir5,LLT,T} and references therein.

\vskip 0.3cm

\hskip -0.6cm {\bf Exercises}
\vskip 0.2cm

\hskip -0.6cm{\bf 1.} {\it Rectangular Catalan and Narayana polynomials,
and MacMahon polytope }

{\bf a.} Define {\it rectangular Catalan} polynomial
\begin{equation}
C(n,m|q)=\frac{(q;q)_{nm}}{\ds\prod_{i=1}^n\prod_{j=1}^m(1-q^{i+j-1})}.
\label{2.11}
\end{equation}
Show that
\begin{equation}
q^{m\scriptsize{\left(\begin{array}{c}n\\ 2\end{array}\right)}}C(n,m|q)=
K_{(n^m),(1^{nm})}(q). \label{2.12}
\end{equation}
Thus, $C(n,m|q)$ is a polynomial of degree $nm(n-1)(m-1)/2$ in the
variable $q$ with non--negative integer coefficients. Moreover,
$$C(n,2|q)=C(2,n|q)=c_n(q)=\frac{1-q}{1-q^{n+1}}\left[\begin{array}{c}
2n\\ n\end{array}\right]_q
$$
coincides with "the most obvious" $q$--analog of the Catalan numbers,
see e.g. \cite{Fur}, p.255, or \cite{St3}.

{\bf b.} It follows from (\ref{2.12}) that the rectangular Catalan
number $C(n,m|1)$ counts the number of {\it lattice} words
$$w=a_1a_2\cdots a_{nm}$$
of weight $(m^n),$ i.e. lattice words in which each $i$ between 1
and $m$ occurs exactly $n$ times. Let us recall that a word
$a_1\cdots a_p$ in the symbols $1,\ldots ,m$ is said to be a {\it
lattice} word, if for $1\le r\le p$ and $1\le j\le m-1$, the
number of occurrences of the symbol $j$ in $a_1\cdots a_r$ is not
less than the number of occurrences of $j+1$:
\begin{equation}
\#\{ i|1\le i\le r~~{\rm and}~~a_i=j\}\ge\#\{ i|1\le i\le r~~{\rm
and}~~a_i=j+1\}. \label{2.12a}
\end{equation}

For any word $w=a_1\cdots a_k$, in which each $a_i$ is a positive
integer, define the major index
$${\rm maj}(w)=\sum_{i=1}^{k-1}i\chi (a_i>a_{i+1}),
$$
and the number of descents~~~
$${\rm des}(w)=\ds\sum_{i=1}^{k-1}\chi (a_i>a_{i+1}).$$

Finally, for any integer $k$ between 0 and $(n-1)(m-1)$, define {\it
rectangular $q$--Narayana} number
$$N(n,m;k|q)=\sum_wq^{{\rm maj}(w)},
$$
where $w$ ranges over all lattice words of weight $(m^n)$ such that ${\rm
des}(w)=k$.

\begin{ex} {\rm Take $n=4$, $m=3$, then
$$\sum_{k=0}^6N(3,4;k|1)t^k=1+22t+113t^2+190t^3+113t^4+22t^5+t^6.
$$
This example shows that if $n,m\ge 3$, it is unlikely that there exists a
simple combinatorial formula for the rectangular Narayana numbers
$N(n,m;k|1)$, but see formula (\ref{2.12b}).}
\end{ex}

$\bullet$ Show that

$i)$ $C(n,m|q)=\ds\sum_wq^{{\rm maj}(w)}$, where $w$
ranges over all lattice words of weight $(m^n)$;

$ii)$  $N(n,m;k|q)=q^{nm((n-1)(m-1)/2-k)}N(n,m;(n-1)
(m-1)-k|q)=N(m,n;k|q),$ for any integer $k$, ~$0\le k\le
(n-1)(m-1)/2$;

$iii)$ $N(2,n;k|q)=q^{k(k+1)}\ds\frac{1-q}{1-q^n}\left[\begin{array}{c}
n\\ k\end{array}\right]_q\left[\begin{array}{c}n\\
k+1\end{array}\right]_q \bdoteq {\rm dim}_qV_{(k,k)}^{\g
l(n-k+1)}$, for any integer $k$, ~$0\le k\le n-1$, where
$V_{(k,k)}^{{\g}l(n-k+1)}$ stands for the irreducible
representation of the Lie algebra ${\g} l(n-k+1)$ corresponding
to the two row partition $(k,k)$; recall that for any finite
dimensional $\g l(N)$--module $V$ the symbol ${\rm dim}_qV$
denotes its $q$--dimension, i.e. the principal specialization of
the character of the module $V$: \vskip -0.5cm
$${\rm dim}_qV=({\rm ch}V)(1,q,\ldots ,q^{N-1});
$$
\vskip 0.3cm
$iv)$ $N(n,m;1|1)=\ds\sum_{j\ge 2}\left(\begin{array}{c}n\\
j\end{array}\right)\left(\begin{array}{c}m\\ j\end{array}\right)
=\left(\begin{array}{c}n+m\\ n\end{array}\right)-nm-1$;

$v)$ (Fermionic formula for $q$--Narayana numbers)

\begin{equation}
q^{m\scriptsize\left(\begin{array}{c}n\\ 2\end{array}\right)}
N(n,m;l|q)=\sum_{\{\nu\}}q^{c(\nu)}\prod_{k,j\ge 1}\left[
\begin{array}{c}P_j^{(k)}(\nu)+m_j(\nu^{(k)})\\ m_j(\nu^{(k)})\end{array}
\right]_q, \label{2.12b}
\end{equation}
summed over all sequences of partitions $\{\nu\} =\{\nu^{(1)},\nu^{(2)},
\ldots ,\nu^{(m-1)}\}$ such that

$\bullet$ $|\nu^{(k)}|=(m-k)n$, $1\le k\le m-1$;

$\bullet$ $(\nu^{(1)})_1'=(m-1)n-l$, i.e. the length of the first column
of the diagram $\nu^{(1)}$ is equal to $(m-1)n-l$;

$\bullet$ $P_j^{(k)}(\nu):=Q_j(\nu^{(k-1)})-2Q_j(\nu^{(k)})+
Q_j(\nu^{(k+1)})\ge 0$, for all $k,j\ge 1$,\\
where by definition we put
$\nu^{(0)}=(1^{nm})$; for any diagram $\ld$ the number
$Q_j(\ld)=\ld_1'+\cdots\ld_j'$ is equal to the number of cells in the
first $j$ columns of the diagram $\ld$, and $m_j(\ld)$ is equal to the
number of parts of $\ld$ of size $j$;

$\bullet$ $c(\nu)=\ds\sum_{k,j\ge
1}\left(\begin{array}{c}(\nu^{(k-1)})_j' -(\nu^{(k)})_j'\\
2\end{array}\right)$.

\begin{con} If $1\le k\le (n-1)(m-1)/2$, then
$$N(n,m;k-1|1)\le N(n,m;k|1),
$$
i.e. the sequence of rectangular Narayana numbers
$\{N(n,m;k|1)\}_{k=0}^{(n-1)(m-1)}$ is symmetric and unimodal.
\end{con}

For definition of unimodal polynomials/sequences see e.g.
\cite{St2}, where one may find a big variety of examples of
unimodal sequences which frequently appear in Algebra,
Combinatorics and Geometry.

{\bf c.} (Volume of the MacMahon polytope) Let ${\M}_{mn}$ be the convex
polytope in $\R^{nm}$ of all points ${\mathbf x}=(x_{ij})_{1\le i\le n,
1\le j\le m}$ satisfying the following conditions
\begin{equation}
0\le x_{ij}\le 1, ~~~x_{ij}\ge x_{i-1,j},~~~x_{ij}\ge x_{i,j-1},
\label{2.14a}
\end{equation}
for all pairs of integers $(i,j)$ such that $1\le i\le n$, $1\le
j\le m$, and where by definition we set $x_{i0}=0=x_{0j}$.

We will call the polytope $\M_{nm}$ by {\it MacMahon polytope}.
The MacMahon polytope is an integral polytope of dimension $nm$
with $\left(\begin{array}{c}m+n\\ n\end{array}\right)$ vertices
which correspond to the set of (0,1)--matrices satisfying
(\ref{2.14a}).

If $k$ is a positive integer, define $i(\M_{nm};k)$ to be the number of
points ${\bf x}\in\M_{nm}$ such that $k{\mathbf x}\in\Z^{nm}$. Thus,
$i(\M_{nm};k)$ is equal to the number of plane partitions of rectangular
shape $(n^m)$ with all parts do not exceed $k$.  By a theorem of
MacMahon (see e.g. \cite{Ma}, Chapter~I, \S 5, Example~13)
\begin{equation}
i(\M_{nm};k)=\prod_{i=1}^n\prod_{j=1}^m\frac{k+i+j-1}{i+j-1}.
\label{2.14b}
\end{equation}
It follows from (\ref{2.14b}) that the Ehrhart polynomial ${\cal
E}(\M_{nm};t)$ of the MacMahon polytope $\M_{nm}$ is completely resolved
into linear factors:
$${\cal E}(\M_{nm};t)=\prod_{i=1}^n\prod_{j=1}^m\frac{t+i+j-1}{i+j-1}.
$$
Hence, the normalized volume
$${\wt{\rm vol}}(\M_{nm})=(nm)!{\rm vol}(\M_{nm})$$
of the MacMahon polytope $\M_{nm}$ is equal to the rectangular
Catalan number $C(n,m|1)$, i.e. the number of standard Young
tableaux of  rectangular shape $(n^m)$. We refer the reader to
\cite{St3}, Section~4.6, and \cite{Hi}, Chapter~IX, for
definition and basic properties of the Ehrhart polynomial ${\cal
E}({\Pe};t)$ of a convex integral polytope ${\Pe}$.

$\bullet$ Show that
\begin{equation}
\sum_{k\ge0}i(\M_{nm};k)z^k=\left(\sum_{j=0}^{(n-1)(m-1)}N(n,m;j)
z^j\right)/(1-z)^{nm+1}, \label{2.14c}
\end{equation}
where
$$N(n,m;j):=N(n,m;j|1)$$
denotes the rectangular Narayana number. Thus, the sequence of
Narayana numbers
$$(1=N(n,m;0),N(n,m;1),\ldots ,N(n,m;(n-1)(m-1))=1)
$$
is the {\it $\delta$--vector} (see e.g. \cite{Hi}, \S 34) of the
MacMahon polytope. In the case $n=2$ (or $m=2$) all these results
may be found in \cite{St3}, Chapter~6, Exercise~6.31.\\ {\bf
Questions.} i) ({\it Higher associahedron}) Does there exist an
$(m-1)(n-1)$--dimensional integral convex (simplicial?) polytope
$Q_{n,m}$ which has $h$--vector
$$h=(h_0(Q_{n,m}),h_1(Q_{n,m}),\ldots ,h_{(n-1)(m-1)}(Q_{n,m}))$$
given by
$$\sum_{i=0}^{(n-1)(m-1)}h_i(Q_{n,m})q^{(n-1)(m-1)-i}=C(n,m|q)~?
$$
We refer the reader to \cite{Hi}, Chapter~I, \S 6 and Chapter~III,
for definitions and basic properties of the $h$--vector of a
simplicial polytope; see also, R.~Stanley (J. Pure and Appl.
Algebra {\bf 71} (1991), 319-331).

ii) Define {\it rectangular} Schr\"oder polynomial
$$S(n,m|q):=C(n,m|1+q),$$
and put
$$S(n,m|q)=\ds\sum_{k\ge0}^{(n-1)(m-1)}S(n,m||k)q^k.
$$
What is a combinatorial interpretation(s) of the numbers
$S(n,m||k)$ and $S(n,m|1)$?

The answers on these questions are known if either $n$ or $m$ is
equal to 2, see e.g. R.~Simion (Adv. in Appl. Math. {\bf 18}
(1997), 149-180, Example~{\bf 4} (the Associahedron)). Note
finally, that MacMahon's polytope and Narayana's  numbers will
appear again in Section~\ref{kfs}, Exercise~3{\bf e}.

{\bf d.} Denote by $f(n,d)$ the number of ways to draw $d$ diagonals in
a convex $(n+2)$--gon, such that no two diagonals intersect in their
interior.

$i)$ Show that
\begin{equation}
f(n,d)=\frac{1}{n+1}\left(\begin{array}{c}n+d+1\\ d+1\end{array}\right)
\left(\begin{array}{c}n-d\\
d\end{array}\right)=K_{((d+1)^2,1^{n-d-1}),(1^{n+d+1})}(1).
\label{2.14d}
\end{equation}
For example, $f(n,n-1)$ is just the Catalan number
$$C_n=\ds\frac{1}{n+1}\left(\begin{array}{c}2n\\ n\end{array}\right).$$
The first equality in the formula (\ref{2.14d}) goes back to
T.K.~Kirkman (1857), E.~Prouhet (1866), and A.~Cayley (1890); for
precise references, see \cite{St3}, p.272.

It follows from (\ref{2.14d}) that $f(n,d)$ is equal to the number of
standard Young tableaux of shape $((d+1)^2,1^{n-d-1})$. It is natural to
ask for a bijection between the polygon  dissections and standard Young
tableaux. Such a bijection was constructed by R.~Stanley (J. Comb.
Theory A, {\bf 76} (1996), 175-177).

$ii)$ Define a statistics $b$ on the set $F(n,d)$ of $d$--tuples of
non--crossing in their interior diagonals in a convex $(n+2)$--gon such
that
$$\sum_{f\in F(n,d)}q^{b(f)}=\wt K_{((d+1)^2,1^{n-d-1}),(1^{n+d+1})}(q).
$$

\hskip -0.6cm{\bf 2.} Let $\mu$ be a composition, $|\mu|=n$.

{\bf a.} Show that
\begin{equation}
{\cal P}_{(1^n),\mu}(q)=\sum_{\eta}K_{\eta\mu}K_{\eta ,(1^n)}(q)=
q^{n(\mu')}\left[\begin{array}{c}n\\ \mu_1,\ldots
,\mu_n\end{array}\right]_q. \label{2.13}
\end{equation}

Formula (\ref{2.13}) goes back to C.~Ehresmann who had computed (Ann.
Math. 35 (1934), 396-443) the Poincare polynomial of the partial flag
variety ${\cal F}_{\mu}$ corresponding to a composition $\mu$ of $n$:
$$t^{\scriptsize\left(\begin{array}{c}n\\ 2\end{array}\right)}
{\cal P}_{(1^n)\mu}(t^{-1})=\sum_{k\ge 0}{\rm dim}H^{2k}({\cal F}_{\mu},
\Q )t^k=\left[\begin{array}{c}n\\ \mu_1,\ldots ,\mu_n\end{array}\right]_t.
$$

{\bf b.} (A generalization of the previous Exercise) Let
$R=\{(\mu_a^{\eta_a})\}^l_{a=1}$ be a sequence of rectangular shape
partitions, $N=\sum_{a=1}^l\mu_a\eta_a$. Show that
\begin{equation}
\sum_{\ld}K_{\ld R}K_{\ld ,(1^N)}(q)=\prod_{a=1}^l
K_{(\mu_a^{\eta_a}),(1^{\mu_a\eta_a})}(q)\left[\begin{array}{c}
N\\ \mu_1\eta_1,\ldots ,\mu_l\eta_l\end{array}\right]_q,
\label{2.14}
\end{equation}
where $K_{\ld R}$ denotes the value of {\it parabolic Kostka
polynomial} $K_{\ld R}(q)$ at $q=1$, see Section~\ref{pkp},
Definition~\ref{d3.6}. In particular, if sequence
$$R=\{\mu_1,\ldots ,\mu_r,1^{\eta_1},\ldots ,1^{\eta_s}\}
$$
consists of only one row or one column partitions, and
$$\ds\sum_{j=1}^r\mu_j+\ds\sum_{j=1}^s\eta_j=N,$$
then
\begin{equation}
\sum_{\ld}K_{\ld R}K_{\ld
,(1^N)}(q)=q^{n(\mu')}\left[\begin{array}{c} N
\\ \mu_1,\ldots ,\mu_r,\eta_1,\ldots ,\eta_s\end{array}\right]_q.
\label{2.15}
\end{equation}

\hskip -0.6cm{\bf 3.} Define  polynomial $\varphi_{\nu\mu}^{\ld}(q)$ to
be the coefficient of $P_{\ld}$ in the product $s_{\nu}P_{\mu}$:
\begin{equation}
s_{\nu}(x)P_{\mu}(x;q)=\sum_{\ld}\varphi_{\nu\mu}^{\ld}(q)P_{\ld}(x;q).
\label{2.16}
\end{equation}
Equivalently,
$$\varphi_{\nu\mu}^{\ld}(q)=\sum_{\eta}K_{\nu\eta}(q)f_{\eta\mu}^{\ld}(q),
$$
where for any three partitions $\ld ,\mu ,\nu$ the polynomial
$f_{\mu\nu}^{\ld}(q)$ is the coefficient of $P_{\ld}$ in the
product $P_{\mu}P_{\nu}$, see e.g. \cite{Ma}, Chapter~III, \S 3.
Recall that $P_{\ld},P_{\mu},\ldots$ denote the Hall--Littlewood
polynomials corresponding to partitions $\ld,\mu,\ldots$.

\begin{con} For any three partitions $\ld ,\mu ,\nu$ the polynomial
$\varphi_{\nu\mu}^{\ld}(q)$ is a polynomial with non--negative integer
coefficients.
\end{con}

For example:
\begin{equation}
\varphi_{(1^n)\mu}^{\ld}(q)=f_{(1^n)\mu}^{\ld}(q)=\prod_{i\ge 1}
\left[\begin{array}{c}\ld_i'-\ld_{i+1}'\\
\ld_i'-\mu_i'\end{array}\right]_q, \label{2.17}
\end{equation}
and therefore $\varphi_{(1^n)\mu}^{\ld}(q)=0$, unless $\ld\setminus\mu$
is a vertical $n$--strip, \cite{Ma}, Chapter~III, (3.2), p.215;

\begin{equation}
\varphi_{(n)\mu}^{\ld}(q)=q^{c(\ld ;\mu)}\prod_{i\ge 1}
\left[\begin{array}{c}\ld_i'-\mu_{i+1}'\\ \ld_i'-\mu_i'\end{array}
\right]_q,~~~c_(\ld ;\mu)=\sum_{i\ge 0}\left(\begin{array}{c}
\ld_i'-\mu_i'\\ 2\end{array}\right), \label{2.18}
\end{equation}
and therefore
$$\varphi_{(n)\mu}^{\ld}(q)=0$$
unless $\mu\subset\ld$, $|\ld\setminus\mu|=n$, \cite{Kir5,HKKOTY}.

Perhaps, it is interesting to mention that if $p$ is a prime number,
then
$$\varphi_{(n)\mu}^{\ld}(p)=p^{n(\ld)-n(\mu)}\alpha_{\ld}(\mu ;p^{-1}),
$$
where for any partitions $\mu\subseteq\ld$ we denote by $\al_{\ld}(\mu
;p)$ the number of subgroups of type $\mu$ in a finite abelian
$p$--group of type $\ld$, see \cite{Bu} and \cite{Kir5} for definitions
and further details.

{\bf a.} Based on (\ref{2.18}), prove the following
combinatorial/fermionic formula for the $q$--weight multiplicity
${\cal P}_{\ld\mu}(q):={\rm dim}_qV_{(\mu)}(\ld)$:
\begin{equation}
{\cal P}_{\ld\mu}(q)=\sum_{\eta}K_{\eta\mu}K_{\eta\ld}(q)=\sum_{\{\nu\}}
q^{c(\nu )}\prod_{k=1}^{r-1}\prod_{i\ge 1}\left[\begin{array}{c}
\nu_i^{(k+1)}-\nu_{i+1}^{(k)}\\ \nu_i^{(k)}-\nu_{i+1}^{(k)}\end{array}
\right]_q, \label{2.19}
\end{equation}
summed over all flags of partitions $\{\nu\} =\{0=\nu^{(0)}\subset
\nu^{(1)}\subset\cdots\subset\nu^{(r)}=\ld'\}$, such that
$|\nu^{(k)}|=\mu_1+\cdots +\mu_k$, for all $k$, $1\le k\le r=l(\mu )$;
$$c(\nu )=\sum_{k=0}^{r-1}\sum_{i\ge 1}\left(\begin{array}{c}
\nu_i^{(k+1)}-\nu_i^{(k)}\\ 2\end{array}\right),
$$
where $\ld'$ denotes the conjugate of the partition $\ld$, see
e.g. \cite{Ma}, Chapter~I, p.2, and for any real number $x$ we put
$$\left(\begin{array}{c}x\\ 2\end{array}\right) :=x(x-1)/2.$$
For a proof and applications of (\ref{2.19}) see
\cite{Kir5,HKKOTY}.

Note also a preprint by G.~Lusztig, {\it Fermionic form and Betti
numbers}, math.QA/0005010, where a conjectural formula for the
Betti numbers of a lagrangian quiver variety (see H.~Nakajima,
Duke Math. J. {\bf 91} (1998), 515--560, for definitions), which
generalizes the first equality in (\ref{2.19}), is suggested. It
seems a very challenging task to find for lagrangian quiver
variety (constructed by H.~Nakajima [{\it ibid}] and G.~Lusztig
[{\it ibid}]) an analog of the second equality in (\ref{2.19}).

{\bf b.} Consider the transition matrix
$$M(e,P):=(M(e,P)_{\ld\mu})$$
which expresses the elementary symmetric polynomials $e_{\ld}$
(see e.g. \cite{Ma}, Chapter~I, \S 2) in terms of
Hall--Littlewood polynomials:
$$e_{\ld}(X_n)=\sum_{\mu}M(e,P)_{\ld\mu}P_{\mu}(X_n;q).
$$
It is well--known ([{\it ibid}], Chapter~III, \S 6) that
$$M(e, P)_{\ld\mu}=\sum_{\eta}K_{\eta'\ld}K_{\eta\mu}(q).
$$

$\bullet$ Use (\ref{2.17}) to derive the following combinatorial
formula for the transition coefficients between the elementary
symmetric and Hall--Littlewood polynomials:
\begin{equation}
M(e,P)_{\ld\mu}=\sum_{\{\nu\}}\prod_{k=1}^{n-1}\prod_{i\ge 1}
\left[\begin{array}{c}\nu_i^{(k+1)}-\nu_{i+1}^{(k+1)}\\
\nu_i^{(k)}-\nu_{i+1}^{(k+1)}\end{array}\right]_q, \label{2.26*}
\end{equation}
summed over all flags of partitions $\{\nu\}: =\{
0=\nu^{(0)}\subset\nu^{(1)}\subset\cdots\subset\nu^{(r)}=\ld'\}$
such that $\nu^{(k)}/\nu^{(k-1)}$ is a horizontal strip of length
$\mu_k$, $1\le k\le r=l(\mu )$.

For a proof and applications of (\ref{2.26*}) see
\cite{Kir5,HKKOTY}.

\begin{rem} {\rm The specialization $q=1$ in the both sides of formula
(\ref{2.26*}) gives rise to a combinatorial formula for the
number of (0,1)--matrices with the prescribed row sums $\ld_i$
and column sums $\mu_j$. The latter set (denoted by ${\cal
R}_{\ld\mu}$) has drawn considerable attention among the
specialists in combinatorics and geometry, see e.g., survey
articles by R.A.~Brualdi (Linear Algebra Appl. {\bf 33} (1980),
159-231), and by G.~Ziegler (in {\it Polytopes--Combinatorics and
Computations}, G.~Kalai and G.~Ziegler (eds.), Birkh\"auser, 2000,
1-43). A problem of finding the precise number of matrices in the
set ${\cal R}_{\ld\mu}$ has been discussed frequently in the
literature, see e.g., R.J.~Ryser (Bull. Amer. Math. Soc. {\bf 66}
(1960), 442-464).

As far as I am aware, for the first time such a formula was
obtained by B.Y.~Wang (Scienta Sinica (Series~A) {\bf 1} (1988),
1-6) and simplified later by B.Y.~Wang and F.~Zhang (Discr. Math.
{\bf 187} (1998), 211-220). The formula for cardinality of the
set ${\cal R}_{\ld\mu}$ obtained by B.Y.~Wang and F.~Zhang,
involves $2^{n-1}-n$ summations and (probably) does not have a
natural $q$--analog which is compatible with the dual Knuth
correspondence (see e.g., D.~Knuth, Pacific J. Math. {\bf 34}
(1970), 709-727) and the Lascoux--Sch\"utzenberger statistics
{\it charge}, see e.g. \cite{Ma}, Chapter~III, \S 6. It is unclear
for the author are there some relations between the Wang--Zhang
formula for the number $|{\cal R}_{\ld\mu}|$ and that
(\ref{2.26*}). Finally, note that yet {\it another} combinatorial
formula for the cardinality of the set ${\cal R}_{\ld\mu}$
follows from Theorem~\ref{t7.9}.}
\end{rem}

{\bf 4.} (Generalized Rogers--Szeg\"o and modified
Hall--Littlewood polynomials)

{\bf a.} The homogeneous Rogers--Szeg\"o polynomial $H_N(x,y;q)$
is defined by
$$H_N(x,y;q)=\sum_{k=0}^N\left[\begin{array}{c}N\\
k\end{array}\right]_qx^ky^{N-k}.
$$
This polynomial can be viewed as a $q$--deformation of
homogeneous Hermite polynomial $H_N(x,y)$. Let us recall the
generating function for the Hermite polynomials
$$\sum_{N\ge 0}\frac{H_N(x,y)}{N!}t^N=\exp (2xt-y^2t^2),
$$
see e.g., G.~Andrews, R.~Askey and R.~Roy, {\it Special
Functions}, Cambridge Univ. Press, New York, 1999, Chapter~6.

$\bullet$ Show that generating function for the Rogers--Szeg\"o
polynomials $H_N(x,y;q)$ is given by
\begin{equation}
\sum_{N\ge 0}\frac{H_N(x,y;q)}{(q;q)_N}t^N=
\frac{1}{(tx;q)_{\infty}(ty;q)_{\infty}}. \label{2.26}
\end{equation}

$\bullet$ Use (\ref{2.26}) to derive the following recurrence
relation for Rogers--Szeg\"o polynomials
$$H_N(x,y;q)=(x+y)H_{N-1}(x,y;q)-(1-q^{N-1})xyH_{N-2}(x,y;q).
$$

{\bf b.} More generally, follow \cite{A}, Chapter~3, Example~17,
let us define the generalized Rogers--Szeg\"o polynomial
\begin{equation}
H_N(X_n;q)=\sum_{k_1+\cdots +k_n=N}\left[\begin{array}{c}N\\
k_1,\ldots ,k_n\end{array}\right]_qx_1^{k_1}\cdots x_n^{k_n}.
\label{2.27}
\end{equation}

$\bullet$ Use (\ref{1.3}) and (\ref{2.13}) to derive that
$$q^{\left(\begin{array}{c}N\\ 2\end{array}\right)}H_N(X_n;q^{-1})
=Q'_{(1^N)}(X_n;q).
$$
Recall that $X_n:=(x_1,\ldots ,x_n)$ stands for the set of
variables, and $Q'_{\mu}(X_n;q)$ denotes the modified
Hall--Littlewood polynomial corresponding to partition $\mu$.

$\bullet$ Show that the generating function for generalized
Rogers--Szeg\"o polynomials is given by
\begin{equation}
\sum_{N\ge 0}\frac{H_N(X_n;q)}{(q;q)_N}t^N
=\prod_{j=1}^n\frac{1}{(tx_j;q)_{\infty}}. \label{2.28}
\end{equation}
Based on (\ref{2.28}), deduce the following recurrence relation
for generalized Rogers--Szeg\"o polynomials
$$H_N(X_n;q)=\sum_{k=1}^n(-1)^{k-1}\frac{(q;q)_{N-1}}{(q;q)_{N-k}}
e_k(X_n)H_{N-k}(X_n;q),
$$
where $e_k(X_n)$ stands for the elementary symmetric function of
degree $k$ in the variables $X_n$.

Perhaps, it is interesting to mention that the specialization
\begin{eqnarray*}
Z_N^{(n)}(q)&:=&q^{-N(N-n)/2n}Q'_{(1^N)}(x_1=1,\ldots ,x_n=1;q)\\
&=& q^{-N(N-n)/2n}\sum_{\ld\vdash N}K_{\ld ,(1^N)}(q)\dim
V_{\ld}^{\g l(N)}
\end{eqnarray*}
is equal to the partition function for the so--called $SU(n)$
Polychronakos--Frahm model, defined by the Hamiltonian
$${\cal H}_{PF}=\sum_{1\le j<k\le N}\frac{P_{jk}}{(t_j-t_k)^2},
$$
where $P_{jk}$ denotes the so--called {\it exchange} operator in
$su(n)$ spin space, and $\{ t_j\}$ are the roots of the $N$-th
order (inhomogeneous) Hermite polynomial
$$H_N(t)=y^{-N}H_N(x,y),~~t=x/y.
$$
For more details, see K.~Hikami, J. Phys. Soc. Japan {\bf
64} (1995), 1047-1050, and the literature quoted therein.

\section{ Parabolic Kostka polynomials}
\label{pkp}
\neweq

Let $X_n=(x_1,\cdots ,x_n)$ be the set of independent variables.
For any sequence of integers $\gamma=(\gamma_1,\ldots ,\gamma_n)$ put
$$x^{\gamma}=x_1^{\gamma_1}x_2^{\gamma_2}\cdots x_n^{\gamma_n}.$$
The symmetric group $S_n$ acts on polynomials and rational
functions in $X_n=(x_1,\ldots ,x_n)$ by permuting  variables.

Define operators
\begin{eqnarray}
J(f)&=&\sum _{w\in S_n}(-1)^{l(w)}w(x^{\delta}f), \label{3.1}\\
\pi_n(f)&=&J(1)^{-1}J(f), \label{3.2}
\end{eqnarray}
where
$$J(1)=\ds\prod_{1\le i<j\le n}(x_i-x_j)$$
is the Vandermond determinant, and $\delta
:=\delta_n=(\hbox{$n-1,n-2,$}\ldots ,1,0)$.

\begin{exs} {\rm $1^0$ For the dominant (i.e. weakly decreasing) integral
weight $\lambda =(\lambda_1\ge\lambda_2\ge\cdots\ge\lambda_n)$, the
character $\ch V_{\lambda}$ of the irreducible highest weight $\lambda$
$\g l(n)$--module $V_{\lambda}$ is given by the Laurent polynomial
$$s_{\lambda}(X_n)=\pi_n(x^{\lambda})=\ch V_{\lambda};
$$
when $\lambda$ is a partition (i.e. $\lambda_n\ge 0$),
$s_{\lambda}(X_n)$ is the Schur function, see e.g. \cite{Ma},
Chapter~I, \S 3.

$2^0$ Let $\lambda =(\lambda_1\ge\lambda_2\ge\cdots\ge\lambda_n\ge 0)$ be
a partition, then
$$\pi_n\left(x^{\lambda}\prod_{1\le i<j\le n}(1-qx_j/x_i)\right)
=v_{\lambda}(q)P_{\lambda}(X_n;q),
$$
where $P_{\lambda}(X_n;q)$ denotes the Hall--Littlewood polynomial,
see e.g. \cite{Ma}, Chapter~III, \S 1, and
$$(1-q)^nv_{\lambda}(q)=\prod_{i\ge 0}(q;q)_{\lambda_i'-\lambda_{i+1}'}.
$$
}
\end{exs}

Here $\lambda'=(\lambda_1',\cdots ,\lambda_m')$ is the conjugate
partition to $\lambda$; by definition we set $\lambda_0'=n$, and
$$(q;q)_m=\ds\prod_{j=1}^m(1-q^j).
$$

Let $\eta= (\eta_1,\ldots ,\eta_p)$ be a composition of $n$. Denote by
$\Phi(\eta)$ the set of ordered pairs $(i,j)\in\Z^2$ such that
$$1\le i \le\eta_1+\cdots +\eta_r<j\le n$$
for some $r$, $1\le r\le p$. For example, if $\eta =(1^n)$, then
$$\Phi(\eta)=\{(i,j)\in\Z^2|1\le i<j\le n\}.
$$

Let $\mu$ be a partition, $l(\mu)\le n$, and $\eta$ be a composition,
$|\eta|=n$. Consider the formal power series $B_{\eta}(X_n;q)$,
$H_{\mu\eta}(X_n;q)$ and $K_{\lambda\mu\eta}(q)$ defined by
\begin{eqnarray}
B_{\eta}(X_n;q)&=&\prod_{(i,j)\in\Phi(\eta)}(1-qx_i/x_j)^{-1}, \label{3.3}\\
H_{\mu\eta}(X_n;q)&=&\pi_n(x^{\mu}B_{\eta}(X_n;q)) \label{3.4}\\
&=&\sum_{\lambda}K_{\lambda\mu\eta}(q)s_{\lambda}(X_n), \label{3.5}
\end{eqnarray}
where the sum runs over the set of dominant weights $\lambda
=(\lambda_1\ge\cdots\ge\lambda_n)$ in $\Z^n$.

It is known \cite{SW} that the coefficients
$K_{\lambda\mu\eta}(q)$ are in fact polynomials in $q$ with
integer coefficients. If partitions $\ld$ and $\mu$, and a
composition $\eta$ such that $l(\mu )\le |\eta |$ are given, we
will call the polynomial $K_{\ld\mu\eta}(q)$ by {\it parabolic}
Kostka polynomial of type $(\ld\mu\eta )$, or simply by {\it
parabolic} Kostka polynomial if no confusion may happen. Note,
that by the same formulae (\ref{3.3})--(\ref{3.5}) the polynomials
$K_{\ld\mu\eta}(q)$ may be defined also in the case when $\mu$ is
a composition. However, in this general case polynomials
$K_{\ld\mu\eta}(q)$ may have {\it negative} coefficients, and
their combinatorial and representation theoretical meanings are
unclear.

\begin{con}\label{c3.2} {\rm (Kirillov--Shimozono \cite{KS})}
Let $\lambda$, $\mu$ be partitions, $|\lambda|=|\mu|$, $l(\mu)=n$,
and $\eta$ be a composition of $n$, then the coefficients
$K_{\lambda\mu\eta}(q)$ are polynomials in $q$ with non--negative integer
coefficients.
\end{con}

The partition $\mu$ and composition $\eta$, $l(\mu)\le |\eta|$,
$l(\eta)=p$, define a sequence of partitions $R=(\mu^{(1)},\ldots
,\mu^{(p)})$, $l(\mu^{(i)})=\eta_i$, by the following rule
$$\mu^{(i)}=(\mu_{\eta_1+\cdots +\eta_{i-1}+1},\cdots ,\mu_{\eta_1+\cdots
+\eta_i}), ~~~1\le i\le p,
$$
where by definition we put $\eta_0:=0$.

\begin{theorem}\label{t3.3} {\rm (\cite{SW})} Let $\lambda ,\mu ,\eta$ be as
above. Then
\begin{equation}
K_{\lambda\mu\eta}(1)={\rm
Mult}[V_{\lambda}:\otimes_{i=1}^pV_{\mu^{(i)}}], \label{3.6}
\end{equation}
i.e. $K_{\lambda\mu\eta}(1)$ is equal to the multiplicity of irreducible
highest weight $\lambda$ $\g l(n)$--module $V_{\lambda}$ in the tensor
product of irreducible highest  weight $\mu^{(i)}$
representations $V_{\mu^{(i)}}$, $1\le i\le p$, of the Lie algebra $\g l(n)$.
\end{theorem}

In the case when all partitions $\mu^{(i)}$ have rectangular shapes,
Theorem~\ref{t3.3} has been proved in \cite{Kir}.

\begin{ex} {\rm If $\eta=(\eta_1,\eta_2)\in\Z_{\ge 0}^2$,
$\eta_1+\eta_2=n$, and $\mu$ is a partition, $l(\mu)\le n$, then
\begin{equation}
K_{\lambda\mu\eta}(q)=q^{c(\lambda ,\mu ,\eta)}{\rm
Mult}[V_{\lambda}:V_{\mu^{(1)}}\otimes V_{\mu^{(2)}}], \label{3.7}
\end{equation}
for some constant $c(\lambda ,\mu, \eta)\in\Z_{\ge 0}$, where
$\mu^{(1)}=(\mu_1,\ldots ,\mu_{\eta_1})$,
$\mu^{(2)}=(\mu_{\eta_1+1},\ldots ,\mu_n)$.}
\end{ex}

\begin{rem} {\rm According to (\ref{3.6}) and Conjecture~\ref{c3.2}, the
parabolic Kostka polynomials $K_{\ld\mu\eta}(q)$ may be
considered as a $q$--analog of the tensor product multiplicities.
Another $q$--analog of the tensor product multiplicities had been
introduced by Lascoux, Leclerc and Thibon \cite{LLT}. Formula
(\ref{3.7}) shows that in general these two $q$--analogs are
different. However, it was conjectured in \cite{KS},
Conjecture~5, and \cite{Kir5}, Conjecture~6.5, that, in fact,
these two $q$--analogs coincide in the case when partition $\mu$
and composition $\eta$ correspond to a dominant sequence of
rectangular shape partitions.}
\end{rem}

In the sequel, we are mainly interested in the case when for a given
partition $\mu$ a composition $\eta$ is chosen in such a way that the
corresponding  sequence of partitions $R=(\mu^{(1)},\ldots ,\mu^{(p)})$
contains only rectangular shape partitions, i.e.
$\mu^{(a)}=(\mu_a^{\eta_a})$ for all $1\le a\le p$. Conversely, any
sequence $R=(R_1,\ldots , R_p)$ of rectangular shape partitions
$R_a=(\mu_a^{\eta_a})$, $1\le a\le p$, such that
$\mu_1\ge\cdots\ge\mu_p$ (we will call such sequence {\it dominant
sequence of  rectangular shape partitions}) defines the partition $\mu$
and composition $\eta$. Namely, the parts of composition
$\eta:=\eta(R)=(\eta_1,\ldots ,\eta_p)$ are equal to the lengths
$\eta_a:=l(R_a)$ of partitions $R_a$. The partition $\mu :=\mu(R)$ is
obtained by concatenating the parts of partitions $R_a$ in order.

\begin{ex} {\rm Take $\lambda=(5,4,4,2,1)$ and
$R=((3,3),(2,2,2),(2),(1,1))$. Then $\eta(R)=(2,3,1,2)$ and
$\mu(R)=(3,3,2,2,2,2,1,1)$. In this case
$$K_{\lambda\mu\eta}(q)=q^5(1+3q+4q^2+2q^3).
$$}
\end{ex}

\begin{de}\label{d3.6} The parabolic Kostka polynomial $K_{\lambda R}(q)$
corresponding to a partition $\ld$ and a dominant sequence of rectangular
shape partitions $R$ is defined by the following formula
\begin{equation}
K_{\ld R}(q)=K_{\ld\mu (R)\eta(R)}(q). \label{3.8}
\end{equation}
\end{de}

We will call the value of parabolic Kostka polynomial $K_{\ld
R}(q)$ at $q=1$ by {\it parabolic} Kostka number and denote it by
$$K_{\ld R}:=K_{\ld R}(1).$$

\begin{theorem} {\rm (\cite{KSS})} Assume that $\ld$ is a partition and $R$
is a dominant sequence of rectangular shape partitions. Then the
parabolic Kostka polynomial $K_{\ld R}(q)$ has non--negative integer
coefficients.
\end{theorem}

\begin{exs} {\rm Let $\ld$ be a partition.

$1^0$ Let $R_a$ be the single row $(\mu_a)$ for all $a$, where
$\mu=(\mu_1,\mu_2,\ldots)$ is a partition of length at most $n$. Then
\begin{equation}
K_{\ld R}(q)=K_{\ld\mu}(q), \label{3.9}
\end{equation}
i.e. $K_{\ld R}(q)$ coincides with the Kostka--Foulkes polynomial
$K_{\ld\mu}(q)$.

Proof of (\ref{3.9})  follows from the following, probably well--known,
identity
\begin{equation}
\pi_n\left(x^{\mu}\prod_{1\le i<j\le n}(1-qx_i/x_j)^{-1}\right)=
\sum_{k\ge 0}e_n(X_n)^{-k}Q_{\mu +(k^n)}'(X_n;q), \label{3.10}
\end{equation}
where $e_n(X_n)=x_1\cdots x_n$; $Q_{\nu}'(X_n;q)$ denotes the
modified Hall--Littlewood polynomial corresponding to partition
$\nu$, see e.g. \cite{Kir5}. Recall that
$$Q_{\nu}'(X_n;q):=Q_{\nu}\left[\frac{X_n}{1-q};q\right]
=\sum_{\mu}K_{\nu\mu}(q)s_{\mu}(X_n).
$$
See Remark~\ref{r3.9} for yet another proof of the equality (\ref{3.9}).

$2^0$ Let $R_a$ be the single column $(1^{\eta_a})$ for all $a$, and
$\eta =(\eta_1,\eta_2,\ldots )$. Then
\begin{equation}
K_{\ld R}(q)={\overline K}_{\ld'\eta^+}(q), \label{3.11}
\end{equation}
the cocharge Kostka--Foulkes polynomial, where $\ld'$ is the conjugate of
the partition $\ld$, and $\eta^+$ is the partition obtained by sorting
the parts of $\eta$ into weakly decreasing order. Formula (\ref{3.11})
follows from that (\ref{3.10}) and the duality theorem for parabolic
Kostka polynomials.
\vskip 0.3cm

\hskip -0.6cm{\bf Duality Theorem} (\cite{Kir6,KS}) {\it Let $\ld$
be a partition, and $R$ be a dominant sequence of rectangular
shape partitions, $R=((\mu_a^{\eta_a}))_{a=1}^p$. Denote by $R'$
a dominant rearrangement of the sequence of rectangular shape
partitions $((\eta_a^{\mu_a}))_{a=1}^p$ obtained by transposing
each of the rectangular in $R$. Then
\begin{equation}
K_{\ld'R'}(q)=q^{n(R)}K_{\ld R}(q^{-1}), \label{3.12}
\end{equation}
where $n(R)=\ds\sum_{1\le a<b\le p}\min (\mu_a,\mu_b)\min (\eta_a,\eta_b)$.}
\vskip 0.3cm

Note that the left hand side of (\ref{3.12}) is computed in $\g
l(m)$, where $m=\sum\mu_a$ is the total number of columns in the
rectangles of $R$, whereas the right hand side of (\ref{3.12}) is
computed in $\g l(n)$, where $n=\sum\eta_a$ is the total number of
rows in the rectangles of $R$.

$3^0$ Let $k$ be a positive integer and $R_a$ be the rectangular
with $k$ columns and $\eta_a$ rows, $|\eta|=n$. Then \cite{W,SW}
$K_{\ld R}(q)$ is the Poincare polynomial of the isotopic
component of the irreducible $GL(n)$--module of highest weight
$(\ld_1-k,\ld_2-k,\ldots ,\ld_n-k)$ in the coordinate ring of the
Zariski closure of the nilpotent conjugacy class whose Jordan
form has diagonal block sizes given by the transpose of the
partition $\eta^+$. In the case $\eta =(1^n)$, the polynomial
$$K_{\ld R}(q)=\ds\sum_{\{e_i\}}q^{e_i}$$
gives the Kostant generalized exponents $\{e_i\}$ of
representation $V_{\ld}$.}
\end{exs}

\vskip 0.2cm \hskip -0.6cm {\bf Exercises} \vskip 0.2cm

\hskip -0.6cm{\bf 1.} Let $\eta =(\eta_1,\ldots ,\eta_n)$ be a
composition, $|\eta|=n$. Show that
$$\pi_n\left(\prod_{(i,j)\in\Phi(\eta)}(1-qx_j/x_i)\right)=
\left[\begin{array}{c} n\\ \eta_1,\ldots ,\eta_n\end{array}\right]_q,
$$
where
$$\left[\begin{array}{c} n\\ \eta_1,\ldots
,\eta_n\end{array}\right]_q:=\ds\frac{(q;q)_n}{\prod_i(q;q)_{\eta_i}}
$$
denotes the $q$--multinomial coefficient. \\
{\bf 2.} Take partitions $\ld =(n+k,n,n-1,\ldots ,3,2)$ and $\mu
=\ld'=(n,n,n-1,\ldots ,2,1^k)$. Show that if $n\ge k$, then
$$K_{\ld\mu}(q)=q^{2k-1}\left(\left(\begin{array}{c}n\\
k-1\end{array}\right)+qR(q)\right)
$$
for some polynomial $R(q)\in \N [q]$. See also Section~\ref{kfs},
Exercise~3{\bf e} for a generalization of this result.
\\
{\bf 3.} i) Let $\ld ,\mu ,\nu$ be partitions,
$|\nu|=|\ld|+|\mu|$, $l(\ld)=p$, $l(\mu)\le s$. Consider
partition
$$\wt\ld =(\ld_1+\mu_1,\ldots ,\ld_1+\mu_s,\ld_1,\ld_2,\ldots ,\ld_p)$$
and a dominant rearrangement $\wt R$ of the sequence of
rectangular shape partitions $R=\{\nu\cup (\ld_1^s)\}$.

$\bullet$ Show that
$$K_{\wt\ld ,\wt R}(q)=q^{Q_{\ld_1}(\nu)-|\ld|}\{ c_{\ld\mu}^{\nu}+
\cdots +q^{n(\nu)-n(\ld)-n(\mu)}\},
$$
where $c_{\ld\mu}^{\nu}$
denotes the {\it Littlewood--Richardson number}, i.e.
$$c_{\ld\mu}^{\nu}={\rm Mult}[V_{\nu}:V_{\ld}\otimes V_{\mu}],
~~{\rm and}~~ Q_{\ld_1}(\nu)=\ds\sum_{j\le\ld_1}\nu_j'.
$$

More generally, let $\lambda\supset\mu$ be partitions such that
the complement $\lambda\setminus\mu$ is a disjoint union of
partitions $\ld^{(1)},\ldots, \ld^{(p)}$, and $l(\mu )=m$. Let
$\nu$ be a partition, define composition $\wt\nu =(\mu,\nu)$ and
partition $\eta=(m,1^{|\nu|})$.

$\bullet$ Show that
$$K_{\lambda\wt\nu\eta}(q)=q^{c(\lambda,\mu,\nu)}
(c^{\nu}_{\ld^{(1)},\ldots ,\ld^{(p)}}+\cdots
+q^{n(\nu)-n(\ld^{(1)})-\cdots -n(\ld^{(p)})}),
$$
where
$$c^{\nu}_{\ld^{(1)},\ldots ,\ld^{(p)}}:={\rm
Mult}[V_{\nu}:V_{\ld^{(1)}}\otimes\cdots\otimes V_{\ld^{(p)}}]
$$
denotes the (multiple) Littlewood--Richardson coefficient, and
$c(\lambda,\mu, \nu)\in\Z_{\ge 0}$.

It seems plausible that polynomial $K_{\lambda\wt\nu\eta}(q)$ has
{\it non--negative} integer coefficients (see
Conjecture~\ref{c3.10*}), and the degree $d(\lambda,\mu, \nu)$ of
parabolic Kostka polynomial $K_{\ld\wt\nu\eta}(q)$ is {\it
linear}, i.e. $d(n\lambda,n\mu,n\nu)=nd(\lambda,\mu, \nu)$ for
any po\-sitive integer $n$ (generalized {\it saturation
conjecture} for multiple Littlewood--Richardson numbers).

Note also, that the order of parts in the definition of
composition $\wt\nu$ is very essential.

See also Exercise~4 to Section~\ref{kfs}, where connections between
this Exercise and the so--called saturation conjecture are explained.

ii) (Ribbon tableaux and parabolic Kostka polynomials)

Let $\lambda\supset\mu$ be partitions such that the complement
$\lambda\setminus\mu$ is a disjoint union of partitions
$\ld^{(1)},\ldots ,\ld^{(p)}$, and $l(\mu)=m$. It is well--known
and goes back to D.E.~Littlewood (Proc. R. Soc. A {\bf 209}
(1951), 333--353) that there exists a unique partition $\Lambda$
such that

$\bullet$ $p$--core$(\Lambda)=\emptyset$;

$\bullet$ $p$--quotient$(\Lambda)=(\ld^{(1)},\ldots ,\ld^{(p)})$,
see e.g. \cite{Ma}, Chapter~I, Example~8.

It follows from the very definition that
$$|\Lambda|=p(|\ld^{(1)}|+\cdots +|\ld^{(p)}|).$$
Now let $\nu$ be a partition such that
$$|\nu|=|\ld^{(1)}|+\cdots +|\ld^{(p)}|,$$
and ${\rm Tab}^{(p)}(\Lambda,\nu)$ denotes the set of semistandard
$p$--{\it ribbon} tableaux of shape $\Lambda$ and weight $\nu$,
see e.g. \cite{LLT}, Section~4, or \cite{Ma}, Chapter~I,
Section~8, Example~8.

Finally, define composition $\wt\nu =(\mu,\nu)$ and partition
$\eta =(m,1^{|\nu|})$, and consider the set $LR(\ld,\wt\nu,\eta)$
of Littlewood--Richardson's tableaux of shape $\ld$ and type
$(\wt\nu,\eta)$, see e.g. \cite{KSS}, Section~2.

$\bullet$ Based on the Stanton--White method (D.~Stanton and
D.~White, J. Comb. Theory A {\bf 40} (1985), 211--247) construct
a bijection
$$\pi :{\rm Tab}^{(p)}(\Lambda,\nu)\to LR(\ld ,\wt\nu,\eta).
$$
In particular,
$$K_{\ld\wt\nu\eta}(1)=|{\rm
Tab}^{(p)}(\Lambda,\nu)|=\epsilon_p(\Lambda)K_{\Lambda,\underbrace{\nu\oplus
\cdots\oplus\nu}_p}(\zeta_p),
$$
where $\zeta_p=\exp(2\pi i/p)$ and $\epsilon_p(\Lambda)=\pm 1$.
The second equality is due to A.~Lascoux, B.~Leclerc and
J.-Y.~Thibon (C.R. Acad. Sci. Paris {\bf 316} (1993), 1--6).

Now follow \cite{LLT}, Section~5, denote by spin$(T)$ the {\it
spin} of a ribbon tableau $T$. Recall that the spin of a ribbon
tableau $T$ is by definition the sum of the spins of its ribbons;
the spin of a ribbon $R$ is by definition $(h(R)-1)/2$, where
$h(R)$ denotes the {\it height} of $R$; the height of a ribbon is
defined to be one less than the number of rows it occupies.

Finally, define the $(q,t)$--ribbon polynomial
$$c^{\nu}_{\ld^{(1)},\ldots ,\ld^{(p)}}(q,t)=\sum_{T\in{\rm
Tab}^{(p)}(\Lambda,\nu)}t^{{\rm spin~}(T)}q^{\wt c(\pi (T))},
$$
where
$$\wt c(\pi(T)):=c(\pi(T))-d(\ld,\mu,\nu).$$

From the very definition polynomials $c^{\nu}_{\ld^{(1)},\ldots
,\ld^{(p)}}(q,t)$ have the following properties

$\bullet$ $c^{\nu}_{\ld^{(1)},\ldots ,\ld^{(p)}}(q,t)\in\N[q,t]$;

$\bullet$ $c^{\nu}_{\ld^{(1)},\ldots ,\ld^{(p)}}(1,1)=|{\rm
Tab}^{(p)}(\Lambda,\nu)|$;
\begin{con}\label{c3.10*} Bijection $\pi$ may be constructed in such
a way that

$\bullet$ $c^{\nu}_{\ld^{(1)},\ldots
,\ld^{(p)}}(q,1)=K_{\ld\wt\nu\eta}(q)$;

$\bullet$ $c^{\nu}_{\ld^{(1)},\ldots
,\ld^{(p)}}(0,t)=c^{\nu}_{\ld^{(1)},\ldots ,\ld^{(p)}}(t), $\\
where $c^{\nu}_{\ld^{(1)},\ldots ,\ld^{(p)}}(t)$ denotes the
$t$--deformation of multiple Littlewood--Richardson numbers
introduced by A.~Lascoux, B.~Leclerc and J.-Y.~Thibon {\rm (see
e.g., B.~Leclerc and J.-Y.~Thibon, {\it Littlewood--Richardson
coefficients and Kazhdan--Lusztig polynomials},
math.QA/9809122,\break 52p.)}.
\end{con}
{\bf 4.} {\bf a.} Let $n\ge 1$ be an integer, consider staircase
partition
$$\delta_n=(n,{n-1,}\ldots, 2,1)$$
and partitions $\mu_n=(\delta_{n-1},1^n)$, and $\eta_n=(n-1,1^n)$.
Show that
$$K_{\delta_n,\mu_n,\eta_n}(q)=q^{n-1}\prod_{j=1}^n\frac{1-q^j}{1-q}.
$$

{\bf b.} Let $n\ge 2$ be an integer.

i) Consider staircase partition $\delta_n=(n,n-1,\ldots ,2,1)$ and
partitions $\mu_n=(\delta_{n-2},1^{2n-1})$, and
$\eta_n=(n-2,1^{2n-1})$. Show that
$$K_{\delta_n,\mu_n,\eta_n}(1)=\frac{2^{2n}(2^{2n}-1)}{2n}|B_{2n}|.
$$
Equivalently,
$$\sum_{n\ge
1}K_{\delta_n,\mu_n,\eta_n}(1)\frac{x^{2n-1}}{(2n-1)!}=\tan x.
$$
Hence, the number of standard Young tableaux of skew shape
$\delta_n\setminus\delta_{n-2}$ is equal to the tangent number
$$T_n:=\ds\frac{2^{2n}(2^{2n}-1)}{2n}|B_{2n}|.$$
Here and further $B_n$ denotes the $n$-th Bernoulli number.
Recall, that
$$\ds\sum_{n\ge 0}B_n\frac{x^n}{n!}=\frac{x}{\exp{x}-1}.$$

ii) Show that the Kostka--Foulkes polynomial
$K_{\delta_n,\mu_n,\eta_n}(q)$ is divisible by the product
$\ds\prod_{j=1}^{n-1}(1+q^j)$, and the ratio is a polynomial with
non--negative integer coefficients.

iii) Show that the number of semistandard Young tableaux of skew
shape $\delta_n\setminus\delta_{n-2}$ and weight $(2^{n-1},1)$ is
equal to $2^{-(n-1)}T_n$.

iv) Show that the number of semistandard Young tableaux of weight
$(2^{n-1})$ and skew shape $\delta_n\setminus (n-1,n-3,\ldots
,2,1)$ is equal to $2^{n-2}|G_{2n-2}|$, where $G_{2n}$ denotes the
$n$--th Genocci number. Recall that
$$\sum_{n=1}^{\infty}G_n\frac{x^n}{n!}=\frac{2x}{1+\exp{x}}
=x(1-\tanh\frac{x}{2}).
$$
It is easy to see that
$$|G_{2n}|=2(2^{2n}-1)|B_{2n}|, ~~
G_{2n+1}=0,~~{\rm if}~~ n>0.$$

v) Consider partitions $\al_n=(n,n-1,\ldots ,3,2)$,
$\beta_n=(\delta_{n-2},1^{2n-2})$, and $\gamma_n=(n-2,1^{2n-2})$.
Show that
$$K_{\al_n,\beta_n,\gamma_n}(1)=|E_{2n-2}|.
$$
Equivalently,
$$\sum_{n\ge
1}K_{\al_n,\beta_n,\gamma_n}(1)\frac{x^{2n-2}}{(2n-2)!}=\sec x.
$$
Hence, the number of standard Young tableaux of skew shape
$(n,{n-1,}\ldots ,3,2)\setminus\delta_{n-2}$ is equal to the
secant (or Euler) number $|E_{2n-2}|$.

Here $E_{2n}$ denotes the $n$--th Euler number. Recall that
$$\ds\sum_{n\ge 0}E_{2n}\frac{x^{2n}}{(2n)!}={\sec{x}}.
$$

Formulae for the numbers $K_{\delta_n,\mu_n,\eta_n}(1)$ from the
items $i)$ and $v)$ of Exercise~{\bf 4b} are well--known, see e.g.
H.O.~Foulkes (Discr. Math. {\bf 15} (1976), 235-252), or I.~Gessel
and G.~Viennot (Adv. Math. {\bf 58} (1985), 300-321).

{\bf c.} i) Denote by ${\cal A}_n$ the set of {\it alternating}
(or {\it zigzag}) permutations of $n$ letters, i.e.
$${\cal A}_n=\left\{(a_1,\ldots ,a_n)\in S_n\vert\begin{array}{l}
a_{2i-1}<a_{2i},~~i=1,\ldots ,\left[\frac{n}{2}\right],\\
a_{2i}>a_{2i+1},~~i=1,\ldots
,\left[\frac{n-1}{2}\right].\end{array}\right\}
$$

Show that
$$\ds\sum_{n\ge 0}|{\cal A}_n|\frac{x^n}{n!}=\tan{x}+\sec{x}.
$$

This result is due to D.~Andr\'e (J. Math. Pures Appl. {\bf 7}
(1881), 167-184).

ii) Denote by ${\cal S}_n$ the set of permutations $w=w_1w_2\cdots
w_n\in S_n$ satisfying the following conditions

$\bullet$ if $n=2k$, then $i$ and $i+1$ must precede $i+k$ for
$1\le i\le k-1$, while $k$ must precede $2k$;

$\bullet$ if $n=2k+1$, then $i$ and $i+1$ must precede $i+k+1$ for
$1\le i\le k$.

Show that
$$\ds\sum_{n\ge 1}|{\cal S}_n|\frac{x^n}{n!}=\tan{x}+\sec{x}.
$$

This result is due to R.~Stanley (Ann. Discr. Math. {\bf 6}
(1980), 333-342).

iii) Construct bijection between the set of alternating
permutations ${\cal A}_n$ and that of ${\cal S}_n$.

{\bf d.} Let ${\cal P}_n$ be a convex integral polytope in $\R^n$
determined by the following inequalities

$\bullet$ $x_i\ge 0$, $1\le i\le n$,

$\bullet$ $x_i+x_{i+1}\le 1$, $1\le i\le n-1$.

i) Show that the number of vertices of the polytope ${\cal P}_n$
is equal to the $(n+1)$--th Fibonacci number $F_{n+1}$.

Recall that Fibonacci's numbers $F_k$ may be defined using the
generating function
$$\sum_{k\ge 0}F_kt^k=\frac{1}{1-t-t^2}.
$$
{\bf Question.} What are the numbers $f_i({\cal P}_n)$, $0\le i\le
n-1$, of $i$--dimensional faces of the polytope ${\cal P}_n$?

ii) Let $\psi_n$ denotes the volume of the polytope ${\cal P}_n$.
Show that
$$\ds\sum_{n\ge 0}\psi_nx^n=\tan{x}+\sec{x}.
$$

This Exercise was proposed by R.~Stanley (Amer. Math. Monthly {\bf
85} (1978), p.197, Problem E2701). For the solution to this
problem given by I.G.~Macdonald and R.B.~Nelsen (independently),
see Amer. Math. Monthly {\bf 86} (1979), p.396.

iii) Let \vskip -0.5cm
$$\delta ({\cal P}_n)=(\delta_0^{(n)},\delta_1^{(n)},
\ldots ,\delta_{n-1}^{(n)})$$ denotes the $\delta$--vector of the
polytope ${\cal P}_n$, i.e.
$$\sum_{k\ge 0}i({\cal P}_n;k)t^k=\left(\sum_{k=0}^n\delta_k^{(n)}
t^k\right)/(1-t)^{n+1},
$$
where $i({\cal P}_n;k)$ denotes the number of points $\al\in{\cal
P}_n$ such that $k\al\in\Z^n$.

$\bullet$ Show that
\begin{eqnarray*}
&&\delta_n^{(n)}=\delta_{n-1}^{(n)}=0,~~
\delta_0^{(n)}=\delta_{n-2}^{(n)}=1,~~
\delta_1^{(n)}=\delta_{n-3}^{(n)}= F_{n+1}-n-1,\\
&&\delta ({\cal P}_n;t)=\ds\sum_{k=0}^n\delta_k^{(n)}t^k=
\sum_{w\in{\cal S}_n}t^{{\rm des}{(w)}}.
\end{eqnarray*}

$\bullet$ Show that $\delta ({\cal P}_n;t)$ is a symmetric
polynomial (in the variable $t$).

\begin{con}\label{c3.10} The $\delta ({\cal P}_n;t)$ is a unimodal
polynomial.
\end{con}

{\bf e.} i) Let $q=\exp\left(\ds\frac{2\pi i}{2k+3}\right)$ be
primitive root of unity. For each integer $m$ between 0 and $k$,
let $V_m$ denote the $(2m+1)$--dimensional irreducible
representation of the quantum universal enveloping algebra
$U_q(sl(2))$. Show that
$$i\left({\cal P}_n;k\right)={\rm
Mult}\left(V_{\left[\frac{k+1}{2}\right]}:
V_{\left[\frac{k+1}{2}\right]}^{\wh\otimes (n+2)}\right),
$$
i.e. the number of integer points $i({\cal P}_n;k)$ of the convex
integral polytope $k{\cal P}_n$ is equal to the multiplicity of
the irreducible representation $V_{\left[\frac{k+1}{2}\right]}$ of
the quantum universal enveloping algebra $U_q(sl(2))$ with
$$q=\ds\exp\left(\frac{2\pi i}{2k+3}\right),$$
in the $(n+2)$--th restricted tensor power of the representation
$V_{\left[\frac{k+1}{2}\right]}$. We refer the reader to
\cite{Kac}, Chapter~13, Exercise~13.34, for definition and basic
properties of the restricted tensor products.

ii) Denote by $A_k=(a_{ij})$ the $(k+1)$ by $(k+1)$ integer matrix
with the following matrix elements:
$$a_{ij}=\cases{1, & if $i+j\ge k+2$;\cr
0, & if $i+j<k+2$.}
$$

$\bullet$ Show that
$$i({\cal P}_n;k)=(A_k^{n+1})_{k+1,k+1},
$$
i.e. the number of integer points of the convex integral
polytope $k{\cal P}_n$ is equal to the $(k+1,k+1)$--entry of the
$(n+1)$--th power of the matrix $A_k$.

$\bullet$ Show that
$$A_k^2=(\min (i,j))_{1\le i,j\le k+1},
$$
and the characteristic polynomial
$$D_k(t):=\det(tI_{k+1}-A_k^2),
$$
where $I_{k+1}$ is the identity matrix of order $k+1$, satisfies
the recurrence relation
\begin{eqnarray*}
D_k(t)&=&(1-2t)D_{k-1}(t)-t^2D_{k-2}(t),~~{\rm if}~~k\ge 2,\\
D_0(t)&=&1-t,~~ D_1(t)=1-3t+t^2.
\end{eqnarray*}

$\bullet$ Show that the roots $t_a$ of the characteristic
polynomial $D_k(t)$ are given by
$$t_a=\frac{1}{4}\sec^2\frac{a\pi}{2k+3},~~1\le a\le k+1.
$$

$\bullet$ Show that
$$A_k^{-2}=(b_{ij})_{1\le i,j\le k+1},
$$
where
$$b_{ij}=\cases{2-\delta_{i,k+1}, & if $i=j$;\cr
-1, & if $|i-j|=1$;\cr 0, & if $|i-j|\ge 2$.}
$$

All statements of Exercise~{\bf e}, ii) may be considered as
well--known, but we are not able to find precise references where
these statements were formulated for the first time.

iii) Let $Z_n$ denotes the $n$--element "zigzag poset" with
elements $\{x_1,\ldots ,x_n\}$ and cover relations
$x_{2i-1}<x_{2i}$ and $x_{2i}>x_{2i+1}$.

$\bullet$ Show that
$$i({\cal P}_{n+1};k)=\Omega (Z_n;k),
$$
where $\Omega (Z_n;k)$ denotes the order polynomial of the poset
$Z_n$. We refer the reader to a book by R.~Stanley \cite{St3},
vol.1, Section~3.11, p.130, for a definition of the {\it order}
polynomial of a finite poset.

iv) Follow \cite{St3}, vol.1, Section~3, p.157, consider the
generating function \vskip -0.5cm
$$G_k(x)=\ds\sum_{n\ge 0}i({\cal P}_n;k)x^n.$$
Show that
$$G_1(x)=1/(1-x)~~{\rm and}~~ G_{k+1}=\frac{1+G_k(x)}{3-x^2-G_k(x)},
~~{\rm if}~~k\ge 1.
$$
This result is due to G.~Ziegler, see e.g. \cite{St3}, vol.1,
Chapter~3, Exercise~23~{\bf d}. The generating function $G_k(x)$
was computed for the first time by R.~Stanley (Ann. Discr. Math.
{\bf 6} (1980), 333-342, Example~3.2).

\begin{prb} Compute the two variables generating function
$$G(x,t)=\sum_{n,k\ge 0}i({\cal P}_n;k)x^nt^k.
$$
\end{prb}

v) Show that the number of integral points $i({\cal P}_n;k)$ of
the polytope $k{\cal P}_n$ is equal to the number of $n$--step
paths when light ray is reflected from $(k+1)$ glass plates (the
so--called {\it multiple reflections} or {\it wave sequences}).

Further details and references to the Exercise {\bf e}, v), may
be found in

$\bullet$ L.~Moser and M.~Wyman, {\it The Fibonacci Quarterly}
{\bf 11} (1973), 302-306;

$\bullet$ B.~Junge and V.~Hoggatt Jr., {\it The Fibonacci
Quarterly} {\bf 11} (1973), 285-291;

$\bullet$ G.~Kreweras, {\it Math. Sci. Hum.} {\bf 14} (1976),
n$^0$53, 5-30;

$\bullet$ J.~Berman and P.~K\"ohler, {\it Mitt. Math. Sem.
Giessen} {\bf 121} (1976), 103-124;\\ and the literature quoted
therein.

\begin{rem} {\rm (Fermionic formula for the numbers $i({\cal
P}_n;k)$) Consider partitions
\begin{eqnarray*}
\ld &:=&\ld_{(n,k)}=\ds\left((n+3)
\left[\frac{k+1}{2}\right],(n+1)\left[\frac{k+1}{2}\right]\right),\\
\mu &:=&\mu_{(n,k)}=\left(\underbrace{2\left[\frac{k+1}{2}\right],
\ldots ,2\left[\frac{k+1}{2}\right]}_{n+2}\right),
\end{eqnarray*}
and set $l=2k+1$. It follows from part i) that $i({\cal P}_n;k)$
is equal to the $l$--restricted Kostka number $K_{\ld\mu}^{(l)}$.

For the reader's convenience, let us recall \cite{Kir4} a
fermionic formula for the $l$--restricted Kostka polynomial
$K_{\ld\mu}^{(l)}(q)$ in the case of two row partition $\ld
=(\ld_1\ge\ld_2\ge 0)$ and the same size partition $\mu$:
\begin{equation}
K_{\ld\mu}^{(l)}(q)=\sum_{\nu}q^{2n(\nu)}\prod_{j\ge 1}
\left[\begin{array}{c}P_{j,l}(\nu;\mu)+m_j(\nu)\\ m_j(\nu)
\end{array}\right]_q, \label{3.13*}
\end{equation}
summed over all partitions $\nu=(\nu_1\ge\nu_2\ge\cdots\ge 0)$
such that

$\bullet$ $|\nu|=\nu_1+\nu_2+\cdots =\ld_2$, $\nu_1\le l$,

$\bullet$ $P_{j,l}(\nu;\mu):=\sum_a\min(j,\mu_a)
-\max(j+\ld_1-\ld_2-l,0)-2Q_j(\nu)\ge 0$ for all $j\ge 1$.

Here and further we are using the following notation: $Q_j(\nu
):=\sum_{k\le j}\nu_k'$, $m_j(\nu)=\nu_j'-\nu_{j+1}'$ and
$n(\nu)=\sum_{i\ge 1}(i-1)\nu_i$.

Thus, we come to a fermionic formula,
\begin{equation}
i({\cal P}_n;k)=K_{\ld_{(n,k)},\mu_{(n,k)}}(1)=\sum_{\nu}
\prod_{j\ge 1}\left(\begin{array}{c}P_j(\nu)+m_j(\nu)\\
m_j(\nu)\end{array} \right), \label{3.14*}
\end{equation}
summed over all partitions $\nu =(\nu_1\ge\nu_2\ge\cdots\ge 0)$
such that

$\bullet$ $|\nu|=\nu_1+\nu_2+\cdots
=(n+1)\ds\left[\frac{k+1}{2}\right]$, ~~$\nu_1\le 2k+1$,

$\bullet$ $P_j(\nu):=(n+2)\min(j,2\left[\frac{k+1}{2}\right])-
\max(j-2\left[\frac{k}{2}\right]-1,0)-2Q_j(\nu)\ge 0$ for all
$j\ge 1$.

In particular, if $k=1$ the formula (\ref{3.14*}) gives a
"fermionic expression" for Fibonacci's numbers. }
\end{rem}

vi) Show that
\begin{equation}
i({\cal P}_n;k)=\sum_{\nu}\prod_{j\ge
1}\left(\begin{array}{c}P_j(\nu;n)+m_j(\nu)\\
m_j(\nu)\end{array}\right), \label{3.15**}
\end{equation}
summed over all partitions $\nu =(\nu_1,\nu_2,\ldots )$ such that

$\bullet$ $\nu_1\le k$;

$\bullet$ $P_j(\nu;n):=(n+1)j-2Q_j(\nu)\ge 0$ for all $j\ge 1$,

where as before $m_j(\nu)=\nu_j'-\nu_{j+1}'$, and
$Q_j(\nu)=\sum_a\min (j,\nu_a)$.

More generally, denote by ${\cal P}_n^{(k)}(r,s)$ the convex
polytope
$$\{x=(x_1,\ldots ,x_n)\in\R_{\ge 0}^n|x\in k{\cal P}_n,~
x_1\le r,~x_n\le s\}.
$$

$\bullet$ Show that
\begin{equation}
\#|{\cal P}_n^{(k)}(r,s)\cap\Z^n|=\sum_{\nu}\prod_{j\ge 1}
\left(\begin{array}{c}P_j(\nu;n,r,s)+m_j(\nu)\\
m_j(\nu)\end{array}\right),\label{3.16**}
\end{equation}
summed over all partitions $\nu =(\nu_1,\nu_2,\ldots )$ such that

$\bullet$ $\nu_1\le k$;

$\bullet$
$P_j(\nu;n,r,s):=(n-1)j+\min(j,r)+\min(j,s)-2Q_j(\nu)\ge 0$ for
all $j\ge 1$.

Note that the RHS(\ref{3.14*}) and the RHS(\ref{3.15**}) have a
different combinatorial structure, but count the same number
($=i({\cal P}_n;k)$). This observation gives rise to an
interesting combinatorial identity:
RHS(\ref{3.14*})=RHS(\ref{3.15**}). It will be interesting to
find a direct combinatorial/bijective proof of the latter.

vii) Consider matrix $M:=M_k(x)$ of size $k+1$ by $k+1$ with
entries
$$m_{ij}=\cases{0, & if $i+j<k+2$;\cr x^{k-j+1}, & if $i+j\ge k+2$.}
$$
For each integer $r$ between 0 and $k$, define the column vector
$$|r\rangle =e_{r+1}+\cdots +e_{k+1}=(\underbrace{\overbrace{0,\ldots
,0}^r,1,\ldots ,1}_{k+1})^t\in{\rm Mat}_{1\times (k+1)}(\Z).
$$
Consider the following product of matrices
$$M(q^{N-1}x)\ldots M(qx)M(x)|r\rangle :=(a_{r,k,N}^{(0)}(x;q),
\ldots ,a_{r,k,N}^{(k)}(x;q))^t.
$$

$\bullet$ Show that for each integer $s$ between 0 and $k$, one
has
\begin{equation}
a_{r,k,N}^{(s)}(x;q)=\sum_{\nu}(qx)^{|\nu|}q^{2n(\nu)-Q_r(\nu)}
\prod_{j=1}^k\left[\begin{array}{c}P_j(\nu;r,s)+m_j(\nu)\\
m_j(\nu)\end{array}\right]_q, \label{3.17**}
\end{equation}
summed over all partitions $\nu =(\nu_1,\nu_2,\ldots )$ such that

$\bullet$ $\nu_1\le k$;

$\bullet$ $P_j(\nu;r,s):=(N-1)j+\min(j,r)+\min(j,s)-2Q_j(\nu)\ge
0$ for all $j\ge 1$;

$\bullet$ $n(\nu)=\sum (i-1)\nu_i$;
~~$m_j(\nu)=\nu_j'-\nu_{j+1}'$, ~~$Q_j(\nu)=\sum_a\min(j,\nu_a)$.

Show that
\begin{equation}
a_{r,k,N}^{(s)}(x;q)=\sum_{(b_1,\ldots ,b_N)\in{\cal
P}_N^{(k)}(r,s)}x^{\sum b_j}q^{\sum jb_j}. \label{3.18**}
\end{equation}
As far as I am aware, the product of matrices
$$M_k(q^{N-1}x)\ldots M_k(qx)M_k(x)|r\rangle$$
was invented by B.~Feigin (private communication, 1992) as a tool
for obtaining a fermionic formula for the level $k$ vacuum
representation of the affine Lie algebra $\widehat{sl}_2$, and
that for the characters of certain unitary representations of the
Virasoro algebra. This product of matrices was studied further by
E.~Frenkel and A.~Szenes (Int. Math. Res. Notices {\bf 2} (1993),
53-60) and by A.N.~Kirillov {\it Fusion algebra and Verlinde
formula}, Preprint of Isaac Newton Institute IN--92019, and
hep--th/9212084. To our knowledge, the fermionic formula
(\ref{3.17**}) for polynomials $a_{r,k,N}^{(s)}(x;q)$ was
obtained for the first time in \cite{Kir4}, p.99.

\vskip 0.2cm

{\bf f.} Consider the sets
$${\cal S}_n^{(i,j)}=\{w=w_1\cdots w_n\in{\cal S}_n|w_1=i,
w_n=j\},~~{\rm and}~~ {\cal S}_n^{(i)}=\bigcup_{j\in [1,n]}{\cal
S}_n^{(i,j)}.
$$

Show that
\vskip 0.3cm

i) $\ds\sum_{w\in{\cal S}_{n+1}^{(1)}}t^{{\rm des}(w)}=\delta
({\cal P}_n;t)=\sum_{u\in{\cal S}_n}t^{{\rm des}(u)}$;

\vskip 0.2cm

ii) if $n\equiv 0({\rm mod}2)$, then
$$\sum_{w\in S_n^{(i,j)}}t^{{\rm des}(w)}=\sum_{u\in
S_n^{(n-j+1,n-i+1)}}t^{{\rm des}(u)},~~{\rm and}~~ S_n^{(n)}=\{
u=wn|w\in S_{n-1}\};
$$

iii) if $n=2k+1$, $k\ge 1$, then
$$\sum_{w\in S_n^{(i,j)}}t^{{\rm des}(w)}=\sum_{u\in
S_n^{(k+2-i,3k+3-j)}}t^{n-2-{\rm des}(u)}.
$$
Let us denote by ${\cal S}_n^*$ the set of all permutations
$w=w_1\cdots w_n$ such that the permutation $\wt w=1(w_1+1)\cdots
(w_n+1)$ belongs to the set ${\cal S}_{n+1}$. It follows from part
i) that $|{\cal S}_n^*|=|{\cal S}_n|$.
\begin{ex}{\rm Consider the case $n=5$. Then
\begin{eqnarray*}
{\cal A}_5&=&\{13254,14253,14352,15243,15342,23154,24153,24351,
25143, 25341,\\&&~~34152,34251,35142,35241,45132,45231\};\\ {\cal
S}_5&=&\{12345,12354,12435,13245,13254,21345,21435,21354,23145,
23154,\\&&~~23514,31245,31254,32514,32154,32145\};\\ {\cal
S}_5^*&=&\{12345,12435,12354,12453,13245,12543,12534,13254,21354,
21345,\\&&~~21435,21453,21543,21534,25134,25143\}.
\end{eqnarray*}
Hence, $\delta ({\cal P}_5;t)=1+7t+7t^2+t^3$;
$$\sum_{w\in{\cal S}_5}t^{{\rm des}(w)}q^{{\rm maj}w}=
1+t(2q+2q^2+2q^3+q^4)+t^2(q^3+2q^4+2q^5+2q^6)+t^3q^7;
$$
$\psi_5:={\rm vol}({\cal P}_5)=\ds\frac{2^6(2^6-1)}{6!}B_6
=\frac{2}{15}$.

For the reader's convenience, we list below (for $n=5$) some
classes of permutations $w\in S_n$ which are equinumerous to the
set of alternating permutations ${\cal A}_n$:

$\bullet$ augmented Andr\'e permutations of the first kind (see
e.g., M.~Purtill, {\it Trans. AMS} {\bf 338} (1993), 77-104, and
the literature quoted therein)
\begin{eqnarray*}
\wt{\cal
A}_5^I&=&\{12345,12435,13425,23415,13245,14235,34125,24135,23145,
21345,\\&&~~41235,31245,21435,32415,41325,31425\};
\end{eqnarray*}

$\bullet$ Andr\'e permutations of the second kind (D.~Foata and
M.-P. Sch\"ut\-zenberger, in {\it A survey of Combinatorial
Theory}, J.N.~Srivastava et al. eds., Amsterdam, 1973, 173-187)
\begin{eqnarray*}
{\cal A}_5^{II}&=&\{12345,12534,14523,34512,15234,14235,34125,
45123,35124,51234,\\&&~~41235,31245,51423,53412,41523,31524\};
\end{eqnarray*}

$\bullet$ augmented simsun permutations (see e.g., R.~Stanley,
{\it Math. Z.} {\bf 216} (1994), p.498)
\begin{eqnarray*}
{\wt{\cal SS}}_5&=&\{12345,12435,13245,13425,14235,21345,21435,
23145,23415,24135,\\&&~~31425,31245,34125,41235,41325,42315\};
\end{eqnarray*}

$\bullet$ Jacobi permutations (G.~Viennot, {\it Journ. Comb. Th.
A} {\bf 29} (1980), p.124, Definition~3)
\begin{eqnarray*}
J_5&=&\{15432,35421,52431,42531,32541,54132,53142,13542,43152,
52143,\\&&~~15243,42153,14253,32154,13254,15432\}.
\end{eqnarray*}

It looks very challenging to construct {\it bijections} between
the sets ${\cal A}_n$, $\wt{\cal A}_n^I$, ${\cal A}_n^{II}$,
${\cal S}_n$, ${\cal S}_n^*$, $\wt{\cal SS}_n$ and $J_n$.}
\end{ex}

{\bf g.} Let $n\ge 2$, denotes by $\wt{\cal P}_n$ the convex
(rational) polytope in $\R^n$ determined by the following
inequalities

$\bullet$ $x_i\ge 0$, $1\le i\le n$,

$\bullet$ $x_i+x_{i+1}\le 1$, $1\le i\le n-1$,

$\bullet$ $x_1+x_n\le 1$.

Show that
$$\ds\sum_{n\ge 1}{\rm vol}(\wt{\cal
P}_n)x^{n-1}=\frac{1}{2}(\tan{x}+\sec{x}),$$
i.e.

$\bullet$ ${\rm vol}({\cal P}_{2k-1})=2{\rm vol}(\wt{\cal
P}_{2k})$, if $k\ge 1$,

$\bullet$ ${\rm vol}({\cal P}_{2k})=2{\rm vol}(\wt{\cal
P}_{2k+1})$, if $k\ge 0$,\\
where by definition we set $\wt{\cal P}_1:=[0,1/2]$.
\vskip 0.3cm

\hskip -0.6cm{\bf 5.} Let us fix $n\ge3$.
i) Consider partitions $\ld_n=(n+1,n,\ldots ,4,3)$,
$\mu_n=(n-2,n-3,\ldots ,2,1^{3n-3})$, and $\eta_n=(n-2,1^{3n-3})$.
Show that
\begin{equation}
K_{\ld_n\mu_n\eta_n}(1)=\frac{(3n-3)!(2^{2n-2}-2)}{(2n-2)!}|B_{2n-2}|,
\label{4.2}
\end{equation}
and thus, the number in the RHS(\ref{4.2}) is equal to the number
of standard (skew) Young tableaux of shape $(n+1,n,\ldots
,3)\setminus (n-2,n-3,\ldots ,1)$. This result is due to I.~Gessel
and G.~Viennot \cite{GV}, \S 11, and A.N.~Kirillov \cite{Kir7}.
It has been shown in \cite{Kir7} that
$$(3g)!\det\left\vert\frac{1}{(3-2i+2j)!}\right\vert_{1\le i,j\le
g}= \frac{(3g)!(2^{2g}-2)}{(2g)!}\vert B_{2g}\vert =(3g)!{\rm
vol}({\cal N}_{g+1}),
$$
where ${\rm vol}({\cal N}_g)$ stands for the volume of the moduli
space ${\cal N}_g$ of stable rank 2 vector bundles over a smooth
complex curve of genus $g$.


ii) Consider partitions $\ld_n=(n+1,n,\ldots ,4,3,2,1)$,\\
$\mu_n={(n-2,}{n-3,}\ldots ,2,1^{3n})$, and $\eta_n=(n-2,1^{3n})$.
Show that
\begin{equation}
K_{\ld_n\mu_n\eta_n}(1)=2^{2n}\frac{(3n)!}{(2n)!}|B_{2n}|.\label{4.3}
\end{equation}
{\bf 6.} Let $m,n\ge 2$ and $N\ge n$ be integers. Consider
partitions $\ld_N=N(n,{n-1,}\ldots ,2,1)+(m,\underbrace{0,\ldots
,0}_{n-1})$, $\mu_N=N(n,n-1,\ldots ,2,1)+(\underbrace{0,\ldots
,0}_n,m)$, and a composition $\eta =(1,n-2,1)$.

Show that if $m\ge n-1$, then
$$K_{\ld_N,\mu_N,\eta}(q)=\sum_{k=0}^m\left(\begin{array}{c}n-2+k\\
k\end{array}\right)q^{k+m}.
$$
{\bf 7.} Let $\al$ and $\beta$ be partitions, and $n,s$ be integer
numbers such that $s\ge |\beta|$. Define $m:=s-|\beta|\ge 0$, and
consider partition
$$\ld :=[\al ,\beta ]_s=(\al_1+s,\ldots ,\al_p+s,\underbrace{s,\ldots ,s}_{n-p-t},
s-\beta_t,\ldots ,s-\beta_1),
$$
where $p=l(\al)$, $t=l(\beta)$, and the dominant sequence of
partitions
$$R:=R_n^{r,s}=\{\underbrace{(1),\ldots
,(1)}_r,\underbrace{(1^{n-1}),\ldots ,(1^{n-1})}_s\},
$$
where
$r:=|\al|-|\beta|+s$. Let us introduce notation $S_m:=R_m^{m,m}$,
and assume additionally that $n\ge m+l(\al )+l(\beta)$.

i) Show that 
\begin{equation}
K_{\ld R}(q)\bdoteq\left[\begin{array}{c}r\\ m\end{array}\right]_q
\left[\begin{array}{c}s\\ m\end{array}\right]_qK_{\al
,(1^{|\al|})}(q) K_{\beta ,(1^{|\beta|})}(q)K_{(m^m),S_m}(q).
\label{3.15a}
\end{equation}

ii) Show that
$$K_{(m^m),S_m}(1)=m!,$$
and $K_{\ld R}(1)$ is equal to the number $c_{r,s}^{[\al ,\beta]}$
of up--down staircase tableaux of shape
$$\{\al ,\beta\}_n:= (\al_1,\ldots
,\al_{p-1},\al_p,\underbrace{0,\ldots ,0}_{n-p-t}, -\beta_t,\ldots
,-\beta_1)$$ and type $\epsilon_{r,s}:=(\underbrace{1,\ldots
,1}_r,\underbrace{-1,\ldots ,-1}_s)$, see \cite{Stm},
Definition~4.4. In particular, the number $a_{2s}^{[\al ,\beta]}$
of {\it alternating} (or {\it oscillating}) tableaux, i.e.
up--down staircase tableaux of shape $\{\al ,\beta\}_n$ and type
$\epsilon =(\underbrace{1,-1,1,-1,\ldots ,1,-1}_{2s})$, [{\it
ibid}], is equal to
$$(s-|\al|)!\left(\begin{array}{c}s\\|\al|\end{array}\right)^2
f^{\al}f^{\beta},~~{\rm if}~~ n\ge s-|\al|+l(\al)+l(\beta).
$$
Formulae for the numbers $c_{r,s}^{[\al,\beta]}$ and
$a_{2s}^{[\al,\beta]}$ are due to J.~Stembridge \cite{Stm},
Proposition~4.8, and in some particular cases to P.~Hanlon
\cite{Ha}. Formula (\ref{3.15a}) may be considered as a
$q$--analog of the above results obtained by J.~Stembridge.


\begin{exs}\label{ex3.10} {\rm $1^0$ Take $\al =\beta =(2,2)$ and $n=s=6$. Then
$\ld :=[\al ,\beta ]_6=(8,8,6,6,4,4)$, $m=2$, and
$R:=R_6^{6,6}=\{\underbrace{(1),\ldots ,(1)}_6,
\underbrace{(1^5),\ldots ,(1^5)}_6\}$. Using the fermionic formula
(\ref{4.1}), one can check that
$$K_{\ld R}(q)=q^{98}(1+q^2)^3\left[\begin{array}{c}6\\
2\end{array}\right]_q^2,~~~K_{\al ,(1^4)}(q)=q^2+q^4=
K_{(2^2),S_2}(q).
$$
$2^0$ Take $\al =\beta =(2,1)$, $n=7$, and $r=s=6$. Then $\ld
:=[\al ,\beta]_7=(8,7,6,6,6,5,4)$, $m=3$, and
$R:=R_7^{6,6}=\{\underbrace{(1),\ldots
,(1)}_6,\underbrace{(1^6),\ldots ,(1^6)}_6\}$. Using the fermionic
formula (\ref{4.1}), one can check that}
\begin{eqnarray*}
&&K_{\ld R}(q)=q^{95}(1+q)^2\left[\begin{array}{c}6\\
3\end{array}\right]_q^2(1+q^2+2q^3+q^4+q^6),\\
&&K_{(2,1),(1^3)}(q)=q+q^2,\\
&&K_{(333),S_3}(q)=q^9(1+q^2+2q^3+q^4+q^6).
\end{eqnarray*}
\end{exs}
\vskip 0.3cm

iii) Show that if $r\ge |\al|$, $s\ge |\beta|$, and
$r-s=|\al|-|\beta|$, then
\begin{equation}
K_{[\al ,\beta]_{n+1},R_{n+1}^{r,s}}(q)\ge K_{[\al
,\beta]_n,R_n^{r,s}}(q) \label{3.16*}
\end{equation}
for all $n\ge l(\al)+l(\beta)$, and if $n\ge
s-|\beta|+l(\al)+l(\beta)$, then the LHS(\ref{3.16*}) is equal to
the RHS(\ref{3.16*}).

\begin{ex}\label{ex3.11} {\rm Take $\al =\beta =(2,1)$, and
$r=s=6$. Using the fermionic formula (\ref{4.1}), one can check
that}
$$K_{[\al ,\beta]_7,R_7^{6,6}}(q)-K_{[\al
,\beta]_6,R_6^{6,6}}(q)=q^{98}(1+q+q^2+q^3+q^4)^2.
$$
\end{ex}

\hskip -0.6cm {\bf Questions.} i) Let $\ld$ be a partition and
$R=\{ R_a=(\mu_a^{\eta_a})\}^l_{a=1}$ be a dominant sequence of
rectangular shape partitions. Under what conditions on $\ld$ and
$R$ the parabolic Kostka number $K_{\ld R}\ne 0$ ?

One obvious necessary condition is the following: $\ld\ge\mu (R)$
with respect to the dominance ordering on the set of partitions.
Let us recall, see Section~\ref{pkp}, that $\mu (R)$ denotes the
partition which is obtained by concatenating the parts of
partitions $R_a$ in order.

ii) When does the parabolic Kostka number $K_{\ld ,R}$ equal to
1?

For example,
$$K_{(8n,8n,8n,8n,8n),((5n,5n),(5n,5n),(5n,5n),(5n,5n))}(q)=q^{24n}.$$

iii) Does there exist an explicit formula for computing the degree
of parabolic Kostka polynomial $K_{\ld R}(q)$ ? Recall that ${\rm
deg~}K_{\ld\mu}(q)=n(\mu)-n(\ld)$, see e.g. \cite{Ma}, Chapter~I.
It looks plausible that the degree $d(\ld,\mu;\eta)$ of parabolic
Kostka polynomial $K_{\ld\mu\eta}(q)$ is a homogeneous degree 1
function, i.e. $d(n\ld,n\mu;\eta)=nd(\ld,\mu ;\eta)$ for any
positive integer $n$.

\begin{ex} {\rm  Below we list some "exotic examples" of
Kostka--Foulkes and parabolic Kostka po\-ly\-no\-mials:
\begin{eqnarray*}
K_{(42),(2,1^4)}(q)&=&q^4(11211),\\
K_{(442),(3,3,2,1,1)}(q)&=&q^2(11311),\\
K_{(4422),(3,3,(2,2),1,1)}(q)&=&q^4(11411),\\
K_{(4431),(3,(2,2),(2,2),1)}(q)&=&2q^4(111),\\
K_{(444222),(3,(3,3),(2,2,2),2,1)}(q)&=&q^4(2332),\\
K_{(664),((3,3),(3,3),2,2)}(q)&=&q^8(21211),\\
K_{(654321),(4,(3,3,3),(2,2),(1,1,1,1))}(q)&=&3q^6(1+q)^3,\\
K_{(7654321),((4,4,4),(3,3,3,3),2,(1,1))}(q)&=&3q^7(1+q)^4,\\
K_{(7654321),((4,4,4,4),4,2,(1,1),(1,1,1,1))}(q)&=&q^8(5,17,24,17,5),\\
K_{(87654321),(5,(4,4,4,4),(3,3,3),(1,1,1,1,1),1)}(q)&=&6q^{11}(1+q)^4(2+3q+2q^2),\\
K_{(2n,2n,2n,2n),((n,n),(n,n),(n,n),(n,n))}(q)&=&q^{4n}\left[\begin{array}{c}n+2\\
2\end{array}\right]_{q^2}.
\end{eqnarray*}
More generally, let us take $\ld =(\underbrace{2n,\ldots
,2n}_{2k})$, and $R=((n^k),(n^k),(n^k),(n^k))$. Then
\begin{equation}
K_{\ld ,R}(q)=q^{2kn}\left[\begin{array}{c}n+k\\
k\end{array}\right]_{q^2}.\label{3.23}
\end{equation}
Indeed, there are exactly $\left(\begin{array}{c}n+k\\
k\end{array}\right)$ admissible configurations of type $(\ld ;R)$,
each of those has zero essential vacancy numbers only. Similarly,
let $n>m\ge 0$ and $k\ge 1$ be integers, consider partition $\ld
=(\underbrace{2n,\ldots ,2n}_{2k},\underbrace{2m,\ldots
,2m}_{2k})$ and sequence of rectangular shape partitions
$$R=(((m+n)^k),((m+n)^k),((m+n)^k),((m+n)^k)).
$$
Then
$$K_{\ld ,R}(q)=q^{2k(n-m)}\left[\begin{array}{c}k+n-m\\ k\end{array}\right]_{q^2}.
$$
$$K_{(9621),((3,3),(3,3),3,3)}(q)=q^{13}(25652)\bdoteq
s_{21}(1,q,q,q^2).
$$
See also Section~\ref{ipsf}, Exercise~14 for a generalization of
the latter example.}
\end{ex}

\section{Fermionic formula for parabolic Kostka po\-ly\-no\-mials}
\label{ffpkp}
\neweq

Let $\ld$ be a partition and $R=((\mu_a^{\eta_a}))_{a=1}^p$ be a sequence
of rectangular shape partitions such that
$$|\ld|=\sum_a|R_a|=\sum_a\mu_a\eta_a.$$

\begin{de}\label{d4.1} The configuration of type $(\ld ;R)$ is a sequence
of partitions $\nu =(\nu^{(1)},\nu^{(2)},\ldots )$ such that
$$|\nu^{(k)}|=\sum_{j>k}\ld_j-\sum_{a\ge 1}\mu_a\max (\eta_a-k,0)
=-\sum_{j\le k}\ld_j+\sum_{a\ge 1}\mu_a\min (k,\eta_a)
$$
for each $k\ge 1$.
\end{de}

Note that if $k\ge l(\ld)$ and $k\ge\eta_a$ for all
$a$, then $\nu^{(k)}$ is empty. In the sequel (except Corollary~\ref{c4.4})
we make the convention that $\nu^{(0)}$ is the empty partition.

For a partition $\mu$ define the number $Q_n(\mu)=\mu_1'+\cdots
+\mu_n'$, which is equal to the number of cells in the first $n$
columns of $\mu$. The {\it vacancy} numbers
$P_n^{(k)}(\nu
):=P_n^{(k)}(\nu ;R)$ of the configuration $\nu$ of type $(\ld
;R)$ are defined by
$$P_n^{(k)}(\nu )=Q_n(\nu^{(k-1)})-2Q_n(\nu^{(k)})+Q_n(\nu^{(k+1)})
+\sum_{a\ge 1}\min (\mu_a,n)\delta_{\eta_a,k}
$$
for $k,n\ge 1$, where $\delta_{i,j}$ is the Kronecker delta.

\begin{de}\label{d4.2} The configuration $\nu$ of type $(\ld ;R)$ is called
admissible, if
$$P_n^{(k)}(\nu ;R)\ge 0~~{\rm for~all}~~ k,n\ge 1.$$
\end{de}

We denote by $C(\ld ;R)$ the set of all admissible configurations
of type $(\ld ;R)$, and call the vacancy number $P_n^{(k)}(\nu
;R)$ {\it essential}, if $m_n(\nu^{(k)})>0$.

Finally, for configuration $\nu$ of type $(\ld ;R)$ let us define its
{\it charge}
$$c(\nu )=\sum_{k,n\ge 1}\left(\begin{array}{c}\alpha_n^{(k-1)}-
\alpha_n^{(k)}+\sum_a\theta (\eta_a-k)\theta (\mu_a-n)\\ 2\end{array}
\right),
$$
and {\it cocharge}
$$\overline c(\nu )=\sum_{k,n\ge 1}\left(\begin{array}{c}\alpha_n^{(k-1)}
-\alpha_n^{(k)}\\ 2\end{array}\right),
$$
where $\alpha_n^{(k)}=(\nu^{(k)})_n'$ denotes the size of the $n$--th
column of the $k$--th partition $\nu^{(k)}$ of the configuration $\nu$;
for any real number $x\in \R$ we put $\theta (x)=1$, if $x\ge 0$, and
$\theta (x)=0$, if $x<0$.

\begin{theorem}\label{t4.3} {\rm (Fermionic formula for parabolic Kostka
polynomials \cite{Kir6,KSS})}\\ Let $\ld$ be a partition and $R$
be a dominant sequence of rectangular shape partitions. Then
\begin{equation}
K_{\ld R}(q)=\sum_{\nu}q^{c(\nu )}\prod_{k,n\ge 1}\left[\begin{array}{c}
P_n^{(k)}(\nu ;R)+m_n(\nu^{(k)})\\ m_n(\nu^{(k)})\end{array}\right]_q,
\label{4.1}
\end{equation}
summed over all admissible configurations $\nu$ of type $(\ld ;R)$;
$m_n(\ld)$ denotes the number of parts of the partition $\ld$ of size
$n$.
\end{theorem}

\begin{cor}\label{c4.4} {\rm (Fermionic formula for Kostka--Foulkes
polynomials \cite{Kir1})} Let $\ld$ and $\mu$ be partitions of the same
size. Then
\begin{equation}
K_{\ld\mu}(q)=\sum_{\nu}q^{c(\nu)}\prod_{k,n\ge 1}\left[\begin{array}{c}
P_n^{(k)}(\nu ,\mu)+m_n(\nu^{(k)})\\ m_n(\nu^{(k)})\end{array}\right]_q,
\label{4.1a}
\end{equation}
summed over all sequences of partitions $\nu
=\{\nu^{(1)},\nu^{(2)},\ldots\}$ such that

$\bullet$ $|\nu^{(k)}|=\sum_{j>k}\ld_j$, $k=1,2,\ldots$;

$\bullet$  $P_n^{(k)}(\nu ,\mu):=Q_n(\nu^{(k-1)})-2Q_n(\nu^{(k)})
+Q_n(\nu^{(k+1)})\ge 0$ for all $k,n\ge 1$, where by definition we
put $\nu^{(0)}=\mu$;

\begin{equation}
\bullet~~c(\nu )=\sum_{k,n\ge
1}\left(\begin{array}{c}(\nu^{(k-1)})_n' -(\nu^{(k)})_n'\\
2\end{array}\right).~~~~~~~~~~~~~~~~~~~~~~~~~~~~~\label{4.3*}
\end{equation}
\end{cor}

\hskip -0.6cm {\bf Exercises} \vskip 0.2cm

\hskip -0.6cm{\bf 1.} Show that
\begin{equation}
\sum_{\nu}q^{c(\nu)}\prod_{k,n\ge 1}\left[\begin{array}{c}
P_n^{(k)}(\nu ;R)+m_n(\nu^{(k)})\\ m_n(\nu^{(k)})\end{array}\right]_q'
=0,\label{4.5}
\end{equation}
where the sum is taken over the set of all {\it non--admissible}
configurations of type $(\ld ;R)$, i.e. configurations $\nu$ of type
$(\ld ;R)$ such that
$$P_n^{(k)}(\nu ;R)+m_n(\nu^{(k)})<0$$
for some $k,n\ge 1$. Here, for any integer $N$ and non--negative
integer $M$ we use the notation $\left[\begin{array}{c}N\\
M\end{array}\right]_q'$ to denote the modified $q$--binomial
coefficient
$$\left[\begin{array}{c}N\\ M\end{array}\right]_q'=
\frac{(1-q^N)(1-q^{N-1})\cdots (1-q^{N-M+1})}{(1-q)(1-q^2)\cdots (1-q^M)}.
$$
{\bf 2.} (KOH's identity) Follow \cite{Z}, denote by $G(n,k)$ the
$q$--binomial coefficient $\left[\begin{array}{c}n+k\\
k\end{array}\right]_q$ (so that $G(n,k)=0$, if either $n<0$ or
$k<0$).

Show that
\begin{equation}
G(n,k)=\sum_{\nu\vdash k}q^{2n(\nu)}\prod_{j=1}^kG((n+2)j-2k+\sum_{l\ge
j}(l-j)d_l,d_j), \label{4.4}
\end{equation}
summed over all partitions $\nu =(1^{d_1},2^{d_2},\ldots
,k^{d_k})$ of $k$;
$$2n(\nu)=\sum_{j\ge 1}jd_j^2+2\sum_{i<j}\min (i,j)d_id_j-k.
$$
Binomial identity (\ref{4.4}) is due to Zeilberger \cite{Z} and
may be considered as an algebraic version of the main result,
Theorem~2.13, of K.~O'Hara's paper \cite{OH}. Note that if one
replaces in the identity (\ref{4.4}) the $q$--binomial
coefficient $G(n,k)$ by the
modified $q$--binomial coefficient $\left[\begin{array}{c}n+k\\
k\end{array}\right]_q'$, see definition in the previous
Exercise~1, then this new version of the identity (\ref{4.4})
becomes "much more elementary", and a short algebraic proof of it
was found by Macdonald \cite{Ma2}.

\begin{rem} {\rm Binomial identity (\ref{4.4}) is a particular
case of the identity (\ref{1.4}) in the case $N=k+n$ and $\ld
=(1^k)$.
Note that the symmetry relation
$$\left[\begin{array}{c}n+k\\ k\end{array}\right]=
\left[\begin{array}{c}n+k\\ n\end{array}\right]
$$
gives rise to a highly nontrivial identity between binomial sums.}
\end{rem}

\begin{prb} Let $\ld$ and $\mu$ be partitions, $l(\mu)\le n$, and $\eta$
be a composition of $n$. Find a fermionic formula for the
parabolic Kostka polynomial $K_{\ld\mu\eta}(q)$ which generalizes
that (\ref{4.1}).
\end{prb}
{\bf 3.} (Domino tableaux and configurations)\\
i) Let $\ld$ and $\mu$ be partitions such that $l(\mu)=s$, and
$C(\ld;\mu)$ be the corresponding set of admissible
configurations. Denote by $C_{(2)}(\ld;\mu)$ the subset of
$C(\ld;\mu)$ consisting of all configurations $\nu$ such that
\begin{equation}
P_n^{(k)}(\nu)m_n(\nu^{(k)})\equiv 0({\rm mod~} 2)~~{\rm for
~all}~~ n,k\ge 1. \label{4.6}
\end{equation}
Conditions (\ref{4.6}) may be rewritten in the following form
$$[P_n^{(k)}(\nu)/2]+[m_n(\nu^{(k)})/2]=[(P_n^{(k)}(\nu)+m_n(\nu^{(k)}))/2],
$$
where for any $x\in\R$, the symbol $[x]$  stands for the nearest
integer which is smaller or equal to $x$. Finally, let $c(\nu)$
denotes the {\it charge} of configuration $\nu$, see (\ref{4.3*})
for definition.

$\bullet$ Show that if $\nu\in C_{(2)}(\ld;\mu)$, then
$$c(\nu)\equiv\epsilon (\ld)({\rm mod~}2),
$$
where $\epsilon (\ld)=\ld_1'+\ld_3'+\ld_5'+\cdots$.

Next, for any configuration $\nu\in C_{(2)}(\ld;\mu)$ define its
{\it spin} by
\begin{equation}
{\rm spin}(\nu)=\frac{1}{2}\sum_{k,n\ge
1}\left\{\frac{(\nu^{(k-1)})_n'-(\nu^{(k)})_n'}{2}\right\},
\label{4.7}
\end{equation}
where for any $x\in\R$, the symbol $\{x\}$, $0\le \{ x\}<1$,
denotes the fractional part of $x$. In other words, the
spin$(\nu)$ of a configuration $\nu$ is equal to the one fourth
of the number of odd numbers appearing among the all differences
$(\nu^{(k-1)})_n'-(\nu^{(k)})_n'$, $k,n\ge 1$. Recall that by
definition $\nu^{(0)}=\mu$.

Finally, for any partition $\tau$ such that $|\tau|=|\ld|+|\mu|$
define the $(q,t)$--domino polynomial
\begin{equation}
c_{\Lambda}^{\tau}(q,t)=\sum_{\nu\in
C_{(2)}(\Lambda;\tau\oplus\tau)} t^{{\rm spin}(\nu)}q^{\wt
c(\nu)}\prod_{k,n\ge 1}\left[
\begin{array}{c}\wh P_n^{(k)}(\nu)+\wh m_n(\nu^{(k)})\\ \wh
m_n(\nu^{(k)})\end{array}\right]_q, \label{4.8*}
\end{equation}
where

$(\al)$ $\Lambda$ is a unique partition such that
2--core$(\Lambda)=\emptyset$ and 2--quotient$(\Lambda)=(\ld,\mu)$;

$(\beta)$ for any configuration $\nu\in
C_{(2)}(\Lambda;\tau\oplus\tau)$, $\wt c(\nu):=(c(\nu)-\epsilon
(\ld))/2$ denotes its normalized {\it charge};

$(\gamma)$ for any non--negative integer $m$ we put $\wh
m:=[m/2]$;

$(\delta)$ as usually, $\left[\begin{array}{c}m\\
n\end{array}\right]_q$ stands for the $q$--binomial coefficient.

$\bullet$ Show that polynomials $c_{\Lambda}^{\tau}(q,t)$ have the
following property:

$$c_{\Lambda}^{\tau}(1,1)=|{\rm Tab}^{(2)}(\Lambda,\tau)|,$$
i.e. $c_{\Lambda}^{\tau}(1,1)$ is equal to the number of
semistandard domino tableaux of shape $\Lambda$ and weight $\tau$.



The formula (\ref{4.8*}) with $q=1$ gives very effective method
for computing the {\it spin} generating function on the set of
domino tableaux.

ii) Denote by $\Ss$ the Sch\"utzenberger transformation acting on
the set of semistandard Young tableaux of the same shape:
$$\Ss : STY(\ld,\eta)\to STY(\ld, \overleftarrow\eta),
$$
where for any composition $\eta =(\eta_1,\ldots ,\eta_s)$ we
denote by $\overleftarrow\eta$ the "reverse" composition
$(\eta_s,\eta_{s-1},\ldots ,\eta_2,\eta_1)$. See e.g. \cite{Ful},
Appendix, where the transformation $\Ss$ is called {\it
evacuation}, or \cite{Kir1,Kir6}.

Let $\Lambda$ and $\mu$ be partitions such that
$|\Lambda|=2|\mu|$, and 2-{\rm core}$(\Lambda)=\emptyset$.

$\bullet$ Show that there exists a natural bijection
\begin{equation}
\theta : STY(\Lambda,\mu\vee\mu)^{\Ss}\simeq{\rm
Tab}^{(2)}(\Lambda,\mu), \label{4.10}
\end{equation}
where for any composition $\mu =(\mu_1,\ldots ,\mu_p)$ the symbol
$\mu\vee\mu$ denotes the composition $(\mu_1,\mu_2,\ldots
,\mu_{p-1},\mu_p,\mu_p,\mu_{p-1},\ldots ,\mu_2,\mu_1)$. In other
words, the set of domino tableaux of shape $\Lambda$ and weight
$\mu$ is in one--to--one correspondence with the set of
semistandard self--evacuating Young tableaux of shape $\Lambda$
and weight $\mu\vee\mu$. Recall, see e.g., J.~Stembridge (Duke
Math. J. {\bf 82} (1996), 585--606) that a tableau $T\in
STY(\ld,\mu)$ is called {\it self--evacuating} if $T$ is
invariant with respect to the action of Sch\"utzenberger's
transformation $\Ss$.

As far as I am aware, the statement that the set of
self--evacuating Young tableaux  of a given shape $\Lambda$ and
weight $\mu\vee\mu$ is equinumerous to that of domino tableaux of
the same shape $\Lambda$ and weight $\mu$ was conjectured by
R.~Stanley and has been proved for the first time by
J.~Stembridge [{\it ibid}]. An "elementary proof" has been found
later by A.D.~Berenstein and A.N.~Kirillov (Proc. of 10th Intern.
Conf. FPSAC, Fields Institute, Toronto, 1998, p.55--66, and
q--alg/9709010). Explicit construction of a bijection between the
sets of self--evacuating and domino tableaux was given by
S.~Fomin (Zap. Nauch. Sem. LOMI {\bf 55} (1986), 156--175) in the
case $\mu =(1^n)$, and was generalized by A.N.~Kirillov,
A.~Lascoux, B.~Leclerc and J.-Y.~Thibon (C.R. Acad. Sci. Paris
{\bf 318} (1994), 395--400) for general $\mu$.

$\bullet$ Show that if a domino tableaux ${\cal T}\in{\rm
Tab}^{(2)}(\Lambda,\mu)$ corresponds to the tableau $T\in
STY(\Lambda,\mu\vee\mu)^{\Ss}$ under the bijection (\ref{4.10})
constructed by A.N.~Kirillov {\it et al.} [{\it ibid}], and the
semistandard tableau $T$ corresponds to the pair $(\nu;J)$ under
the rigged configurations bijection, then
\begin{equation}{\rm spin}({\cal T})={\rm
spin}(\nu),
\end{equation}
where the spin$(\nu)$ of a configuration $\nu$ is defined by
(\ref{4.7}).

\begin{ex} {\rm Take $\ld =(321)$, $\mu=(3)$ and $\nu =(4321)$.
Then $\Lambda=(77411)$, $\wt\ld =(7321)$, $\wt R=(43321)$, and
$K_{\Lambda ,\nu\oplus\nu}(1)=2642$, $K_{\wt\ld \wt
R}(q)=q^3(1367531)$. It is not very difficult to check that
\begin{eqnarray*}
&&K_{\Lambda ,\nu\oplus\nu}(-1)=|{\rm Tab}^{(2)}(\Lambda ,\nu)|
=K_{\wt\ld \wt R}(1)=26,\\
&&|C(\Lambda ;\nu\oplus\nu)|=47, ~~~|C_{(2)}(\Lambda
;\nu\oplus\nu)|=19,
\end{eqnarray*}
and $(q,t)$--domino polynomial is equal to
\begin{eqnarray*}
&&q^4t+[q^4+2q^5+2q^6+q^6(1+q)]t^2+
[3q^5+q^5(1+q)+3q^6+q^6(1+q)\\
&&~~~~+q^9(1+q)]t^3+[q^7(1+q)+2q^8(1+q)]t^4\\
&&=q^4t+q^4t^2[1+q(1+q)(2+q)]
+q^5t^3(1+q)[4+q+q^4]\\
&&~~~~+q^7t^4(1+q)[1+2q].
\end{eqnarray*}
}
\end{ex}

iii) Let $\ld$ and $\mu$ be partitions such that 2-{\rm
core}$(\ld)=\emptyset$, and $|\ld|=2|\mu|$.

$\bullet$ Show that for any admissible configuration $\nu\in
C(\ld;\mu\oplus\mu)$
$$2{\rm spin}(\nu)\in\Z_{\ge 0}~~{\rm and}~~ 2{\rm
spin}(\nu)\equiv c(\nu)({\rm mod~}2).
$$

Assume additionally that $\ld$ has a form $2\nu\oplus 2\nu$ for
some partition $\nu$, or equivalently, that
2--quotient$(\ld)=(\nu,\nu)$.

$\bullet$ Show that if $\nu\in C_{(2)}(\ld;\mu\oplus\mu)$, then
$${\rm spin}(\nu)\in\Z_{\ge 0}.
$$

iv) (Sch\"utzenberger transformation and rigged configurations)\\
Let $\ld$ be a partition, $R$ be a sequence of rectangular shape
partitions, and $T\in LR(\ld,R)$ be a Littlewood--Richardson
tableau. Denote by $\Ss (T)\in LR(\ld,{\overleftarrow R})$ the
image of tableau $T$ under the Sch\"utzenberger transformation,
see e.g. \cite{KSS}.

$\bullet$ Show that if the pair $(\nu,\{ J^{(k)}_{n,\al}\})$
corresponds to the tableau $T$ under the rigged configurations
bijection, then the pair
\begin{equation}
(\nu,\{P_n^{(k)}(\nu;R)-J^{(k)}_{n,m_n(\nu^{(k)})-\alpha +1}\})
\end{equation}
corresponds to the tableau $\Ss(T)$.

This result was discovered and proved by the author in some
particular cases (including the case of semistandard tableaux). A
proof in general rectangular case has been obtained by
A.N.~Kirillov, A.~Schilling and M.~Shimozono \cite{KSS}.

v) Let $\ld$ and $\mu$, $l(\mu)=p$, be partitions, and $T\in
STY(\ld,\mu)$ be a semistandard Young tableau of shape $\ld$ and
weight $\mu$. With tableau $T$ one can associate a word $w:=w(T)$
which is given by reading the numbers in $T$ from right to left
in successive rows starting with the top row. With the word
$w:=w(T)$ one can associate a sequence of subwords $w_1,w_2,\ldots
,w_p$ such that weight $(w_i)=(1^{|w_i|})$ for $i=1,\ldots ,p$,
see e.g. \cite{Ma}, Chapter~III, Section~6, p.242.

Next, let $\nu$ be an admissible configuration corresponding to
the tableau $T$ under the rigged configurations bijection, see
e.g. \cite{Kir1}.

$\bullet$ Show that ~
$$(\nu^{(1)})_k'=\partial (w_k),~~ 1\le k\le p,$$
where for any {\it standard} word $w$ the symbol $\partial (w)$
denotes the number of integers $j$ such that $j$ and $j+1$ both
appear in the word $w$, but $j+1$ lies to the right of $j$. For
example, if $w=5761324$, then $\partial (w)=3$. If $w:=w(T)$,
then the numbers $\partial (w_k)$ are related to a structure of
{\it descent} set Des$(T)$ of tableau $T$, see e.g. \cite{Kir1}.

Recall \cite{Kir,Kir1,Kir6,KSS} that the rigged configurations
bijection establishes a one--to--one correspondence between the
set of Littlewood--Richardson tableaux $LR(\ld,R)$ and that of
rigged configurations $RC(\ld,R)$. This bijection was discovered
and constructed by the author in the middle of 80's, and since
that time many interesting combinatorial properties (see e.g.
\cite{Kir1,Kir2,Kir6} and the present paper) of rigged
configurations bijection were found by the author, G.-N.~Han,
S.~Fishel, M.~Kleber, A.~Schilling and M.~Schimozono.
\vskip 0.2cm

\hskip -0.6cm{\bf 4.} Let $p\ge 2$ be an integer, $\ld$ be a
partition such that $p$--core$(\ld)=\emptyset$, and $R$ be a
dominant sequence of rectangular shape partitions such that
$|\ld|=p|R|$. Denote by $R^{(p)}$ a dominant rearrangement of the
sequence of rectangular shape partitions
$\{\underbrace{R,R,\ldots ,R}_p\}$.

Denote by $C_{(p)}(\ld;R^{(p)})$ the subset of the set of
admissible configurations $C(\ld;R^{(p)})$ consisting of all
configurations $\nu\in C(\ld;R^{(p)})$ such that
$$\left[\frac{P_n^{(k)}(\nu;R)+m_n(\nu^{(k)})}{p}\right]=
\left[\frac{P_n^{(k)}(\nu;R)}{p}\right]+
\left[\frac{m_n(\nu^{(k)})}{p}\right]
$$
for all integers $k,n\ge 1$. Recall that for any $x\in\R$ the
symbol $[x]$ stands for the nearest integer which is smaller or
equal to $x$.

If $\nu\in C(\ld;R^{(p)})$, denote by $K_{\nu}(q)$ the following
product
\begin{equation}
K_{\nu}(q)=q^{c(\nu)}\prod_{k,n\ge 1}\left[\begin{array}{c}
P_n^{(k)}(\nu;R)+m_n(\nu^{(k)})\\
m_n(\nu^{(k)})\end{array}\right]_q, \label{4.11}
\end{equation}
and if $\nu\in C_{(p)}(\ld;R^{(p)})$, denote by
$K_{\nu}^{(p)}(q)$ the product
\begin{equation}
K_{\nu}^{(p)}(q)=q^{c(\nu)}\prod_{k,n\ge 1}\left[\begin{array}{c}
\wh P_n^{(k)}(\nu;R)+\wh m_n(\nu^{(k)})\\
\wh m_n(\nu^{(k)})\end{array}\right]_q, \label{4.12}
\end{equation}
where for any $m\in\Z$ we denote by $\wh m:=[m/p]$.

It follows from Theorem~\ref{t4.3} that
\begin{equation}
K_{\ld R}(q)=\sum_{\nu\in C(\ld;R^{(p)})}K_{\nu}(q). \label{4.13}
\end{equation}

$\bullet$ Show that if $\nu\in C_{(p)}(\ld;R^{(p)})$, then
$$K_{\nu}(\zeta_p)=\pm K_{\nu}^{(p)}(1)\in\Z.
$$
It looks a challenging task to describe explicitly the sign in
the latter formula.

$\bullet$ Show that if $p=2$ or 3, then the integer numbers
$K_{\nu}(\zeta_p)$ have the same sign for {\it all}
configurations $\nu\in C_{(p)}(\ld;R^{(p)})$. Note that this is
not true if $p\ge 4$.
\begin{con} \label{c4.7} Let $\ld$ and $\mu$ be partitions such
that $p$--core$(\ld )=\emptyset$ and $|\ld|=p|\mu|$. For any
configuration $\nu\in C_{(p)}(\ld;\mu^{(p)})$ one can define the
integer number {\rm spin}$(\nu)$ which is called by the {\it spin}
of configuration $\nu$, such that
$$\sum_{\nu\in C_{(p)}(\ld;\mu^{(p)})}t^{{\rm
spin}(\nu)}K_{\nu}^{(p)}(1)=\epsilon_p(\ld)\sum_{T\in
Tab^{(p)}(\ld,\mu)}t^{{\rm spin}(T)},
$$
where $\epsilon_p(\ld)=\pm 1$, and
$\mu^{(p)}=\mu\oplus\cdots\oplus\mu$ ($p$ times).
\end{con}
\vskip 0.2cm

\hskip -0.6cm {\bf 5.} Let $\ld$, $\mu
=(\mu_1\ge\cdots\ge\mu_p>0)$ and $\nu$ be partitions such that
$|\nu|=|\ld|+|\mu|$. In the sequel we will identify a partition
$\ld$ with its diagram, see e.g. \cite{Ma}, Chapter~I, p.2.

Let $x_1$ be  the rightmost square in the top row of $\ld$, and
let $x_2$ be the leftmost square in the bottom row of $\mu$. Let
$\mu^v$ (respectively, $\mu^h$) denote the diagram obtained from
$\mu$ by a shift sending $x_2$ to the square immediately above
$x_1$ (respectively, immediately right of $x_1$) and let $\ld
*\mu^v$ (respectively, $\ld *\mu^h$) denote the diagram
$\ld\cup\nu^v$ (respectively, $\ld\cup\mu^h$). Finally, let $\ld
*\mu$ denote the partition $((\ld_1^p)+\mu)\oplus\ld$ and
$\wt\nu$ (respectively, $\wt\nu_h$ and $\wt\nu_v$) denote a
dominant rearrangement of the sequence of rectangular shape
partitions $((\ld_1^p),\nu)$ (respectively, $(((\ld_1-1)^p),\nu)$
and $((\ld_1^{p-1}),\nu)$.

Show that
\begin{equation}
K_{\ld *\mu ,\wt\nu}(q)=K_{\ld *\mu^h ,\wt\nu_h}(q)+q^{\al (\ld
,\mu ,\nu )}K_{\ld *\mu^v ,\wt\nu_v}(q), \label{4.15}
\end{equation}
where $\al (\ld ,\mu ,\nu )$ is a non--negative integer which may
be computed explicitly.

The latter relation (\ref{4.15}) between parabolic Kostka
polynomials may be considered as a $q$--analog of the formula (1)
from \cite{Ma}, Chapter~I, Section~5, Exercise~21(a), in the case
when both $\theta$ and $\varphi$ (see [{\it ibid}]) are
partitions.

Note finally, see Section~\ref{pkp}, Exercise~3 ii), that
$$K_{\ld *\mu ,\wt\nu}(1)=|{\rm Tab}^{(2)}(\Lambda ,\nu)|,$$
where $\Lambda$ is a unique partition such that
2-core$(\Lambda)=\emptyset$ and 2-quotient$(\Lambda)=(\ld ,\mu)$.
It will be interesting to find an explanation of the relation
(\ref{4.15}) based on consideration domino tableaux only.

{\bf 6.} (Examples of fewnomial parabolic Kostka polynomials)

{\bf a.} Let $n\ge 2$ be an integer. Consider almost staircase
partition
$$\ld =(2n,2n-1,\ldots ,n,\wh{n-1},n-2,\ldots ,2,1),
$$
and sequence of rectangular shape partitions
$$R=((2n-1),(n^n),((n-1)^{n-1}),(1)).$$

i) Show that \vskip -0.5cm
$$K_{\ld R}(q)=(n-1)q^{\left(\begin{array}{c}n\\
2\end{array}\right)+1}(1+q)^2.
$$

ii) Let $k\ge 1$ be an integer. Consider partition
$$\ld^{(k)}=\ld\oplus\underbrace{(2n-1)\oplus\cdots\oplus
(2n-1)}_k,
$$
and sequence of rectangular shape partitions
$$R^{(k)}=(\underbrace{(2n-1),\ldots ,(2n-1)}_{k+1},(n^n),
((n-1)^{n-1}),(1)).
$$
Show that
\begin{eqnarray*}
K_{\ld^{(k)},R^{(k)}}(q)&=&q^{\left(\begin{array}{c}n\\
2\end{array}\right)+1}\left[
(n-1)+3(n-1)q+3n\frac{q^2-q^{k+1}}{1-q}\right.\\
&+&\left.(3n-4)q^{k+1}+(n-1)q^{k+2}\right].
\end{eqnarray*}
In particular, $\lim_{k\to\infty}K_{\ld^{(k)},R^{(k)}}(q)\bdoteq
\ds\frac{(n-1)+2(n-1)q+3q^2}{1-q}$.

{\bf b.} Let $n\ge 1$ be an integer. Consider almost staircase
partition
$$\ld =(2n+1,2n,\ldots ,n+1,{\wh n},n-1,\ldots ,2,1),
$$
and sequence of rectangular shape partitions
$$R=((2n),((n+1)^n),(n^{n-1}),(1)).$$

Show that \vskip -0.5cm
$$K_{\ld R}(q)=q^{\left(\begin{array}{c}n\\
2\end{array}\right)+1}(1+q)(n+(n-1)q).
$$

{\bf c.} Let $n\ge 1$ be an integer, and $\ld$ be the almost
staircase partition from {\bf b}. Consider a sequence of
rectangular shape partitions $R=((2n),(n^n),(n^n),(1)).$

Show that
$$K_{\ld R}(q)=q^{\left(\begin{array}{c}n\\
2\end{array}\right)+1}(n+(2n-1)q+nq^2).
$$

{\bf d.} Let $n\ge 2$ be an integer. Consider almost staircase
partition
$$\ld =(2n,\wh{2n-1},2n-2,\ldots ,n,\wh{n-1},n-2,\ldots
,2,1)
$$
and sequence of rectangular shape partitions
$$R=((2n-1),((n-1)^{n-1}),((n-1)^{n-1}),(1)).$$

i) Show that
$$K_{\ld R}(q)=q^{\left(\begin{array}{c}n\\
2\end{array}\right)+1}(n-1,2n-3,n-1).
$$

ii) Let $k\ge 1$ be an integer. Consider partition
$$\ld^{(k)}=\ld\oplus\underbrace{(2n-1)\oplus (2n-1)\cdots\oplus (2n-1)}_k,
$$
and sequence of rectangular shape partitions
$$R^{(k)}=(\underbrace{(2n-1),\ldots ,(2n-1)}_{k+1},((n-1)^{n-1}),
((n-1)^{n-1}),(1)).
$$

Show that
\begin{eqnarray*}
K_{\ld^{(k)},R^{(k)}}(q)&=&q^{\left(\begin{array}{c}n\\
2\end{array}\right)+1}\left[(n-1)+(3n-4)q+
3n\frac{q^2-q^{k+1}}{1-q}\right.\\
&+&\left.(3n-4)q^{k+1}+(n-1)q^{k+2}\right].
\end{eqnarray*}
In particular, $\lim_{k\to\infty}K_{\ld^{(k)},R^{(k)}}(q)\bdoteq
\ds\frac{(n-1)+(2n-3)q+4q^2}{1-q}.$

{\bf e.} Let $n\ge 3$ be an integer. Consider staircase partition
$$\ld =({2n+1,} 2n,\ldots ,2,1)
$$
and a sequence of rectangular shape partitions
$$R=((2n),(2n-1,2n-1),((n+1)^{n-1}),((n-1)^{n-2}),(2)).$$

Show that
$$K_{\ld R}(q)=2(n-2)q^{\left(\begin{array}{c}n-2\\
2\end{array}\right)+2}(1+q)^3.
$$

{\bf f.} i) Let $\rho_n=(2n-1,2n-3,\ldots ,3,1)$ be a partition
with parts which are consecutive odd numbers, starting from 1 till
$2n-1$. Consider partition $\ld =\rho_n\oplus (2n)\oplus (2n)$,
and a sequence of rectangular shape partitions
$$R=((2n-1),(\underbrace{n,\ldots
,n}_{\left[\frac{n+3}{2}\right]}), (\underbrace{n,\ldots
,n}_{\left[\frac{n+2}{2}\right]}),(1)).
$$
Show that
$$K_{\ld R}(q)=q^{\left[\frac{n+2}{2}\right]\left[\frac{n+3}{2}\right]}
\left(\left[\frac{n+2}{2}\right],n,\left[\frac{n+2}{2}\right]\right).
$$

ii) Consider partition $\ld =\rho_n\oplus (2n+1)\oplus (2n)$, and
a sequence of rectangular shape partitions
$$R=((2n),(\underbrace{n,\ldots
,n}_{\left[\frac{n+3}{2}\right]}), (\underbrace{n,\ldots
,n}_{\left[\frac{n+2}{2}\right]}),(1)).
$$
Show that
$$K_{\ld R}(q)=q^{\left[\frac{n+2}{2}\right]\left[\frac{n+3}{2}\right]}
\left(\left[\frac{n+2}{2}\right],n+1,\left[\frac{n+2}{2}\right]\right).
$$

iii) Let $\wt\rho_n$ denotes partition
$\rho_n+(\underbrace{0,\ldots ,0}_{n-1},1)$. Consider partition
$\ld =\wt\rho_n\oplus (2n)\oplus (2n)$, and a sequence of
rectangular shape partitions
$$R=((2n-1),(\underbrace{n,\ldots
,n}_{\left[\frac{n+3}{2}\right]}), (\underbrace{n,\ldots
,n}_{\left[\frac{n+2}{2}\right]}),(2)).
$$
Show that if $n\ge 3$, then
$$K_{\ld R}(q)=q^{\left[\frac{n+2}{2}\right]\left[\frac{n+3}{2}\right]}
\left(\left[\frac{n+2}{2}\right],\left[\frac{3n}{2}\right],
\left[\frac{3n}{2}\right],\left[\frac{n+2}{2}\right]\right).
$$

iv) Let $\ld$ be the same partition as in iii). Consider a
sequence of rectangular shape partitions
$$R=((2n-1),(\underbrace{n,\ldots
,n}_{\left[\frac{n+3}{2}\right]}), (\underbrace{n,\ldots
,n}_{\left[\frac{n+2}{2}\right]}),(1,1)).
$$
Show that if $n\ge 2$, then
$$K_{\ld R}(q)=q^{\left[\frac{n+2}{2}\right]\left[\frac{n+3}{2}\right]}
\left(\left[\frac{n}{2}\right],\left[\frac{3n-2}{2}\right],
\left[\frac{3n-2}{2}\right],\left[\frac{n}{2}\right]\right).
$$

\section{Generalized exponents and mixed tensor rep\-re\-sen\-tations}
\label{gemtr}
\neweq

Let $\g =sl(N,\C )$ denote the Lie algebra of all $N\times N$ complex
matrices of trace 0, and $G=GL(N,\C )$ denote the Lie group of all
invertible $N\times N$ complex matrices. Let
$${ad} : G\to{\rm Aut}(\g )$$
denote the adjoint representation of $G$, defined  by
$$(adX)(A)=XAX^{-1},$$
where $X\in G$, and $A\in\g$.

The adjoint action of $GL(N,\C)$ extends to the action on the symmetric
algebra
$$S^{\bullet}(\g )=\oplus_{k\ge 0}S^k(\g ),$$
where $S^k(\g)$ denotes the $k$--th symmetric power. It is well
known \cite{K} that the ring
$$I=S^{\bullet}(\g)^G=\{f\in S^{\bullet}(\g) |X\cdot f=f, \forall
x\in G\}
$$
of invariants of this action is a polynomial ring in
$N-1$ variables $f_2,\ldots ,f_{N-1}$, where $f_i\in S^i(\g )^G$.
By a theorem of Kostant \cite{K},
$$S^{\bullet}(\g) =I\otimes H$$
is a free module over $G$--invariants $I$ generated by harmonics
$H$. Moreover,
$$H=\oplus_{p\ge 0}H^p$$
is a graded (so $H^p=H\cap S^p(\g)$), locally finite
$\g$--representation. The graded character ${\rm ch}_q$ of the
symmetric algebra of adjoint representation is given by the
following formal power series
$${\rm ch}_q(S^{\bullet}(\g ))=\sum_{k\ge 0}q^k{\rm ch}(S^k(\g ))=
\prod_{1\le i,j\le n}(1-qx_i/x_j)^{-1}.
$$

For any finite dimensional $\g$--representation $V$ let us define
\begin{equation}G_N(V):=G_N(V;q)=\sum_{p\ge 0}\langle V, H^p\rangle q^p,
\label{5.1}
\end{equation}
where
$$\langle V_1,V_2\rangle ={\rm dim}_{\g}{\rm Hom}(V_1,V_2)$$
is the standard pairing on the representation ring of the Lie
algebra $\g$. By a theorem of Kostant \cite{K},
$$G_N(V)|_{q=1}={\rm dim}V(0),$$
where $V(0)$ denotes the zero weight subspace of representation
$V$. Hence, $G_N(V)$ is a polynomial in $q$ with  non--negative
integer coefficients. Follow Kostant \cite{K}, the integers
$e_1,\ldots ,e_s$ with
$$G_n(V)=\sum_{i=1}^sq^{e_i}$$
are called {\it generalized exponents} of the representation $V$.
Kostant's problem $[ibid]$ is to determine/compute these numbers
for a given representation $V$.

Let $V_{\ld}:=V_{\ld}^{[N]}$ denotes the irreducible highest weight
$\ld$ representation of the Lie algebra $\g := sl(N,\C )$.
Theorem~\ref{t5.0} below together with the fermionic formula (\ref{4.1a}) for the
Kostka--Foulkes polynomials, gives an effective method for computing the
generalized exponents of {\it irreducible} representation of the Lie
algebra $\g =sl(N)$.

\begin{theorem}\label{t5.0} {\rm (\cite{Gup2})} Let $\ld$ be a partition,
then
$$ G_N(V_{\ld})=\left\{\begin{array}{ll} K_{\ld
,\left(\left(\frac{|\ld|}{N}\right)^N\right)}(q),& ~{\rm if}~~|\ld|\equiv 0({\rm
mod}N),\\ 0,& ~{\rm otherwise}.\end{array}\right.
$$
\end{theorem}

For an "elementary" proof of Theorem~\ref{t5.0}, which is based
only on the theory of symmetric functions, see \cite{DLT}.

Using Theorem~\ref{t5.0}, one can compute, in principal, the
generalized exponents for any finite dimensional $\g
l(N)$--module $V$. What seems to be very interesting is that for
certain representations, see below, there exist alternative
expressions for the generalized exponents polynomials which have
independent interest and more convenient for computations. Before
turning to our main results of this Section, it is useful to
recall a few definitions and results from \cite{Gup,Br2,Gup2}.

Let $\al ,\beta$ be partitions, and $V_{\al}^{(N)}$,
$V_{\beta}^{(N)}$ be the highest weight $\al$ and $\beta$
(respectively) irreducible representations of the Lie algebra $\g
l(N)$. For any finite dimensional $\g l(N)$--module $V$ let $V^*$
denote its dual. If $l(\al)+l(\beta)\le N$, denote by $V_{\al
,\beta}^{(N)}$ the Cartan piece in the tensor product
$V_{\al}^{(N)}\otimes V_{\beta}^{(N)*}$, i.e. the irreducible
submodule generated by the tensor product of highest weight
vectors of each component. Follow \cite{Gup,Gup2}, we call a
representation obtained in this way by {\it mixed tensor
representation}. Clearly, $V_{\al ,\beta}^{(N)}$ is the dual of
$V_{\beta ,\al}^{(N)}$.

Since it is irreducible, $V_{\al ,\beta}^{(N)}$ is equal to
$V_{\ld}^{(N)}$ for a unique partition $\ld$ of less than $N$ rows. Let
us write $[\al ,\beta ]_N$ for this $\ld$. It is well--known
\cite{Gup,Gup2}, and goes back to D.E.~Littlewood \cite{Lt}, that
$$[\al,\beta]_N=(\al_1+\beta_1,\ldots
,\al_s+\beta_1,\underbrace{\beta_1,\ldots
,\beta_1}_{N-s-r},\beta_1-\beta_r,\ldots ,\beta_2-\beta_1,0), $$
where $s=l(\al)$, $r=l(\beta)$, and we assume that $s+r\le N$. For
example,
$$V_{0,0}^{(N)}=V_0^{(N)}\simeq\C,~~
V_{(1),(1)}^{(N)}=V_{(21^{N-2}) }^{(N)}=\g$$ is the adjoint
representation.

\begin{theorem}\label{t5.2} {\rm (\cite{Gup,Gup2})}
i) Fix $N\ge 1$. Then the representations
$V_{\alpha\beta}^{[N]}$, where partitions $\alpha$ and $\beta$ satisfy
$l(\alpha)+l(\beta)\le N$, and $|\alpha|=|\beta|$, form a complete,
repetition--free list of the irreducible finite--dimensional
representations of the group $PGL_N$;\\
ii) $G_N(V_{\alpha ,\beta})=K_{[\alpha ,\beta]_N,(\beta_1^N)}(q)$;\\
iii) $G_N(V_{\alpha}\otimes V_{\beta}^*):=G_N(V_{\alpha}^{[N]}
\otimes(V_{\beta}^{[N]})^*)\\
=\ds\sum_{\mu}K_{\alpha\mu}(q)K_{\beta\mu}(q)
\left[\begin{array}{c}N\\ N-\mu_1',\mu_1'-\mu_2',\ldots ,\mu_N'
\end{array}\right]_q$, where $\mu_i'$, $1\le i\le N$, denote the column
lengths of partition $\mu$, and
$$\left[\begin{array}{c}N\\ m_1,\ldots ,m_N
\end{array}\right]_q=\ds\frac{(q;q)_N}{(q;q)_{m_1}\cdots (q;q)_{m_N}}
$$
denotes the $q$--multinomial coefficient if $m_i\ge 0$, $1\le i\le N$,
$N=m_1+\cdots +m_N$, and 0 otherwise.
\end{theorem}

Our first result in this Section is Theorem~\ref{t5.3} below,
which connects the generalized exponents polynomial
$G_N(V_{\al}\otimes V_{\beta}^*)$ with a certain parabolic Kostka
polynomial, and gives, via Corollary~\ref{c5.4}, the first real
means for computing the $G_N(V_{\al}\otimes V_{\beta}^*)$.

\begin{theorem} \label{t5.3} Let $\al ,\beta$ be partitions,
$|\al|=|\beta|$, $l(\al )\le r$. Then
\begin{equation}
G_N(V_{\al}\otimes V_{\beta}^*)\bdoteq K_{[\al ,\beta ]_{N+r},R_N}(q),
\label{5.2}
\end{equation}
where $R_N=\{\underbrace{(\beta_1),\ldots
,(\beta_1)}_N,(\beta_1^r)\}$.
\end{theorem}

\begin{cor}\label{c5.4} {\rm (Fermionic formula for the generalized exponents
polynomial $G_N(V_{\alpha}\otimes V_{\beta}^*)$)} Let $\al ,\beta$ be
partitions, $|\al|=|\beta|$, $l(\al )\le r$. Then
\begin{equation}
q^{|\al|}G_N(V_{\al}\otimes V_{\beta}^*)=\sum_{\nu}q^{c(\nu )}
\prod_{k,j\ge 1}\left[\begin{array}{c}P_j^{(k)}(\nu)+m_j(\nu^{(k)})+
k\delta_{j,\beta_1}\theta (r-k)\\ P_j^{(k)}(\nu)\end{array}\right]_q,
\label{5.3}
\end{equation}
summed over all admissible configurations $\nu$ of type $([\al ,\beta
]_N;(\beta_1^N))$.
\end{cor}

See Definitions~\ref{d4.1} and \ref{d4.2} for explanation of symbols in
the RHS(\ref{5.3})

\begin{rem} {\rm Let $\alpha ,\beta$ be partitions such that
$l(\alpha)+l(\beta)\le N$, $l(\alpha)\le r$.
Then
$$G_N(V_{\alpha}\otimes V_{\beta}^*)\ne 0$$
if and only if $|\alpha|\equiv |\beta|{\rm mod}N$ and
$\wt\beta_1=\beta_1+\ds\frac{|\alpha|-|\beta|}{N} \ge 0$;\\
if so, then
\begin{equation}
G_N(V_{\alpha}\otimes V_{\beta}^*)=K_{[\alpha ,\beta]_{N+r},\wt R_N}(q),
\label{5.4}
\end{equation}
where $\wt R_N=\{\underbrace{(\wt\beta_1),\ldots
,(\wt\beta_1)}_N,(\beta_1^r)\}$.

Also we have
\begin{equation}
G_N(V_{\alpha}\otimes V_{\beta}^*)=\sum_{\mu}K_{\wt\alpha ,\mu}(q)
K_{\wt\beta ,\mu}(q)\left[\begin{array}{c}N\\ N-\mu_1',\mu_1'-\mu_2',
\ldots ,\mu_N'\end{array}\right]_q, \label{5.5}
\end{equation}
where $\wt\alpha =\alpha$, and $\wt\beta
=\beta\oplus\left(\left(\frac{|\alpha|-|\beta|}{N}\right)^N\right)$, if
$|\alpha|\ge|\beta|$, and $\wt\alpha =\alpha\oplus\left(\left(
\frac{|\beta|-|\alpha|}{N}\right)^N\right)$, and $\wt\beta =\beta$, if
$|\alpha|\le |\beta|$.}
\end{rem}

\begin{cor} \label{c5.6} Let $\al ,\beta$ be partitions, $|\al|=|\beta|$,
then
\begin{equation}
\sum_{\mu}\frac{K_{\al\mu}(q)K_{\beta\mu}(q)}{b_{\mu}(q)}=
\frac{K_{\beta\al}(q,q)}{H_{\al}(q)}, \label{5.6}
\end{equation}
where $K_{\beta\al}(q,q)$ denotes the specialization $q=t$ of the double
Kostka polynomial $K_{\beta\al}(q,t)$ introduced by Macdonald \cite{Ma},
Chapter~VI, \S 8;
$$b_{\mu}(q):=\prod_{i\ge 1}(q;q)_{\mu_i'-\mu_{i+1}'}.$$
\end{cor}

{\it Proof.} It follows from Theorem~\ref{t5.2}, $iii)$ that the
$${\rm LHS(\ref{5.3})}=\ds\lim_{N\to\infty}G_N(V_{\al}\otimes
V_{\beta}^*).
$$
If we compare the RHS(\ref{5.3}) with that of (\ref{6.10}), we
immediately see that
$$\lim_{N\to \infty}G_N(V_{\al}\otimes V_{\beta}^*)=q^{-|\al|}
\lim_{N\to\infty}s_{\al}*s_{\beta}(q,\ldots ,q^{N-1})=
s_{\al}*s_{\beta}(1,q,q^2,\ldots ),
$$
where $s_{\al}*s_{\beta}(1,q,\ldots )$ denotes the principal
specialization $x_i=q^{i-1}$, $1,2,3,\ldots$, of the internal
product of Schur functions $s_{\al}$ and $s_{\beta}$, see
Section~\ref{ipsf} for definition. It follows from \cite{Ma},
Chapter~VI, \S 8, Example~3, that
\begin{equation}
s_{\al}*s_{\beta}(1,q,q^2,\ldots )=\frac{K_{\al\beta}(q,q)}{H_{\beta}(q)}
=\frac{K_{\beta\al}(q,q)}{H_{\al}(q)}. \label{5.6a}
\end{equation}

\qed

\begin{rem}\label{r5.7} {\rm The identity (\ref{5.6}) is essentially
equivalent to the definition of Kostka--Foulkes polynomials
$K_{\ld\mu}(q)$, see  Section~\ref{kfpwm}, (\ref{2.1}). Indeed, by
definition,
$$s_{\beta}(x)=\sum_{\mu}K_{\beta\mu}(q)P_{\mu}(x;q)=
\sum_{\mu}\frac{K_{\beta\mu}(q)}{b_{\mu}(q)}Q_{\mu}(x;q).
$$
Therefore, in the $\ld$--ring notation, we have
$$s_{\beta}\left[\frac{X}{1-q}\right]=\sum_{\mu}
\frac{K_{\beta\mu}(q)}{b_{\mu}(q)}Q_{\mu}\left[\frac{X}{1-q};q\right].
$$
By definition,
$$Q_{\mu}\ds\left[\frac{X}{1-q};q\right]=Q_{\mu}'(x;q)$$
is the modified Hall-Littlewood polynomial, and
$$Q_{\mu}'(x;q)=\sum_{\al}K_{\al\mu}(q)s_{\al}(x).
$$
On the other hand, by definition (see \cite{Ma}, Chapter~VI, (8.11))
$$J_{\mu}(x;q,t)=\sum_{\ld}K_{\ld\mu}(q,t)s_{\ld}[X(1-t)],
$$
or
\begin{equation}
J_{\mu}\left[\frac{X}{1-t};q,t\right]=\sum_{\ld}K_{\ld\mu}(q,t)s_{\ld}(x).
\label{5.7}
\end{equation}
Now, if we put $q=t$ in (\ref{5.7}), we obtain
$$H_{\mu}(q)s_{\mu}\left[\frac{X}{1-q}\right]=\sum_{\ld}
K_{\ld\mu}(q,q)s_{\ld}(x).
$$
Hence,
$$\sum_{\al}\frac{K_{\al\beta}(q,q)}{H_{\beta}(q)}s_{\al}(x)=
s_{\beta}\left[\frac{X}{1-q}\right]=\sum_{\al}\sum_{\mu}
\frac{K_{\beta\mu}(q)K_{\al\mu}(q)}{b_{\mu}(q)}s_{\al}(x).
$$
These equalities imply (\ref{5.6}).}
\end{rem}

We conclude this Section with a generalization of
Theorem~\ref{t5.0}.

By a theorem of Kostant \cite{K} the symmetric algebra
$S^{\bullet}(\g )$ of the adjoint representation $\g$ of the Lie
algebra $\g l(N)$ is isomorphic as the graded $\g l(N)$--modules
to the direct sum of modules $V_{\mu}\otimes V_{\mu}^*$ when
$\mu$ ranges over the set ${\cal P}_N$ of all partitions whose
length do not exceed $N$:
$$S^{\bullet}(\g )\simeq\oplus_{\mu\in{\cal P}_N}(V_{\mu}\otimes
V_{\mu}^*).
$$
Let $l\ge 2$ be an integer, we are going to introduce and study
{\it $l$--restricted generalized exponents}. Consider the algebra
of $l$--restricted representations of the $q$--analog of universal
enveloping algebra $U_q(sl(N))$, when $q=\ds\exp\left(\frac{2\pi
i}{l+N}\right)$ is the primitive $(l+N)$--th root of unity. This
algebra, known as {\it Verlinde algebra} ${\cal V}(l,N)$, has
generators $V_{\ld}$, where $\ld$ ranges over the set ${\cal
P}_N^{(l)}$ of $l$--restricted partitions, i.e. partitions $\ld$
such that $l(\ld)<N$, $l(\ld')\le l$. The multiplication in the
algebra ${\cal V}(l,N)$ is given by the so--called {\it
$l$--restricted tensor product}, which we will denote by
$\wh\otimes$. We refer the reader to \cite{Kac}, Chapter~13,
Exercise~13.34, for definition and basic properties of the
Verlinde algebra.

\begin{de}\label{d5.8} Let $l\ge 2$ be an integer, define
\begin{eqnarray*}
S^{\bullet}_l(\g )&=&\bigoplus_{k\ge 0}S_l^k(\g ),~~~{\rm where}\\
S_l^k(\g )&=&\bigoplus_{|\mu|=k}(V_{\mu}\wh\otimes V_{\mu}^*),
\end{eqnarray*}
and $\mu$ ranges over the set ${\cal P}_N^{(l)}$ of $l$--restricted
partitions, and $|\mu|=k$.
\end{de}

\begin{de}\label{d5.9} {\rm ($l$--restricted generalized exponents)} Let
$V$ belongs to the Verlinde algebra ${\cal V}(l,N)$. The $l$--restricted
generalized exponents polynomial $G_N^{(l)}(V)$ is defined by the
following formula
$$G_N^{(l)}(V):=G_N^{(l)}(V;q)=\sum_{k\ge 0}\left[ V:S_l^{k}(\g )\right]q^k,
$$
where for any two elements $V$ and $W$ of the Verlinde algebra ${\cal V}(l,N)$
the symbol $[V:W]$ has the following meaning:

let $V=\sum_{\ld}a_{\ld}V_{\ld}$ and $W=\sum_{\mu}b_{\mu}V_{\mu}$ be the
decompositions of the elements $V$ and $W$ in terms of the generators
$V_{\ld}$ of the Verlinde algebra ${\cal V}(l,N)$. By
definition,
$$[V:W]=\sum_{\ld}a_{\ld}b_{\ld}.$$
\end{de}

\begin{de}\label{d510} {\rm (Verlinde polynomials)} Let $l\ge 2$ and
$L\ge 1$ be integers, and $V\in{\cal V}(l,N)$. Denote by $(S_l^{\bullet}(\g
)^{\wh\otimes L})^{(k)}$ the degree $k$ component of the $L$--th
$l$--restricted tensor power of the element $S_l^{\bullet}(\g )\in
{\cal V}(l,N)$. Polynomial
$$G_N^{(l,L)}(V):=\sum_{k\ge 0}\left[ V:(S_l^{\bullet}(\g )^{\wh\otimes L}
)^{(k)}\right]q^k
$$
is called by Verlinde polynomial.
\end{de}

Let us denote by $\wh{\cal P}_N^{(l)}$ the set of all partitions $\ld$
such that $l(\ld )\le N$ and $\ld_1-\ld_N\le l$.

\begin{prb} Let $\ld ,\al ,\beta$ be partitions from the set
$\wh{\cal P}_N^{(l)}$. Compute polynomials $G_N^{(l,L)}(V_{\ld})$ and
$G_N^{(l,L)}(V_{\al}\wh\otimes V^*_{\beta})$.
\end{prb}

\begin{theorem}\label{t5.12} Let $\ld\in\wh{\cal P}_N^{(l)}$, then
$$G_N^{(l)}(V_{\ld})=\left\{\begin{array}{ll}K^{(l)}_{\ld ,\left(\left(
\frac{|\ld |}{N}\right)^N\right)}(q),& {\rm if}~~~|\ld|\equiv 0({\rm
mod}N),\\ 0,& {\rm otherwise}.\end{array}\right.
$$
\end{theorem}

Here for a partition $\ld\in\wh{\cal P}_n^{(l)}$ and any partition
$\mu$ the symbol $K_{\ld ,\mu}^{(l)}(q)$ denotes {\it the
$l$--restricted Kostka--Foulkes polynomial}. We refer the reader
to \cite{HKKOTY}, Section~4.1, see also references therein, for a
definition of the so--called {\it $l$--restricted one dimensional
sums} $X_{\mu}^{(l)}(\ld)$ which we identify with $l$--restricted
Kostka--Foulkes polynomials
$$K_{\ld\mu}^{(l)}(q):=X_{\mu}^{(l)}(\ld).$$
In my opinion, this is one of the most "natural" definitions of
the $l$--restricted Kostka--Foulkes polynomials. Other
definitions appeal either to the Representation Theory or Quantum
Groups at roots of unity [folklore], or that of the
Hecke--Iwahori algebras (Goodman~F. and Wenzl~H., Adv. Math. {\bf
82} (1990), 244-265). Restricted Kostka--Foulkes polynomials have
many interesting combinatorial properties (restricted Young
tableaux, restricted Littlewood--Richardson rule, $\ldots$), and
appear to be connected with the characters of Virasoro and affine
Lie algebras, see e.g. \cite{Kir4}, p.101-105; \cite{HKKOTY},
Section~5, and references therein.







\begin{prb}\label{pb5.15} Let $\al$ and $\beta$ be partitions of the
same size from the set $\wh{\cal P}_N^{(l)}$. Define {\it
$l$--restricted double--Kostka} polynomial $K_{\al\beta}^{(l)}(q,t)$
with the following properties

i)  $K_{\al\beta}^{(l)}(q,t)$ is a polynomial with non--negative integer
coefficients;

ii) $K_{\al\beta}^{(l)}(0,t)$ is equal to the $l$--restricted
Kostka--Foulkes polynomial $K_{\al\beta}^{(l)}(t)$;


iii) if $l$ is big enough, then $K_{\al\beta}^{(l)}(q,t)$
coincides with the double Kostka polynomial $K_{\al\beta}(q,t)$
introduced by Macdonald \cite{Ma}, Chapter~VI.
\end{prb}
See recent preprint by L.~Lapointe, A.~Lascoux and J.~Morse, {\it
A filtration of the symmetric function space and refinement of the
Macdonald positivity conjecture}, math.QA/0008073, where certain
polynomials $K_{\ld\mu}^{(k)}(q,t)\le K_{\ld\mu}(q,t)$ which
might be satisfy the above conditions i) and ii), have been
introduced and studied.

\vskip 0.5cm
 \hskip -0.6cm {\bf Exercises} \vskip 0.2cm

\vskip 0.3cm \hskip -0.6cm {\bf 1.} Show that if $\beta =(n)$, then the
identity (\ref{5.6}) takes the following form
\begin{equation}
\sum_{\mu\vdash n}\frac{q^{n(\mu)}K_{\al ,\mu}(q)}{b_{\mu}(q)}=
\frac{K_{\al',(1^n)}(q)}{(q;q)_n}. \label{5.9}
\end{equation}
In the particular case $\al =(n)$, the identity (\ref{5.9}) becomes
\begin{equation}
\sum_{\mu\vdash n}\frac{q^{2n(\mu)}}{\prod_i(q;q)_{m_i(\mu)}}=
\frac{1}{(q;q)_n}. \label{5.10}
\end{equation}
This identity is due to Philip Hall (Comm. Math. Helv. {\bf 11}
(1938), 126-129).

\hskip -0.6cm {\bf 2.} Let $\al$ be a partition of $n$. Show that
$$\sum_{\mu\vdash n}\frac{q^{n(\mu)}K_{\al ,\mu}(q)
\prod_{j=1}^{\mu_1'}(1+zq^{1-j})}{b_{\mu}(q)}=\prod_{(i,j)\in\al}
(1+zq^{j-i})\frac{K_{\al',(1^n)}(q)}{(q;q)_n}.
$$
In the particular case $\al =(n)$, we obtain a generalization of
the Philip Hall identity (\ref{5.10}):
$$\sum_{\mu\vdash n}\frac{q^{2n(\mu)}\prod_{j=1}^{\mu_1'}
(1+zq^{1-j})}{b_{\mu}(q)}=\frac{(-z;q)_n}{(q;q)_n}.
$$
{\bf 3.} Fusion rules and Catalan numbers.

Let $q^{n+2}=1$, $n\equiv 1({\rm mod}2)$ be a primitive root of 1
of odd order. For each integer $m$ such that $0\le 2m+1\le n$,
denote by $V_m$ the $(2m+1)$--dimensional irreducible
representation of $U_q(sl(2))$. Show that
\begin{eqnarray*}
&&{\rm Mult}[V_0:V_{\frac{n-1}{2}}^{\wh\otimes 2k}]=C_k,~~{\rm
if}~~k\le n-2,\\ &&{\rm Mult}[V_0:V_{\frac{n-1}{2}}^{\wh\otimes
2k}]<C_k,~~{\rm if}~~k\ge n-1.
\end{eqnarray*}
Here $C_k$ denotes the $k$--th Catalan number. In other words, if
$n\ge k+2$, then the Catalan number $C_k$ is equal to the
$(n-2)$--restricted Kostka number $K_{\ld\mu}^{(n-2)}$, where
$\ld=(k(n-1)/2,k(n-1)/2)$ and $\mu =(((n-1)/2)^{2k})$.

Using the fermionic formula (\ref{3.13*}) for restricted Kostka
polynomials, show that
$$C_k=\sum_{\nu}\prod_{j\ge 1}\left(\begin{array}{c}P_j(\nu)
+m_j(\nu)\\ m_j(\nu)\end{array}\right),
$$
summed over all partitions $\nu$ such that

$\bullet$ $|\nu|=k(n-1)/2$, $\nu_1\le n-2$,

$\bullet$ $P_j(\nu):=k\min(2j,n-1)-2Q_j(\nu)\ge 0$.

Thus, if $n\ge k+2$, then $(n-2)$--restricted Kostka polynomial
$K_{\ld\mu}^{(n-2)}(q)$ with $\ld =((k(n-1)/2)^2)$ and
$\mu=(((n-1)/2)^{2k})$, see (\ref{3.13*}), may be considered as a
$q$--analog of the Catalan number $C_k$.

\hskip -0.6cm{\bf 4.} Generalized exponents for the Lie algebras
of type $A_2$ and $A_3$.

{\bf a.} Let $\ld =(\ld_1\ge\ld_2\ge 0)$ and $\mu =(l^3)$ be
partitions, $|\ld|=|\mu|$.

Show that
$$G_{\ld}^{(3)}(q)=K_{\ld ,(l^3)}(q)=q^c\left[\begin{array}{c}
\min (\ld_1-\ld_2,\ld_2)+1\\ 1\end{array}\right]_q,
$$
where $c=\max (\ld_1-\ld_2,\ld_2)$.

{\bf b.} Let $\ld =(\ld_1\ge\ld_2\ge\ld_3\ge 0)$ and $\mu =(l^4)$
be partitions, $|\ld|=|\mu|$. Show that
\begin{eqnarray*}
&&q^{2l(l-1)-n(\ld)}G_{\ld}^{(4)}(q^{-1})=\wt K_{\ld ,(l^4)}(q)=
\left[\begin{array}{c} \min(\ld_1-\ld_2,\ld_3)+1\\
1\end{array}\right]_q  \\ &&\times\sum_{0\le
a<\frac{\ld_2-\ld_3}{2}}q^{2a}\left[\begin{array}{c}
\max(\ld_1-\ld_2,\ld_3)+4\min(l-\ld_3-a,\ld_2-l-a,0)+1\\
1\end{array}\right]_q\\ &&+\epsilon q^{\ld_2-\ld_3}\left[
\begin{array}{c}\min(\ld_1-\ld_2,\ld_3)+2\\ 2\end{array}\right]_q,
\end{eqnarray*}
where $\epsilon =\ds\frac{1+(-1)^{\ld_2-\ld_3}}{2}$.

{\bf c.} Show that if $\ld =(\ld_1\ge\ld_2\ge\cdots\ge\ld_N\ge 0)$
is a partition such that $\ld_N\ge l$, and $|\ld|=(N+l)l$, then
$$G_{\ld}^{(N+1)}(q)=q^{n(\al)+l}\left[\begin{array}{c}N\\
\al'\end{array}\right],
$$
where $\al =(\ld_1-l,\ld_2-l,\ldots ,\ld_N-l)$.

{\bf d.} Let $\ld =(\ld_1\ge\cdots\ge\ld_N\ge 0)$ be partition
such that $l\ge\ld_2$, and $|\ld|=(N+1)l$. Show that
$$G_{\ld}^{(N+1)}(q)=q^{l\left(\begin{array}{c}N+1\\
2\end{array}\right)-N|\beta|+n(\beta)}\left[\begin{array}{c} N\\
\beta'\end{array}\right],
$$
where $\beta =(\ld_2,\ld_3,\ldots ,\ld_N)$.

{\bf e.} Let $\ld =(\ld_1\ge\ld_2\ge\ld_3\ge 0)$ be partition such
that $\ld_1=\ld_2+\ld_3$, and $|\ld|=4l$. Show that
$$\wt K_{\ld ,(l^4)}(q)=\frac{(1-q^{\ld_2-\ld_3+1})
(1-q^{\ld_3+1})^2}{(1-q)^2(1-q^2)}+\epsilon q^{\ld_2+1}
\frac{1-q^{\ld_3+1}}{1-q^2},
$$
where $\epsilon =\ds\frac{1+(-1)^{\ld_1}}{2}$.

\section{Internal product of Schur functions}
\label{ipsf}
\neweq

The irreducible characters $\chi^{\ld}$ of the symmetric group $S_n$ are
indexed in a natural way by partitions $\ld$ of $n$. If $w\in S_n$, then
define $\rho (w)$ to be the partition of $n$ whose parts are the cycle
lengths of $w$. For any partition $\ld$ of $m$ of length $l$, define the
power--sum symmetric function
$$p_{\ld}=p_{\ld_1}\ldots p_{\ld_l},$$
where $p_n(x)=\sum x_i^n$. For brevity write $p_w:=p_{\rho (w)}$.
The Schur functions $s_{\ld}$ and power--sums $p_{\mu}$ are
related by a famous result of Frobenius
\begin{equation}
s_{\ld}=\frac{1}{n!}\sum_{w\in S_n}\chi^{\ld}(w)p_w. \label{6.1}
\end{equation}
For a pair of partitions $\al$ and $\beta$, $|\al|=|\beta|=n$, let us
define the internal product  $s_{\al}*s_{\beta}$ of Schur functions
$s_{\al}$ and $s_{\beta}$:
\begin{equation}
s_{\al}*s_{\beta}=\frac{1}{n!}\sum_{w\in S_n}\chi^{\al}(w)\chi^{\beta}(w)p_w.
\label{6.2}
\end{equation}
It is well--known that
$$s_{\al}*s_{(n)}=s_{\al},~~ s_{\al}*s_{(1^n)}=s_{\al'},
$$
where $\al'$ denotes the conjugate partition to $\al$.

Let $\al ,\beta ,\gamma$ be partitions of a natural number $n\ge 1$,
consider the following numbers
\begin{equation}
g_{\al\beta\gamma}=\frac{1}{n!}\sum_{w\in S_n} \chi^{\al}(w)
\chi^{\beta}(w)\chi^{\gamma}(w). \label{6.3}
\end{equation}
The numbers $g_{\al\beta\gamma}$ coincide with the structural constants
for multiplication of the characters $\chi^{\al}$ of the symmetric group
$S_n$:
\begin{equation}
\chi^{\al}\chi^{\beta}=\sum_{\gamma}g_{\al\beta\gamma}\chi^{\gamma}.
\label{6.4}
\end{equation}
Hence, $g_{\al\beta\gamma}$ are non--negative integers. It is clear that
\begin{equation}
s_{\al}*s_{\beta}=\sum_{\gamma}g_{\al\beta\gamma}s_{\gamma}. \label{6.5}
\end{equation}
Let us introduce polynomials $L_{\al\beta}^{\mu}(q)$ via the
decomposition of the internal product of Schur functions
$s_{\al}*s_{\beta}(x)$ in terms of Hall--Littlewood functions:
\begin{equation}
s_{\al}*s_{\beta}(x)=\sum_{\mu}L_{\al\beta}^{\mu}(q)P_{\mu}(x;q).
\label{6.6}
\end{equation}

It follows from (\ref{2.1}) and (\ref{6.5}) that
\begin{equation}
L_{\al\beta}^{\mu}(q)=\sum_{\gamma}g_{\al\beta\gamma}K_{\gamma\mu}(q).
\label{6.7}
\end{equation}
Thus, the polynomials $L_{\al\beta}^{\mu}(q)$ have non--negative integer
coefficients, and
$$L_{\al\beta}^{\mu}(0)=g_{\al\beta\mu}.$$
The polynomials $L_{\al\beta}^{\mu}(q)$ can be considered as a
generalization of the Kostka--Foulkes polynomials. Indeed, if
partition $\beta$ consists  of one part, $\beta =(n)$, then
$$L_{\al \beta}^{\mu}(q)=K_{\al ,\mu}(q).
$$

\begin{ex}\label{e6.1} {\rm Take partitions $\al =(4,2)$ and $\beta
=(3,2,1)$, then
$$s_{\al}*s_{\beta}=s_{51}+2s_{42}+2s_{41^2}+s_{3^2}+3s_{321}+2s_{31^3}
+s_{2^3}+2s_{2^21^2}+s_{21^4}.
$$
Therefore,
\begin{eqnarray*}
&&L_{\al\beta}^{(51)}(q)=1,~~L_{\al\beta}^{(42)}(q)=(21),~~
L_{\al\beta}^{(411)}(q)=(231),~~L_{\al\beta}^{(33)}(q)=(121),\\
&&L_{\al\beta}^{(321)}(q)=(3531),~~
L_{\al\beta}^{(31^3)}(q)=(257531),\\
&&L_{\al\beta}^{(2^3)}(q)=(135531),~~
L_{\al\beta}^{(2^21^2)}=(2699631),\\
&&L_{\al\beta}^{(21^4)}(q)=(1,4,8,12,14,13,10,6,3,1)\\
&&~~~~~~~~~~~\,=(1+q)^2(1+q^2)(1+q+q^2)(1+q+q^3),\\
&&L_{\al\beta}^{(1^6)}(q)=q(1,3,6,10,14,18,20,20,18,14,10,6,3,1)\\
&&~~~~~~~~~~\,= q(1+q)^3(1+q^2)(1+q^2+q^4)^2.
\end{eqnarray*}
}
\end{ex}

The value of polynomial $L_{\al\beta}^{\mu}(q)$ at $q=1$ admits
the following combinatorial interpretation. First of all, for any
skew shapes $A$ and $B$ consider the set ${\bnu}(A,B)$ of all
$B$--tableaux of shape $A$. Let $b_1,b_2,\ldots$ be the row
lengths of skew shape $B$, recall that a semistandard Young
tableau $T$ of (skew) shape $A$ and weight $(b_1,b_2,\ldots)$ is
called {\it $B$--tableau}, if the word $w(T)$ associated to $T$
(and given by reading the numbers in $T$ from right to left in
successive rows, starting with the top row) is a $B$--lattice
one.

Recall that for a given skew shape $B=\ld\setminus\mu$, a word
$w=w_1w_2\cdots w_p$, $w_i\in\Z_{>0}$, is called $B$--lattice, if for any
$j$, $1\le j\le p$, such that $w_j=k>1$, the following inequality holds
\begin{eqnarray}
&\#\{ i|i<j~~{\rm and}~~w_i=w_j-1\} +\min
(\ld_k-\mu_k,\mu_{k-1}-\mu_k)\nonumber\\
&\ge \#\{ i|i\le j~~{\rm and}~~w_i=w_j\}. \label{6.8a}
\end{eqnarray}
In particular, if $B$ is a {\it diagram} (i.e. $\mu =0$) then the
definitions (\ref{6.8a}) and (\ref{2.12a}) of {\it $B$--lattice} words and
{\it lattice} words coincide.

Summarizing,
$${\bnu}(A,B)=\{ T\in STY(A,(b_1,b_2,\ldots ))|w(T)~~{\rm is~~a}~~
\hbox{\rm $B$--lattice~~word}\}.
$$
Let $\nu (A,B)$ be the number of $B$--tableaux of shape $A$. It follows
from \cite{Do} and \cite{Zl} that
$$\nu (A,B)=\langle s_A,s_B\rangle ={\rm dim}_{\g l(n)}(V_A,V_B).
$$

Note that if $A=\ld\setminus\mu$ and $B$ is a {\it diagram}, say
$B=\eta$, then $\nu(A,B)$ is equal to the Littlewood--Richardson number
$c_{\mu\eta}^{\ld}$, i.e. the number of semistandard Young tableaux $T$
of the (skew) shape $A$ and weight $\eta$ such that the corresponding
word $w(T)$ is a lattice word.

\begin{ex}\label{e6.2} {\rm Take $A=(4332)\setminus (211)$,
$B=(5321)\setminus (21)$. Then
$$K_{A,(3221)}(q)=\sum_{\pi}c_{211,\pi}^{4332}K_{\pi ,3221}(q)=
\sum_{T\in STY(A,(3221))}q^{c(T)}=2+4q+3q^2+q^3.
$$
On the other hand, $\nu (A,B)=8$. It is easy to check that two
semistandard Young tableaux $T_1,T_2$ of shape $A$ and weight (3221),
which correspond to the words $w(T_1)=11214233$ and $w(T_2)=41213231$,
do not belong to the set ${\bnu}(A,B)$.}
\end{ex}

\begin{pr}\label{p6.3} {\rm (\cite{Do})} Let $\al ,\beta$ be partitions,
$\mu$ be composition and $|\al|=|\beta|=|\mu|$. Then
\begin{equation}
L_{\al\beta}^{\mu}(1)=\sum_{\scriptsize\begin{array}{c}
\al=\coprod_iA_i\\ \beta =\coprod_iB_i\\ |A_i|=|B_i|=\mu_i\end{array}}
\prod_{i\ge 1}\nu (A_i,B_i), \label{6.7a}
\end{equation}
\end{pr}
Here the symbol $\ds\coprod_iA_i$ stands for the disjoint union of
sets $A_i$.

\begin{de}\label{d6.1} Let $\al ,\beta$ be partitions, $\mu$ be
composition, and $|\al|=|\beta|=|\mu|$, $l(\mu)=r$. A sequence of skew
semistandard Young tableaux $(T_1,\ldots ,T_r)$ is called {\it
Littlewood--Richardson sequence of type $(\al ,\beta ;\mu)$}, if there
exist two sequences of partitions
$$\{0=\al_0\subset\al_1\subset\al_2\subset\cdots\subset\al_r=\al\}~~{\rm
and}~~\{0=\beta_0\subset\beta_1\subset\cdots\subset\beta_r=\beta\},
$$
such that

i) ~$|\al_i|-|\al_{i-1}|=|\beta_i|-|\beta_{i-1}|=\mu_i$, $1\le i\le r$;

ii) ~$T_i$ is a $(\beta_i\setminus\beta_{i-1})$--tableau of (skew) shape
$\al_i\setminus\al_{i-1}$, $1\le i\le r$.
\end{de}

It is clear from (\ref{6.7a}) that $L_{\al\beta}^{\mu}(1)$ is
equal to the number of Littlewood--Richardson sequences of
tableaux $(T_1,\ldots ,T_r)$ of type $(\al ,\beta ;\mu)$. We will
denote this set by ${\bnu}^{\mu}(\al ,\beta)$ and its cardinality
by $\nu^{\mu}(\al,\beta)$. From the very definition,
$$L_{\al\beta}^{\mu}(1)=\nu^{\mu}(\al,\beta)=
\sum_{\ld}g_{\al\beta\ld}K_{\ld\mu}.
$$
It is not difficult to see that
$$\nu^{(1^n)}(\al,\beta)=f^{\al}f^{\beta},$$
where $f^{\al}$ denotes the number of standard Young tableaux of
shape $\al$. For a $q$--analog of the latter formula, see
Exercise~1 to this Section.

\begin{prb}\label{prb6.5} Find a statistics $d$ on the set
$\bnu^{\mu}(\al ,\beta)$ with the generating function
$L_{\al\beta}^{\mu}(q)$.
\end{prb}

Let $N\ge 2$, consider the principal specialization $x_i=q^i$, $1\le i\le
N-1$, and $x_i=0$, if $i\ge N$, of the internal product of Schur functions
$s_{\al}$ and $s_{\beta}$:
\begin{equation}
s_{\al}*s_{\beta}(q,q^2,\ldots ,q^{N-1})=\frac{1}{n!}\sum_{w\in S_n}
\chi^{\al}(w)\chi^{\beta}(w)\prod_{k\ge 1}\left(
\frac{q^k-q^{kN}}{1-q^k}\right)^{\rho_k(w)}. \label{6.9a}
\end{equation}

By a result of R.-K.~Brylinski \cite{Br2}, Corollary~5.3, the
polynomials
$$s_{\al}*s_{\beta}(q,\ldots ,q^{N-1})
$$
admit the following interpretation. Let $P_{n,N}$ denote the
variety of $n$ by $n$ complex matrices $z$ such that $z^N=0$.
Denote by
$$R_{n,N}:=\C [P_{n,N}]$$
the coordinate ring of polynomial functions on $P_{n,N}$ with
values in the field of complex numbers $\C$. This is a graded
ring:
$$R_{n,N}=\ds\oplus_{k\ge 0}R_{n,N}^{(k)},$$
where $R_{n,N}^{(k)}$ is a finite dimensional $\g l(n)$--module
with respect to the adjoint action. Let $\al$ and $\beta$ be
partitions of common size. Then \cite{Br2}
$$s_{\al}*s_{\beta}(q,\ldots ,q^{N-1})=\sum_{k\ge 0}\langle V_{[\al ,
\beta]_n},P_{n,N}^{(k)}\rangle q^k,
$$
as long as $n\ge\max (Nl(\al),Nl(\beta),l(\al)+l(\beta))$.

Our main result of this Section is Theorem~\ref{t6.1} below, which
connects the principal specialization of the internal product of
Schur functions with certain parabolic Kostka polynomials, and
gives, via Corollary~\ref{c6.7}, an effective method for
computing the polynomials $s_{\al}*s_{\beta}(q,\ldots ,q^{N-1})$
which for the first time does not use the character table of the
symmetric group $S_n$.

\begin{theorem}\label{t6.1} i) Let $\al ,\beta$ be partitions, $|\al|=|\beta|$,
$l(\al)\le r$, and $l(\al)+l(\beta)\le Nr$. Consider the sequence 
of rectangular shape partitions
$$R_N=\{\underbrace{(\beta_1^r),\ldots ,(\beta_1^r)}_N\}.
$$
Then
\begin{equation}
s_{\al}*s_{\beta}(q,\ldots ,q^{N-1})\bdoteq K_{[\al
,\beta]_{Nr},R_N}(q). \label{4.8}
\end{equation}

ii) (Dual form) Let $\al ,\beta$ be  partitions such that $|\al
|=|\beta |$, $l(\al)\le r$ and $l(\beta )\le k$. For given
non--negative integer $N$ such that $\al_1+\beta_1\le Nr$,
consider partition
$$\ld_N:=(rN-\beta_k,rN-\beta_{k-1},\ldots ,rN-\beta_1,\al)
$$
and a sequence of rectangular shape partitions
$$R_N:=(\underbrace{(r^k),\ldots ,(r^k)}_N).$$
Then
\begin{equation}
K_{\ld_NR_N}(q)\bdoteq s_{\al}*s_{\beta}(q,\ldots ,q^{N-1})
\label{6.11*}
\end{equation}
\end{theorem}

\begin{cor}\label{c6.7} {\rm (Fermionic formula for the principal specialization
of the internal product of Schur functions, \cite{Kir2,Kir6})} Let
$\al ,\beta$ be partitions, $|\al|=|\beta|$, $l(\al)=r$. Then
\begin{eqnarray}
&&s_{\al}*s_{\beta}(q,\ldots ,q^{N-1})\label{6.9}\\
&&=\sum_{\nu}q^{c(\nu)}\prod_{k,j\ge 1}\left[\begin{array}{c}
P_j^{(k)}(\nu)+m_j(\nu^{(k)})+N(k-1)\delta_{j,\beta_1}\theta (r-k)\\
P_j^{(k)}(\nu )\end{array}\right]_q,\nonumber
\end{eqnarray}
where the sum runs over all admissible configurations $\nu$ of type $([\al
,\beta ]_N,(\beta_1^N))$.
\end{cor}

See Definitions~\ref{d4.1} and \ref{d4.2} for explanation of symbols in
the RHS(\ref{6.9}).

\begin{cor}\label{c6.8} {\rm (\cite{Kir2})} Let $\al ,\beta$ be partitions,
$|\al|=|\beta|$. Then the principal specialization
$$s_{\al}*s_{\beta}(q,\ldots ,q^{N-1})
$$
is a symmetric and unimodal polynomial.
\end{cor}

At the end of this Section we discuss a connection between the structural
constants $g_{\al\beta\gamma}$ and the
decomposition into irreducible parts of the exterior algebra
$$\Lambda^{\bullet}(\g )=\ds\oplus_{k\ge 0}\Lambda^k(\g )$$
of the adjoint representation $\g$ of the Lie algebra $\g l(n)$.

To start, let us recall the well--known formula for the graded character
${\rm ch}_q(\Lambda^{\bullet}(\g ))$ of the exterior algebra of adjoint
representation:
$${\rm ch}_q(\Lambda^{\bullet}(\g ))=\sum_{k\ge 0}q^k{\rm ch}
(\Lambda^k(\g ))=\prod_{1\le i,j\le n}(1+qx_i/x_j).
$$
For any finite dimensional $\g l(n)$--module $V$, let us
consider the following polynomial
$$E_n(V):=E_n(V;q)=\sum_{k\ge 0}{\rm Mult}[V:\Lambda^k(\g )]q^k.
$$

\begin{prb} Let $V_{\ld}$ and $V_{\mu}$ be irreducible finite dimensional
representations of the Lie algebra $\g l(n)$. Compute the
polynomials $E_n(V_{\ld})$ and $E_n(V_{\ld}\otimes V_{\mu}^*)$.
\end{prb}

By a theorem of Berenstein and Zelevinsky \cite{BZ}, it is known
that

i) $E_n(V_{\ld})\vert_{q=1}=2^n{\rm
Mult}[V_{\ld}:V_{\delta}\otimes V_{\delta}]$;

ii) ("Kostant conjecture")

$\bullet$ $E_n(V_{\ld})\ne 0$, if and only if $2\delta\ge\ld$ with respect to
the dominance ordering on the set of partitions;

$\bullet$ $E_n(V_{\ld})\vert_{q=1}=2^n$ if and only if there exists a
permutation $w\in S_n$ such that $\ld =\delta +w(\delta )$. Recall
that $\delta :=\delta_n=(n-1,n-2,\ldots ,2,1,0)$.

\begin{pr}{\rm (\cite{STW}, \S 4)} Let $\ld$ be partition, $2\delta\ge\ld$.
Consider partition $\wt\ld =\ld +(1^n)$, then
\begin{equation}
E_n(V_{\ld})=(1+q)\left\{\sum_{\mu
=(a|b)}g_{n^2,\wt\ld ,\mu}q^{b}\right\}, \label{6.10}
\end{equation}
summed over all hook partitions $\mu =(a+1,1^b)$, $a+b=n^2-1$.
\end{pr}

{\it Proof.} Consider the $q$--discriminant of order $n$ by $n$:
$$\Delta_q(X_n)=\prod_{1\le i,j\le n}(x_i-qx_j).
$$
It is a symmetric function in $X_n=(x_1,\ldots ,x_n)$, and we can
consider the decomposition of $q$--discriminant $\Delta_q(X_n)$
in terms of Schur functions $s_{\ld}(X_n)$:
$$\Delta_q(X_n)=\sum_{\mu\vdash
n^2}er_{\mu}(q)s_{\mu}(X_n).
$$
It is not difficult to see that
$$E_n(V_{\ld})=er_{\wt\ld}(-q),$$
where $\wt\ld =\ld +(1^n)$. On the other hand,
\begin{eqnarray*}
\Delta_q(X_n)&=&s_{n^2}((1-q)X_n)=\sum_{\wt\ld\vdash
n^2}(s_{n^2}*s_{\wt\ld})(1-q)s_{\wt\ld}(X_n)\\
&=&\sum_{\wt\ld\vdash n^2} \left(\sum_{\mu\vdash n^2}g_{n^2,\wt\ld
,\mu}s_{\mu}(1-q)\right) s_{\wt\ld}(X_n),
\end{eqnarray*}
so that
$$er_{\wt\ld}(q)=\sum_{\mu\vdash n^2}g_{n^2,\wt\ld
,\mu}s_{\mu}(1-q).
$$
It remains to use the well--known fact, see e.g. \cite{Ma}, that

$$s_{\mu}(1-q)=\left\{\begin{array}{ll}(1-q)(-q)^b,& {\rm
if}~~\mu =(a+1,1^b),\\ 0,& {\rm otherwise}.\end{array}\right.
$$
\qed

Thus, the problem of computation of  polynomials $E_n(V_{\al})$ is
equivalent to finding the structure constants
$g_{\al\beta\gamma}$ of the multiplication of the characters of
the symmetric group $S_{n^2}$ in a {\it particular} case when $\al
=(n^2)$, $\beta\vdash n^2$ and $\gamma =(a|b)$ is a hook
partition, $a+b=n^2-1$.

\begin{pr}\label{p6.11} Let $\al$ and $\beta$ be partitions of the same
size $n^2-n$ whose length do not exceed $n$. Then
\begin{equation}E_n(V_{\al}\otimes
V_{\beta}^*)=\sum_{\ld\vdash n^2}K_{\ld\wt\al}^{(-1)}(-q)K_{\ld\wt\beta}
^{(-1)}(-q)\prod_{i\ge 0} (-q;-q)_{\ld_i'-\ld_{i+1}'}, \label{6.14a}
\end{equation}
where $K_{\ld\mu}^{(-1)}(q)$ denotes the so--called {\it inverse}
Kostka--Foulkes polynomial, i.e. the $(\ld ,\mu)$--entry of the inverse
Kostka--Foulkes matrix $(K_{\ld ,\mu}(q))^{-1}$; in the product on the
RHS(\ref{6.14a}) we made the convention that $\ld_0'=n$; the sum
(\ref{6.14a}) is taken over partitions $\ld$ of size $n^2$ such that
$\ld\le (2n+1,2n-1,\ldots ,3,1)$ with respect to the dominance ordering
on the set of partitions. Note that $K_{\ld\mu}^{(-1)}(q)=0$, unless
$\ld\ge\mu$.
\end{pr}
Formula (\ref{6.14a}) is a special case of a result obtained by
M.~Reeder \cite{Re}. Follow \cite{Re}, we give below a proof of this
formula which is based only on the theory of symmetric functions.

{\it Proof of Proposition~\ref{p6.11}.} Consider the ring $\Lambda_n(t)$
of symmetric polynomials in $n$ independent variables $x_1,\ldots ,x_n$
with coefficients in the field $\Q (q)$ of rational functions in $q$.
There exist two natural scalar products $(~,~)$ and $\langle ~,~\rangle$
on the ring $\Lambda_n(t)$ with values in the field $\Q (t)$. The first
scalar product $(~,~)$ is characterized uniquely by the condition that
the set of all Schur functions
$$\{ s_{\ld}(X_n),~~ \ld ~~{\rm runs~through~all~
partitions~of~length}~~ \le n\}
$$
forms an {\it orthonormal} basis of $\Lambda_n(t)$ with respect
to the scalar product $(~,~)$. Similarly, the second scalar
product $\langle ~,~\rangle$ is defined uniquely by the condition
that the set of all Hall--Littlewood polynomials
$$\{P_{\ld}(X_n;q), ~\ld ~
{\rm runs~through~all~partitions~of~length}~ \le n\}
$$
forms an {\it orthogonal} basis in $\Lambda_n(t)$ with respect to
the scalar product
$$\langle ~,~\rangle ~:~ \langle P_{\ld},P_{\mu}\rangle
=\ds\frac{1}{w_{\ld}(q)}\delta_{\ld ,\mu},
$$
where
$$(1-q)^nw_{\ld}(q)=b_{\ld}(q)(q;q)_{n-\ld_1'}.$$

It is well--known, see e.g. \cite{Ma}, Chapter~VI, \S 9, that
these two scalar products are connected by the following
relations:
$$(f,g)=\langle f\Delta_n(q),g\rangle ,~~~\langle f,g\rangle
=(f/\Delta_n(q),g),
$$
where
$$\Delta_n(q)=\ds\prod_{1\le i\neq j\le n}(1-qx_i/x_j).$$

Now we are ready to prove Proposition~\ref{p6.11}. First of all, if
$$s_{\beta}(x)\Delta_n(q)=\sum_{\mu}c_{\mu ,\beta}(q)P_{\mu}(x;q)w_{\mu}(q),
$$
then
$$c_{\mu ,\beta}(q)=\langle s_{\beta}\Delta_n(q),P_{\mu}\rangle
=(s_{\beta},P_{\mu}).
$$
Thus,
\begin{eqnarray*}
E_n(V_{\al}\otimes V_{\beta}^*;q)&=&(s_{\al},s_{\beta}\Delta_n(q))
=\sum_{\mu}c_{\mu ,\beta}(q)w_{\mu}(q)(s_{\al},P_{\mu})\\
&=&\sum_{\mu}(s_{\al},P_{\mu})(s_{\beta},P_{\mu})w_{\mu}(q).
\end{eqnarray*}
It remains to observe that since
$$s_{\al}(x)=\sum_{\eta}K_{\al\eta}(q)P_{\eta}(x;q),$$
the scalar product $(s_{\al},P_{\mu})$ is equal to the inverse
Kostka--Foulkes polynomial $K_{\mu\al}^{(-1)}(q)$.

\qed

\begin{rem} {\rm Similar arguments may be used to prove the statement
$iii)$ of Theorem~\ref{t5.2}. Namely,}
\begin{eqnarray*}G_n(V_{\al}\otimes
V_{\beta}^*)&=&\ds\frac{(q;q)_n}{(1-q)^n}
(s_{\al},s_{\beta}/\Delta_n(q))=[n]!\langle
s_{\al},s_{\beta}\rangle\\
&=&\sum_{\mu}\frac{K_{\al\mu}(q)K_{\beta\mu}(q)}{b_{\mu}(q)}
\frac{(q;q)_n}{(q;q)_{n-\mu_1'}}.
\end{eqnarray*}
\end{rem}
\begin{rem}\label{r6.12} {\rm Combinatorial interpretation of the
inverse Kostka--Foulkes polynomials\break $K_{\ld\mu}^{(-1)}(q)$
was obtained by Egecioglu and Remmel \cite{ER}. See also an
interesting paper by J.O.~Carbonara in Discrete Math. {\bf 193}
(1998), 117-145.}
\end{rem}

\vskip 0.3cm
\hskip -0.6cm{\bf Exercises} \vskip 0.2cm

\hskip -0.6cm{\bf 1.} Let $\al ,\beta$ be partitions, $|\al|=|\beta|=n$.
Show that

$$L_{\al\beta}^{(1^n)}(q)=K_{\beta'\al}(q,q){\wt K}_{\al ,(1^n)}(q)=
K_{\al'\beta}(q,q){\wt K}_{\beta ,(1^n)}(q),
$$
where
$$K_{\al\beta}(q,q):=K_{\al\beta}(q,t)\vert_{q=t},$$
and $K_{\al\beta}(q,t)$ stands for the double Kostka polynomial
introduced by I.~Macdonald \cite{Ma}, Chapter~VI, (8.11).

\hskip -0.6cm{\bf 2.} (Two variable generalization of polynomials
$L_{\al\beta}^{\mu}(q)$) Let $\al ,\beta ,\mu$ be partitions,
$|\al|=|\beta|=|\mu|=n$, define
$$L_{\al\beta}^{\mu}(q,t)=\sum_{\gamma}g_{\al\beta\gamma}K_{\gamma\mu}(q,t).
$$
Polynomials $L_{\al\beta}^{\mu}(q,t)$ may be considered as a
generalization of the double Kostka polynomials $K_{\al\mu}(q,t)$.
Indeed, if $\beta =(n)$, then
$$L_{\al\beta}^{\mu}(q,t)=K_{\al\mu}(q,t),~~
{\rm and}~~ L_{\al (1^n)}^{\mu}(q,t)=K_{\al'\mu}(q,t).
$$
Polynomials $L_{\al\beta}^{\mu}(q,t)$ have properties similar to
those of $K_{\al\mu}(q,t)$.

Show that

i) $L_{\al\beta}^{\mu}(0,t)=L_{\al\beta}^{\mu}(t)$;

ii) $L_{\al\beta}^{\mu}(0,0)=g_{\al\beta\mu}$,
$L_{\al\beta}^{\mu}(1,1)=f^{\al}f^{\beta}$, where $f^{\al}$ denotes the
number of standard (i.e. weight $(1^{|\al|})$) Young tableaux of shape
$\al$;

iii) $L_{\al\beta}^{\mu}(q,t)=L_{\al'\beta'}^{\mu'}(t,q)$;

iv) $L_{\al\beta}^{\mu}(q,t)=q^{n(\mu')}t^{n(\mu)}
L_{\al'\beta}^{\mu}(q^{-1},t^{-1})$;

v) $L_{\al\beta}^{1^n}(q,t)=K_{\al'\beta}(t,t){\wt K}_{\beta ,(1^n)}(t)
=K_{\beta'\al}(t,t){\wt K}_{\al ,(1^n)}(t)$.

vi) let $\al,\beta$ and $\mu$ be partitions of the same size, then
$$\sum_{\eta}K_{\mu\eta}(q)L_{\al\beta}^{\eta}(q)/b_{\eta}(q)=
L_{\al\beta}^{\mu}(q,q)/H_{\mu}(q),
$$
where $H_{\mu}(q)$ denotes the hook polynomial
$\ds\prod_{x\in\mu}(1-q^{h(x)})$, and for any partition $\eta$ we
put $b_{\eta}(q)=\ds\prod_{i\ge 1}(q;q)_{\eta_i'-\eta_{i+1}'}$,
c.f. Corollary~\ref{c5.6}.

{\hskip -0.6cm{\bf 3.} Let $\ld$ and $\mu$ be partitions and
$\ld\ge\mu$ with respect to the dominance ordering on the set of
partitions. Denote by $d(\al\beta\mu)$ the degree of polynomial
$L_{\al\beta}^{\mu}(q)$, and consider polynomial
$$\wt L_{\al\beta}^{\mu}(q)=q^{d(\al\beta\mu)}
L_{\al\beta}^{\mu}(q^{-1}).
$$

Show that
$$\wt L_{\al\beta}^{\mu}(q)\ge \wt L^{\ld}_{\al\beta}(q),$$
i.e. the difference $\wt L_{\al\beta}^{\mu}(q)-\wt
L_{\al\beta}^{\ld}(q)$ is a polynomial with non--negative
coefficients. More generally, let $I:=[\mu_1,\mu_2]$ be an
interval in the Young graph $(\Y ,\le )$, then
$$\sum_{\tau\in [\mu_1,\mu_2]}\mu(\tau ;I)\wt L_{\al\beta}^{\tau}(q)\ge
0,
$$
where $\mu (\tau ;I)$ denotes the M\"obius function of the interval $I$.

If partition $\alpha$ (or $\beta$) consists of only one part, the
latter inequality was discovered and proved by Lascoux and
Sch\"utzenberger \cite{LS2}; see also \cite{Kir1} for yet another
proof.

$\bullet$ Construct a natural embedding of the sets
$$\bnu^{\mu}(\al ,\beta)\hookrightarrow\bnu^{\ld}(\al ,\beta).
$$
Hence, for any partition $\mu$ there exists a natural embedding
(standardization map)
\begin{equation}
i_{\mu}~:~\bnu^{\mu}(\al ,\beta)\hookrightarrow STY(\al)\times STY(\beta).
\label{6.14}
\end{equation}

\begin{prb}\label{pb6.12} Describe the image of the standardization map
(\ref{6.14}).
\end{prb}
\begin{prb}\label{pb6.13} Let $\al$ and $\beta$ be partitions,
define coefficients $r_{\al\beta k}^{(N)}$ via the decomposition
$$s_{\al}*s_{\beta}(q,\ldots ,q^{N-1})=\sum_{k=0}^{N|\al|}r_{\al\beta
k}^{(N)}q^k.
$$
It follows from Corollary~\ref{c6.8} that if $1\le
k\le\ds\frac{1}{2}N|\al |$, then
$$g_{\al\beta k}^{(N)}:=r_{\al\beta , k}^{(N)}-r_{\al\beta ,k-1}^{(N)}\ge
0.
$$
Give a combinatorial, algebra--geometrical and representation
theoretical interpretations of the non--negative integer numbers
$g_{\al\beta k}^{(N)}$.
\end{prb}

\hskip -0.6cm{\bf 4.} Let $N\ge 1$, $M\ge 1$ be integer numbers.
For any pair of partitions $\al$ and $\beta$ of the same size
$n$, consider the following polynomial (cf. (\ref{6.9a})):
\begin{eqnarray}
&&S_{\al\beta ;N,M}(q)=\label{6.15}\\
&&\frac{1}{n!}\sum_{w\in S_n} \chi^{\al}(w)
\chi^{\beta}(w)\prod_{k\ge
1}\left(\frac{q^k-(1+(-1)^k)q^{kN}+(-1)^kq^{k(N+M-1)}}{1-q^k}
\right)^{\rho_k(w)}.\nonumber
\end{eqnarray}
It is clear that
\begin{eqnarray*}
S_{\al\beta ;N,1}(q)&=&s_{\al}*s_{\beta}(q,\ldots ,q^{N-1}),\\
S_{\al\beta ;2,2}(q)&=&(1+q)\sum_{\mu =(a|b)}g_{\al\beta\mu}q^b,
\end{eqnarray*}
summed over all hook partitions $\mu =(a+1,1^b)$, $a+b=n-1$.

$\bullet$ Show that $S_{\al\beta ;N,M}(q)$ is a polynomial with
non--negative integer coefficients.

\begin{ex} {\rm Take partitions $\al =(31)$ and $\beta =(22)$. Using
the character table for the symmetric group $S_4$, one can easily
find the following expression for the internal product of Schur
functions in question:
$$s_{\al}*s_{\beta}=\frac{1}{4!}(6p_1^4-6p_2^2),$$
and therefore,
\begin{eqnarray*}
S_{\al\beta ;n,m}(q)&=&q^5\left[\begin{array}{c}n+m-2\\
1\end{array}\right]_q\left[\begin{array}{c}2n+2m-5\\
1\end{array}\right]_q+q^7\left[\begin{array}{c}3\\
1\end{array}\right]_q\left[\begin{array}{c}n+m-1\\
4\end{array}\right]_q\\
&+&q^9\left[\begin{array}{c}3\\
1\end{array}\right]_q\left[\begin{array}{c}n+m-2\\
4\end{array}\right]_q+q^{2n+2}\left[\begin{array}{c}n-1\\
1\end{array}\right]_{q^2}\left[\begin{array}{c}m-1\\
1\end{array}\right]_{q^2}.
\end{eqnarray*}
In particular,
$$S_{\al\beta ;2,2}(q)=q^5(1+q+q^2)+q^6=q^5(1+q)^2.$$}
\end{ex}

\begin{prb} Find a fermionic formula for polynomials $S_{\al\beta
;N,M}(q)$ which generalizes that (\ref{6.9})
\end{prb}

\hskip -0.6cm{\bf 5.} Let $\al$ be partition of size $n\ge 2$.
Consider the new partition
$$\ld_{\al}=\al +((n-2)^n)=(\al_1+n-2,\al_2+n-2,\ldots ,\al_n+n-2).
$$
Show that
\begin{equation}
E_n(V_{\ld_{\al}})=\prod_{(i,j)\in \al}(1+q^{2(i-j)+1})K_{\al ,(1^n)}(q^2).
\label{6.16}
\end{equation}
In particular, if $\al =(1^n)$, then
$$E_n(\ld_{\al})=\ds\prod_{j=1}^n(1+q^{2j-1})
$$
is the Poincare polynomial of the unitary group $U(n)$; if
$\al =(2,1^{n-2})$, i.e. the highest weight of the adjoint
representation, then
$$E_n(V_{\ld_{\al}})=\prod_{j=1}^{n-1}(1+q^{2j-1})\frac{1-q^{2(n-1)}}{1-q}.
$$
Formula (\ref{6.16}) was stated as a conjecture by Gupta and
Hanlon \cite{GH} and has been proved by Stembridge \cite{Stm2}.

\hskip -0.6cm{\bf 6.} Assume that partition $\ld$ has the form
$\ld =\delta_n+w(\delta_n)$ for some permutation $w\in S_n$. Show
that
\begin{equation}
E_n(V_{\ld})=q^{l(ww_0})\prod_{k\ge 1}(1+q^{2k-1})^{a_k(w)-a_{k-1}(w)},
\label{6.17}
\end{equation}
where $a_k(w)$ is equal to the length of the longest subsequence
of $w$ that can be written as union of $k$ increasing
subsequences.

\hskip -0.6cm{\bf 7.} (Explicit description of partitions $\ld$
such that $E_n(V_{\ld})|_{q=1}=2^n$)

By a theorem of B.~Kostant (see \cite{BZ}, Proposition~13) and
that of A.~Be\-ren\-stein and A.~Zelevinsky \cite{BZ}, Theorem~14,
the partitions $\ld$ such that $E_n(V_{\ld})|_{q=1}=2^n$ are in a
bijective correspondence with subsets $I\subset [1,n-1]$:
partition $\ld (I)$ corresponding to $I$ is equal to
$\delta_n+w_I(\delta_n)$, where $w_I$ is the element of the
maximal length in the Weyl group generated by the simple
reflections $s_i$, $i\in I$. The main purpose of the Exercises~6
and 7 is to give a direct construction of the  partition
$\ld(I)$, and describe explicitly the RHS(\ref{6.17}) in terms of
the subset $I$.

Let $I$ be a subset of the interval $[1,n-1]$. For each $j\in
[1,n]\setminus I$ we define two numbers $k_I^+(j)$ and $k_I^-(j)$:

$\bullet$ Let us put $k_I^{+}(n-1)=k_I^+(n)=0$, and if $1\le j\le n-2$,
then
$$k_I^+(j)=\max\{ k|[j+1,\ldots ,j+k]\subset I\}.
$$
Thus, we have $k_I^+(j)=0$, if $j+1\not\in I$.

$\bullet$ Let us put $k_I^-(1)=0$, and if $2\le j\le n$, then
$$k_I^-(j)=\max\{ k|[j-k,\ldots ,j-1]\subset I\}.$$
Thus, we have
$k_I^-(j)=0$, if $j-1\not\in I$.

Finally, let us define partition $\ld (I)=(\ld_1(I),\ldots ,\ld_n(I))$
by the rule
$$\ld (I)_i=\sum_{j\ge i}(2+k_I^-(j)+k_I^+(j))\chi (j\not\in I).
$$

$\bullet$ Show that \vskip -0.5cm
$$E(V_{\ld})|_{q=1}=2^n
$$
if and only if there exists a subset $I\subset [1,n-1]$, such that
$\ld =\ld (I)$. Thus, the number of partitions $\ld$ such that
$E(V_{\ld})|_{q=1}=2^n$ is  equal to $2^{n-1}$.

\hskip -0.6cm {\bf 8.} Let $I=\{ 1\le j_1<j_2\le\cdots <j_p\le
n-1\}$ be a subset of the interval $[1,n-1]$. Define by induction
the composition $m=(m_1,\ldots ,m_s)$ of the number $p=\# |I|$.
First of all, put $m_0=1$. If $i\ge 1$ and the numbers
$m_0,\ldots ,m_{i-1}$ have been already defined, then put $m_i=$
$$\max\{ k\in [1,p+1-m_0-\cdots -m_{i-1}]~|~j_{m_0+\cdots
+m_{i-1}+k-1}=j_{m_0+\cdots +m_{i-1}}+k-1\} .
$$
Let $s:=\max\{ j|m_j\ne 0\}$, it is clear that $m_1+\cdots +m_s=p$.

Show that
$$E_n(V_{\ld (I)})=q^{c(I)}(1+q)^{n-p}\prod_{j=1}^s\prod_{r=1}^{m_j}
(1+q^{2r+1}),
$$
where $c(I)=\left(\begin{array}{c}n\\ 2\end{array}\right)
-\sum_{j=1}^s\left(\begin{array}{c}m_j+1\\ 2\end{array}\right)$.

In particular, if $I=\{1,2,\ldots ,\wh j,\ldots ,n-1\}$, then
$$E_n(V_{\ld (I)})=q^{j(n-j)}(1+q)^2\prod_{r=1}^{j-1}(1+q^{2r+1})
\prod_{r=1}^{n-j-1}(1+q^{2r+1}).
$$

\hskip -0.6cm{\bf 9.} Let $\ld =(\ld_1,\ld_2,\ldots ,\ld_n)$ be a
partition. Show that

$\bullet$ (Reduction formulas)\\
$i)$ if $\ld\ge 2\delta_{n+1}$, then
$E_{n+1}(V_{\ld})=q^n(1+q)E_n(V_{\wt\ld})$,

where $\wt\ld
=(\ld_1-2,\ldots ,\ld_n-2)$;\\
$ii)$ if $\ld\ge 2\delta_n$, then
$E_{n+1}(V_{\wh\ld})=q^n(1+q)E_n(V_{\ld})$,

where $\wh\ld
=(2n,\ld_1,\ldots ,\ld_n)$.

$\bullet$  (Duality symmetry)\\
$E_n(V_{\ld})=E_n(V^*_{\ld})=E_n(V_{\ld^+})$,

where $\ld^+=(2n-2-\ld_n,2n-2-\ld_{n-1},\ldots ,2n-2-\ld_1)$.\\
In the case $q=1$ these formulae may be found in \cite{STW}, \S\S
5,6.

$\bullet$ (Poincare duality)\\ $E_n(V;q)=q^{n(n-1)}E_n(V;q^{-1})$.
\vskip 0.3cm

\hskip -0.6cm{\bf 10.} Show that
$$\sum_{\ld\le 2\delta_n}E_n(V_{\ld};-q)G_n(V_{\ld};q)=1,
$$
where $G_n(V;q)$ is the generalized exponents polynomial, see
(\ref{5.1}).

\hskip -0.6cm{\bf 11.} Let $\ld
=(\ld_1\ge\ld_2\ge\cdots\ge\ld_k>0)$ be a partition of $n$. Show
that
\begin{equation}
K_{\ld ,(1^n)}^{(-1)}(q)=(-1)^{|\ld|-k}q^{d(\ld)}\left[\begin{array}{c}
\ld_1'\\\ld_1'-\ld_2',\ldots
,\ld_{k-1}'-\ld_k',\ld_k'\end{array}\right]_q, \label{6.19}
\end{equation}
where $d(\ld)=n(\ld)+|\ld|+k(k+1)/2$.

Formula (\ref{6.19}) was conjectured by J.O.~Carbonara and proved
by I.G.~Macdonald, see Discr. Math. {\bf 193} (1998), 117-145,
Proposition~3.

\hskip -0.6cm{\bf 12.} Let $\al$ and $\beta$ be partitions of the
same size $n$. Show that
$$s_{\al}*s_{\beta}(1-q)=\sum_{\mu}K_{\mu\al}^{(-1)}(q)K_{\mu\beta}^{(-1)}(q)
b_{\mu}(q)=H_{\beta}(q)K_{\beta\al}^{(-1)}(q,q),
$$
where $K_{\beta\al}^{(-1)}(q,q)$ denotes the $(\beta ,\al
)$--entry of the inverse matrix $\left((K_{\ld\mu}(q,q))_{\ld
,\mu\vdash n} \right)^{-1}$ of size $p(n)$ by $p(n)$. Here $p(n)$
stands for the number of partitions of size $n$. Recall, see
\cite{Ma}, Chapter~VI, Example~6, p.364, that
$$\det |(K_{\ld\mu}(q,q))_{\ld ,\mu\vdash n}|=\prod_{\ld\vdash
n}\frac{H_{\ld}(q)}{b_{\ld}(q)}.
$$

\hskip -0.6cm{\bf 13.} (Generalization of Corollary~\ref{c6.8})
Let $\al$ and $\beta$ be partitions of the same size and
$(k_1,k_2,\ldots ,k_r)$ be a unimodal and symmetric sequence  of
non--negative integer numbers.

Show that the specialization
\begin{equation}
s_{\al}*s_{\beta}(\underbrace{q,\ldots ,q}_{k_1},
\underbrace{q^2,\ldots ,q^2}_{k_2},\ldots ,\underbrace{q^r,\ldots
,q^r}_{k_r}) \label{6.21*}
\end{equation}
is a symmetric and unimodal polynomial.

This result is due essentially to A.~Kerber, see e.g., R.~Stanley
(Studies in App. Math. {\bf 72} (1985), 263--281), or \cite{St2},
Theorem~14.

\hskip -0.6cm{\bf Question.} Does there exist a fermionic formula
for polynomials (\ref{6.21*}), which generalizes that (\ref{6.9})
?

\hskip -0.6cm{\bf 14.} Let $\ld$ be a partition of size $n$.
Consider partition $\wt\ld =(3n,2n,\ld )$ and a sequence of
rectangular shape partitions $R:=((n,n),(n,n),(n),(n))$. Let
$s_{\ld}(1,q,q,q^2)$ denotes the specialization $x_1=1,~~
x_2=x_3=q,~~ x_4=q^2,$ $x_k=0,$ if $k\ge 5$, of the Schur function
$s_{\ld}(x)$. Show that
$$s_{\ld}(1,q,q,q^2)\bdoteq K_{\wt\ld R}(q).
$$

\section{Liskova semigroup}
\label{kfs}
\neweq

\subsection{Realizable polynomials and product formula}
\label{rp}

\begin{de} A polynomial $P(q)\in\N[q]$ with non--zero constant term
is called $K$--realizable if there exist two partitions $\ld$ and
$\mu$, $|\ld |=|\mu|$, a composition $\eta$, $|\eta|=l(\mu)$, and
an integer $L$ such that
$$P(q)=q^LK_{\ld\mu\eta}(q).
$$
\end{de}

The set of all $K$--realizable polynomials will be denoted by
${\L}$.

\begin{de} A polynomial $P(q)\in\N [q]$ with non--zero constant term
is called PK--rea\-li\-zable if there exist a partition $\ld$, a
dominant sequence of rectangular shape partitions $R$, and an
integer $L$ such that
$$P(q)=q^LK_{\ld R}(q).
$$
\end{de}

The set of all PK--rea\-lizable polynomials will be denoted by
${\L}^+$.

\begin{de} A polynomial $P(q)\in 1+q\N [q]$ is called KF--realizable, if
there exist partitions $\ld$ and $\mu$ such that
$$P(q)={\wt K}_{\ld\mu}(q).
$$
\end{de}

The set of all KF--realizable polynomials will be denoted by
${\L^{++}}$. It is clear that ${\L}^{++}\subset{\L}^+\subset{\L}$.

\begin{prb} {\rm Describe the sets ${\L}$, ${\L}^+$ and ${\L}^{++}$.}
\end{prb}

\begin{con} If $n\ge 3$ and $a\in\Z_{\ge 1}$, then $1+aq^n\not\in{\L}$.
\end{con}

\begin{con} Let $K_{\ld R}(q)=q^{\bullet}(a_0,a_1,\ldots ,a_l)$ be
parabolic Kostka polynomial corresponding to a partition $\ld$
and a dominant sequence of rectangular shape partitions $R$. Then
both sequences $(a_0,a_2,a_4,\ldots )$ and $(a_1,a_3,a_5,\ldots )$
are unimodal.
\end{con}

What one can say about the structure of the sets ${\L}$, ${\L}^+$,
and ${\L}^{++}$? First of all, it follows from Duality Theorem
that if a polynomial $P(q)$ belongs to the set ${\L}^+$, then the
polynomial $q^{{\rm deg}(P)}P(q^{-1})$ also belongs to the set
${\L}^+$. Another interesting property of the set ${\L}^+$ follows
from the following result.

\begin{theorem}\label{t7.6} Let $\al$ and $\beta$ be partitions, and $R$,
$S=\{(\mu_a^{\eta_a})\}$, $1\le a\le p$, be two dominant sequences of
rectangular shape partitions. Let $k=\sum\mu_a$, and $n$ be such that
$n\ge l(\al )+l(\beta )$, and $n\ge\eta_a$, $1\le a\le p$. Consider
partition
$$[\al ,\beta ]_{n,k}:=(\al_1+k,\al_2+k,\ldots
,\al_s+k,\underbrace{k,\ldots k}_{n-s-r},k-\beta_r,\ldots ,k-\beta_1),
$$
and a dominant rearrangement $Q$ of the sequence of rectangular
shape partitions $R\cup \{(\mu_a^{n-\eta_a})\}^p_{a=1}$. Then
$$K_{[\al ,\beta ]_{n,k}Q}(q)=K_{\al R}(q)K_{\beta S}(q).
$$
\end{theorem}

It follows from Theorem~\ref{t7.6} that the set of PK--realizable
polynomials ${\L}^+$ is closed under multiplication, and forms a
semigroup called by {\it Liskova} semigroup.

\hskip -0.6cm{\bf Question.} Denote by $\K_n^+$ the set of all
degree $n$ parabolic Kostka polynomials $K_{\ld R}(q)$
corresponding to a partition $\ld$ and a dominant sequence of
rectangular shape partitions $R$. Is it true that the set
${\K}_n^+$ is {\it finite} for all $n\ge 0$ ? If so, what does the
generating function $\sum_{n\ge 0}|{\K}_n^+|t^n$ looks like?

\begin{con}\label{con7.8} Let $P(q)$ be a polynomial with
non--negative integer coefficients and non--zero constant term.
There exists a non--negative integer $N$, depending on $P(q)$,
such that the product $(1+q)^NP(q)$ belongs to the Liskova
semigroup $\L^+$.
\end{con}

The least non--negative integer $N$ such that
$(1+q)^NP(q)\in\L^+$ will be denoted by $IL(P)$. For example, one
can show that if $n$ and $m$ are non--negative integers then
$IL(n)\le 2$, see e.g., Section~\ref{ffpkp}, Exercise~6, and
$IL(n+mq)\le |n-m|+\delta_{n,m}$.

\begin{prb} Given an integer $n\ge 3$, compute the number
$IL(1+q^n)$.
\end{prb}

\subsection{ Generalized $q$--Gaussian coefficients}
\label{ggc}

Let $\ld$ be a partition, and
$$\left[\begin{array}{c}n\\ \ld\end{array}\right]:=\prod_{x\in\ld}
\frac{1-q^{c(x)}}{1-q^{h(x)}}
$$
denotes the generalized $q$--Gaussian coefficient, \cite{Ma}, Chapter~I,
\S 3, Example~1. Then the polynomial $\left[\begin{array}{c}n\\
\ld\end{array}\right]$ is KF--realizable. This statement follows from
\cite{Kir3}, Lemma~1. Namely, let $\ld$ be a partition, consider
partitions
$$\ld_N=(N|\ld|,\ld ):=(N|\ld|,\ld_1,\ld_2,\ldots
),~~{\rm and}~~ \mu_N=(|\ld|^{N+1})=(\underbrace{|\ld|,\ldots
,|\ld|}_{N+1}).
$$
Then \vskip -0.5cm
\begin{equation}
{\wt K}_{\ld_N\mu_N}(q)=\left[\begin{array}{c}N\\ \ld'\end{array}\right].
\label{7.1}
\end{equation}
See \cite{Kir3,Kir6} for proof and further details. It has been shown in
\cite{Kir3,Kir6} that the identity (\ref{7.1}) can be considered as a
generalization of the KOH identity (\ref{4.4}), see \cite{OH,Z,Kir3}.

Identity (\ref{7.1}) admits a generalization. To start, let us
recall the well--known facts (see, e.g. \cite{Ma}, Chapter~I) that
\begin{equation}
s_{\ld}(1,\ldots, q^{N-1})=q^{n(\ld )}\left[\begin{array}{c}N\\
\ld'\end{array}\right], ~~~s_{\ld}*s_{(n)}(1,\ldots ,q^{N-1})=
s_{\ld}(1,\ldots ,q^{N-1})\label{7.2}
\end{equation}
and hence
\begin{equation}
s_{\ld}*s_{(n)}(1,\ldots ,q^{N-1})\bdoteq
K_{\ld_N\mu_N}(q), \label{7.3}
\end{equation}
where $\ld_N=(N|\ld|,\ld )$, and $\mu_N=(|\ld|^{N+1})$. We will give
below a generalization of (\ref{7.3}) when partition $(n)$ is replaced by
an arbitrary rectangular shape partition.

\begin{theorem}\label{t7.5} Let $\ld$ be a partition, and
$R=(l^m)$ be a
rectangular shape partition, $|\ld|=lm$. Consider the partition
$\ld_N=(\underbrace{Nl,\ldots Nl}_m,\ld )$ and the sequence $Q_N$ of
rectangular shape partitions $Q_N=(\underbrace{R,\ldots ,R}_{N+1})$. Then
\begin{equation}
s_{\ld}*s_R(q,\ldots ,q^{N-1})\bdoteq K_{\ld_NQ_N}(q). \label{7.4}
\end{equation}
\end{theorem}

\begin{rem}\label{r7.8} {\rm In Section~\ref{ipsf}, Theorem~\ref{t6.1}, we
gave another generalization of the identity (\ref{7.3}) when partition
$(n)$ is replaced by {\it arbitrary} partition $\beta$, $|\beta|=n$. As
a corollary, we obtain the following identity
$$K_{[\al ,R]_{Nr},R_N}(q)\bdoteq K_{\al_N,Q_N}(q),
$$
where \begin{eqnarray*} [\al ,R]_{Nr}&=&(\al_1+l,\ldots
,\al_r+l,\underbrace{l,\ldots ,l}_{Nr-r-m}), \\
\al_N&=&(\underbrace{Nl,\ldots ,Nl}_m,
\al_1,\ldots ,\al_r),\\
R_N&=&\{ (l^r)^N\},~~Q_N=\{ (l^m)^{N+1}\},
\end{eqnarray*}
and $r=l(\al)$.}
\end{rem}

\subsection{Multinomial coefficients}
\label{mc}

\begin{lem}\label{l7.7} {\rm (\cite{Ma}, Chapter~I)} Let $\ld =(n,1^m)$ be a hook,
and $\mu$ be a partition. Then
$$\wt K_{\ld\mu}(q)=\left[\begin{array}{c}\mu_1'-1\\ m\end{array}\right]_q.
$$
\end{lem}

It follows from Lemma~\ref{l7.7} and Theorem~\ref{t7.6} that the
$q$--multinomial coefficient
$$\left[\begin{array}{c} N\\ m_1,\ldots ,m_k\end{array}\right]_q:=
\frac{(q;q)_N}{(q;q)_{m_1}\cdots (q;q)_{m_k}},
$$
where $m_i$, $1\le i\le k$, are non--negative integers and
$N=m_1+\cdots +m_k$, belongs to the Liskova semigroup ${\L}^+$.
Indeed,
$$\left[\begin{array}{c} N\\ m_1,\ldots ,m_k\end{array}\right]_q=
 \left[\begin{array}{c} N\\ m_1\end{array}\right]_q
 \left[\begin{array}{c} N-m_1\\ m_2\end{array}\right]_q\cdots
 \left[\begin{array}{c} N-m_1-\cdots -m_{k-1}\\ m_k\end{array}\right]_q.
$$

\vskip 0.2cm \hskip -0.6cm {\bf Question.} Is it true that the
$q$--multinomial coefficients belong to the set ${\L}^{++}$?
\vskip 0.2cm

Note, that according to Section~\ref{pkp}, (\ref{3.23}), the
$q^2$--binomial coefficient $\left[\begin{array}{c}n\\
m\end{array}\right]_{q^2}$ also belongs to the Liskova semigroup
$\L^+$. Thus, the $q^2$--multinomial coefficients
$\left[\begin{array}{c}N\\ m_1,\ldots
,m_k\end{array}\right]_{q^2}$ belong to the Liskova semigroup as
well.

\vskip 0.2cm \hskip -0.6cm {\bf Questions.} i) Is it true that
for any partition $\ld$ the generalized $q^2$--binomial
coefficient $\left[\begin{array}{c}n\\ \ld
\end{array}\right]_{q^2}$ belongs to the Liskova semigroup ?

ii) Assume that some symmetric and unimodal polynomial $P(q)$
belongs to the Liskova semigroup $\L^+$. Is it true that the
polynomial $P(q^2)$ also belongs to the semigroup $\L^+$ ?

\subsection{Kostka--Foulkes polynomials and transportation matrices}
\label{kfptm}

Let $\ld ,\mu$ be partitions, $|\ld|=|\mu|$, and $n,N$ be natural
numbers such that $l(\ld)\le n$, $l(\mu)\le n$, and $N\ge
\ld_1+\mu_1$. Let us define partitions $\al_N=(N^n)$ and
$$\beta_N=(N-\ld_n,N-\ld_{n-1},\ldots ,N-\ld_1,\mu_1,\mu_2,\ldots ,\mu_n).
$$

 \begin{theorem}\label{t7.8} Let $\ld ,\mu ,n,N,\al_N$ and $\beta_N$ be as
 above. Then

i)~ $K_{\al_N\beta_N}(q)\le K_{\al_{N+1}\beta_{N+1}}(q)$;
\begin{equation}
ii)\hskip 0.27cm
If~~N\ge|\ld|,~~then~~K_{\al_N\beta_N}(q)\bdoteq\sum_{\eta}K_{\eta\ld}(q)
K_{\eta\mu}(q).~~~~~~~~~\label {7.5}
 \end{equation}
 \end{theorem}

 Identity (\ref{7.5}) admits the following combinatorial interpretation.
 First of all, note that the value of the LHS(\ref{7.5}) at $q=1$ is
 equal to the number of plane partitions $\pi$ of shape $(N^n)$ whose
 diagonal sums are fixed. More precisely, $\pi =(\pi_{ij})_{1\le i<j\le
 n}$ is a $n$ by $n$ matrix of non--negative integers such that
 \begin{eqnarray*}
&&\pi_{ij}\ge\pi_{i+1j},~~~1\le i\le n-1,~~~1\le j\le n;~~~
\pi_{ij}\ge\pi_{ij+1},\\
&&1\le i\le n,~~~1\le j\le n-1;~~~
\sum_{j=1}^k\pi_{j,k-j+1}=\sum_{j=1}^k\mu_{n-k+j};\\
&&\sum_{j=1}^{n-k+1}\pi_{n-j+1,k+j-1}=\sum_{j=k}^n\ld_j,~~~1\le
k\le n.
 \end{eqnarray*}
 Indeed, let $T$ be a semistandard Young tableau of shape $(N^n)$ and
 weight $(N-\ld_n,N-\ld_{n-1},\ldots ,N-\ld_1,\mu_1,\ldots ,\mu_n)$.
 Denote by $d_{ij}$, $1\le i\le n$, $1\le j\le 2n$, the number of cells
 in the $i$--th row of $T$ occupied by numbers that do not exceed $j$. Consider
 the diamond $\wt\pi =(\wt\pi_{ij})$, where
 $$\wt\pi_{ij}=N-d_{ij},~~ 1\le
 i\le n,~~ 1\le j\le 2n-1,~~{\rm and}~~ 0\le j-i\le 2n-1.
 $$
 It is clear that
 $\wt\pi_{ij}\ge\wt\pi_{i+1,j}$, $\wt\pi_{ij}\ge\wt\pi_{i,j+1}$ and
 $$\sum_{j=1}^k\wt\pi_{jk}=\sum_{j=1}^k\mu_{n-k+j},~~~
 \sum_{j=k}^n\wt\pi_{j,n+k-1}=\sum_{j=k}^n\ld_j,~~~1\le k\le n.
 $$
 If we put $\pi_{ij}=\wt\pi_{i,j-i+1}$,
 we will obtain the plane partition
 $\pi =(\pi_{ij})_{1\le i,j\le n}$ as it was defined above.

 Clearly, the diamond $\wt\pi$ defines two semistandard Young tableaux $P$
 and $Q$ of the same shape: $P$ of weight
 ${\overleftarrow\ld}=(\ld_n,\ld_{n-1},\ldots ,\ld_1)$ and $Q$ of weight
 ${\overleftarrow\mu}=(\mu_n,\mu_{n-1},\ldots ,\mu_1)$. The statement ii)
 of Theorem~\ref{t7.10} is equivalent to the following pure combinatorial
 one: the difference of charges
 $$c(T)-c(\Ss (P))-c(\Ss (Q))$$
 is a constant depending only on $n,\ld$ and $\mu$. Here $\Ss$ denotes the
 Sch\"utzenberger transformation
 $$\Ss :STY(\ld ,\mu)\to STY(\ld,{\overleftarrow\mu})
 $$
 on the set of semistandard Young tableaux
 of shape $\ld$. See e.g. \cite{LS2}, or \cite{Ful}, Appendix,
 p.184, where the transformation $\Ss$ is called {\it evacuation},
 or \cite{Kir1,Kir6}.

 \begin{theorem}\label{t7.9} Let $\ld ,\mu$ be partitions, $|\ld|=|\mu|$,
 $l(\ld)\le n$, $N\ge\ld_1$. Let us define the rectangular shape
 partition $\al_N=(n^N)$, and dominant sequence of rectangular shape partitions
 $R_N=\{\mu ,(1^{N-\ld_n}),\ldots ,(1^{N-\ld_1})\}$. Then

i)~ $K_{\al_NR_N}(q)\le K_{\al_{n+1}R_{N+1}}(q)$;
\begin{equation}
ii)~~~  If~~ N\ge |\ld|,~~ then ~~
K_{\al_NR_N}(q)\bdoteq\sum_{\eta}
 {\overline K}_{\eta\ld}(q)K_{\eta'\mu}(q).~~~~~
 \label{7.6}
 \end{equation}
 \end{theorem}

 More generally, let $\beta$ be partition, $l(\beta)\le n$, and
 $R,S=\{(\mu_a^{\eta_a})\}_{a=1}^p$ be two dominant sequences of
 rectangular shape partitions such that $\sum\eta_a=n$, and
 $|S|=|R|+|\beta|$. Consider partitions
 $\al_N=(N-\beta_n,N-\beta_{n-1},\ldots ,N-\beta_1)$, and a dominant
 rearrangement $Q_N$ of the sequence of rectangular shape partitions
 $R\cup\{(N-\mu_a)^{\eta_a}\}^p_{a=1}$.

 \begin{theorem}\label{t7.10} {~~~}

i) ~$K_{\al_NQ_N}(q)\le K_{\al_{N+1}Q_{N+1}}(q)$;

 \begin{equation}
 ii)~~~ If ~~N\ge |S|,~~ then ~~ K_{\al_NQ_N}(q)\bdoteq\sum_{\gamma}
 K_{\gamma\setminus\beta ,R}(q)K_{\gamma ,S}(q).~~\label{7.7}
 \end{equation}
 \end{theorem}

Similarly, let $\beta$ be partition, $l(\beta')\le k$, and
$R,S=\{(\mu_a^{\eta_a})\}^p_{a=1}$ be two dominant sequences of
rectangular shape partitions such that $\sum\mu_a=k$, and
$|S|=|R|+|\beta|$. Consider partition
$\al_N=(N-\beta_k,N-\beta_{k-1},\ldots ,N-\beta_1)'$, and a dominant
rearrangement $Q_N$ of the sequence of rectangular shape partitions
$R\cup\{(\mu_a^{N-\eta_a})\}^p_{a=1}$.

\begin{theorem}\label{t7.11} {~~}
i) ~$K_{\al_NQ_N}\le K_{\al_{N+1}Q_{N+1}}(q)$;
\begin{equation}
ii)~~~ If ~~N\ge |S|,~~ then ~~K_{\al_NQ_N}(q)\bdoteq\sum_{\gamma}
{\overline K}_{\gamma'\setminus\beta',R}(q)K_{\gamma
,S}(q).~\label{7.8}
\end{equation}
\end{theorem}

\subsection{Gelfand--Tsetlin polytope and volume of weight subspaces}
\label{vws}

Let $\ld$ be a partition and $\mu$ be a composition of the same
size and whose lengths do not exceed $n$. Let $GT(\ld ,\mu)$ be
the convex polytope in the space $\R^{\frac{n(n+1)}{2}}$ of all
points ${\bf x}=(x_{ij})_{1\le i\le j\le n}$ satisfying the
following conditions

$\bullet$ $x_{in}=\ld_i$, $1\le i\le n$;

$\bullet$ $x_{i,j+1}\ge x_{ij}$, $x_{ij}\ge x_{i+1,j+1}$, for all $1\le
i\le j\le n-1$;

$\bullet$ $x_{11}=\mu_1$, and if $2\le j\le n$, then
$\sum_{i=1}^jx_{ij}-\sum_{i=1}^{j-1}x_{ij-1}=\mu_j$.

The polytope $GT(\ld ,\mu)$ is called {\it Gelfand--Tsetlin
polytope}. The following result goes back to the original paper by
Gelfand and Tsetlin (Doklady Akad. Nauk SSSR (N.S.) {\bf 71}
(1950), 825-828), and was rediscovered many times.
\begin{pr} Dimension $\dim V_{\ld}(\mu)$ of the weight $\mu$ subspace
$V_{\ld}(\mu)$ of the irreducible highest weight $\ld$
representation $V_{\ld}$ of the general linear algebra $\g l(n)$
is equal to the number of integer points in the Gelfand--Tsetlin
polytope $GT(\ld ,\mu)$:
$$\dim V_{\ld}(\mu)=\#|GT(\ld ,\mu)\cap\Z^{n(n+1)/2}|.
$$
\end{pr}
\begin{con}\label{c7.15} {\rm (\cite{KB})} Let $\ld$ and $\mu$ be partitions.
All vertices of the Gelfand--Tsetlin polytope $GT(\ld ,\mu)$ have
integer coordinates, i.e. $GT(\ld ,\mu)$ is a convex integral polytope.
\end{con}

If Conjecture~\ref{c7.15} is true, then
\begin{equation}
\sum_{l\ge 0}K_{l\ld ,l\mu}t^l=\frac{P_{\ld\mu}(t)}{(1-t)^{d+1}},
\label{7.9b}
\end{equation}
where $d=\dim GT(\ld ,\mu)$, and $P_{\ld\mu}(t)$ is a polynomial
with {\it non--negative} integer coefficients such that
$P_{\ld\mu}(1)$ is equal to the (normalized) volume of the
Gelfand--Tsetlin polytope $G(\ld ,\mu )$.

Another consequence of Conjecture~\ref{c7.15}, which can be
proved using the specialization $q=1$ of {\it fermionic} formula
(\ref{4.1a}), is that the Kostka number $K_{l\ld ,l\mu}$ is a
polynomial ${\cal E}_{\ld\mu}(l)$ in $l$ with integer
coefficients. We will call the latter by {\it Ehrhart polynomial}
of the weight subspace $V_{\ld}(\mu)$, and denote it by ${\cal
E}_{\ld\mu}(t)$. Hence, if $l\in\Z_{\ge 0}$, then
$${\cal E}_{\ld\mu}(t)|_{t=l}=K_{l\ld ,l\mu}.$$
It is well--known, see e.g., Exercise~6 to Section~\ref{kfs}, that
in general the Ehrhart polynomial of a convex integral polytope
may have negative integer coefficients. Nevertheless, based on
examples, we make a conjecture that the Ehrhart polynomial ${\cal
E}_{\ld\mu}(t)$ of the weight subspace $V_{\ld}(\mu )$ has in
fact {\it non--negative} integer coefficients.

\begin{ex}\label{e7.16} {\rm Let $n\ge 2$ be a positive integer. Consider
partitions $\ld =(n^n)$ and $\mu =((n-1)^n,1^n)$. Then the Kostka number
$K_{l\ld ,l\mu}$ is equal to the number of transportation matrices of
size $n$ by $n$ with row and column sums both equal to $(l^n)$. Thus,
${\cal E}_{\ld\mu}(t)$ is the Ehrhart polynomial of the Birkhoff
polytope
\begin{equation}
{\cal B}_n=\left\{(x_{ij})\in\R^{n^2}|\sum_{i=1}^nx_{ij}=1,~
\sum_{j=1}^nx_{ij}=1\right\}, \label{7.8a}
\end{equation}
and $P_{\ld\mu}(1)$ is equal to the normalized volume ${\wt{\rm
vol}}({\cal B}_n)$ of the polytope ${\cal B}_n$.}
\end{ex}

To our knowledge, the Ehrhart polynomial and volume of the Birkhoff
polytope ${\cal B}_n$ are known only up to $n\le 8$, \cite{DLT}. The
(normalized) leading coefficient of Ehrhart's polynomial ${\cal
E}_{\ld\mu}(t)$ of a weight subspace $V_{\ld}(\mu)$ is equal to the
(normalized) volume of Gelfand--Tsetlin's polytope $G(\ld ,\mu)$, and
called by {\it volume} of the weight subspace $V_{\ld}(\mu)$. It plays
an important role in the theory of unitary representations of semisimple
Lie groups (Harish--Chandra, Duistermaat~J.J. and Heckman~G.J.).

\vskip 0.3cm
\hskip -0.6cm {\bf Exercises}\vskip 0.2cm

\hskip -0.6cm{\bf 1.} {\bf a.} Let $\ld$ be a partition, $|\ld|=n$. Show
that
\begin{equation}
\sum_{l(\ld)\le p}\left(K_{\ld
,(1^n)}(q)\right)^2=\sum_{\nu}q^{c(\nu)}\prod_{k,j\ge 1}\left[
\begin{array}{c}P_j^{(k)}(\nu)+m_j(\nu^{(k)})\\ m_j(\nu^{(k)})\end{array}
\right]_q, \label{7.9}
\end{equation}
summed over all admissible configurations $\nu\in C((n^n),
((n-1)^n,1^n))$ of the square shape $(n^n)$ and weight
$((n-1)^n,1^n)$, such that  $l(\nu^{(p)})\le n-p$, i.e. the
number of parts of the partition $\nu^{(p)}$ does not exceed
$n-p$.

For the reader's convenience, we recall the meaning of notions
and symbols which were used in this Exercise (see
Section~\ref{ffpkp}, Definitions~\ref{d4.1} and \ref{d4.2}).

First of all, summation in the RHS(\ref{7.9}) is taken over all
sequences of partitions $\nu =\{\nu^{(1)},\nu^{(2)},\ldots\}$
such that

$\bullet$ $|\nu^{(k)}|=n(n-k)$, $1\le k\le n-1$; by definition,
$\nu^{(0)}=\emptyset$;

$\bullet$ $P_j^{(k)}(\nu ):=n\min
(n,j+1)\delta_{k,1}+Q_j(\nu^{(k-1)})-2Q_j(\nu^{(k)})+Q_j(\nu^{(k+1)})\ge
0$, for all $k,j\ge 1$.

By definition, these two conditions mean that the configuration
$\nu$ is {\it admissible} of type $((n^2);((n-1)^n,1^n))$, i.e.
belongs to the set of admissible configurations
$C((n^2);((n-1)^n,1^n))$. Recall, that notation $((n-1)^n,1^n)$
stands for the partition $(\underbrace{n-1,\ldots
,n-1}_n,\underbrace{1,\ldots ,1}_n)$.

$\bullet$ By definition, for any partition $\nu$,
$Q_j(\nu)=\sum_a\min (j;\nu_a)$ is equal to the number of boxes
in the first $j$ columns of the diagram corresponding to
partition $\nu$, and $m_j(\nu)=\nu_j'-\nu_{j+1}'$ is equal to the
number of parts of partition $\nu$ of size $j$.

Secondly, the symbol $c(\nu )$ stands for the {\it charge} of
configuration $\nu$, i.e.
$$c(\nu)=\sum_{j\ge
1}\left(\begin{array}{c}n(1+\delta_{j,1})-(\nu^{(1)})_j'\\
2\end{array}\right)+\sum_{k\ge 2,j\ge 1}\left(\begin{array}{c}(\nu^{(k-1)})_j'-(\nu^{(k)})_j'\\
2\end{array}\right).
$$

{\bf b.} Deduce from identity (\ref{7.9}) yet another combinatorial
formula for Catalan numbers:
\begin{eqnarray*}C_n&=&\frac{(2n)!}{n!(n+1)!}=K_{(n,n),(1^n)}(1)\\
&=&\sum_{\beta\vdash (n^2-n)}
\prod_{k=1}^{n-1}\left(\begin{array}{c}(2k+1)n
-2k-2(\beta_1+\cdots +\beta_k)+\beta_k-\beta_{k+1}\\
\beta_k-\beta_{k+1}\end{array}\right),
\end{eqnarray*}
summed over all partitions $\beta
=(\beta_1\ge\beta_2\ge\cdots\ge\beta_n)$ of $n^2-n$, such that for all
$k$, $1\le k\le n$, the following inequalities hold
\begin{equation}
k(n-1)\le\beta_1+\cdots +\beta_k\le\frac{1}{2}((2k+1)n-2k-n\delta_{k,n})).
\label{7.9a}
\end{equation}

The sum $\ds\sum_{l(\ld)\le 2}(K_{\ld ,(1^n)}(q))^2$ gives a new
$q$--analog of the Catalan numbers which is different from those of
Carlitz and Riordan~\cite{Car}.

If $p=3$, the RHS(\ref{7.9}) with $q=1$ gives rise to a
combinatorial formula for the number of vexillary permutations
(see, e.g. \cite{Ma1}, Definition~(1.27)) in the symmetric group
$S_n$.




\hskip -0.6cm{\bf 2.} (Kostant partition function)

{\bf a.} Let
$$\Sigma\subset\Phi (1^n)=\{(i,j)\in\Z^2~|~1\le i<j\le n\}$$
be a subset of the set of all positive roots of type $A_{n-1}$.
The $q$--analog of Kostant partition function $K_{\Sigma}(\gamma
;q)$ is defined via the decomposition
\begin{equation}
\prod_{(i,j)\in\Sigma}(1-qx_i/x_j)^{-1}=\sum_{\gamma}
K_{\Sigma}(\gamma ;q)x^{\gamma}, \label{3.13}
\end{equation}
summed over all sequences of integers $\gamma =(\gamma_1,\gamma_2,\ldots
,\gamma_n)\in\Z^n$, $|\gamma|=0$. In other words, $K_{\Sigma}(\gamma ;q)$
is the generating function for the number of ways
of forming $\gamma$ as linear combinations of elements of $\Sigma$ with
non--negative integer coefficients:
$$K_{\Sigma}(\gamma ;q)=\ds\sum_{\{m_{i,j}\}}q^{\parallel m_{ij}
\parallel},
$$
where the sum is taken over all sequences of non--negative integers $\{m_{ij}\}$,
$(i,j)\in\Sigma$, such that
$$\ds\sum_{(i,j)\in\Sigma}m_{ij}e_{ij}=\gamma,$$
where the vector $e_{ij}\in\Z^n$ has the following components
$(e_{ij})_k=\delta_{ik}-\delta_{jk}$, $1\le k\le n$; $\parallel
m\parallel =\ds\sum_{(i,j)\in\Sigma}m_{ij}$.

$i)$ Show that $K_{\Phi(1^n)}(\gamma ;1)\ne 0$, if and only if
$$\gamma\in Y_n=\left\{(y_1,\ldots ,y_n)\in\Z^n\left\vert\begin{array}{l}
y_1+\cdots +y_j\ge 0,~~{\rm if}~~1\le j\le n,\\
y_1+\cdots +y_n=0 \end{array}\right.\right\}.
$$

$ii)$ Show that
\begin{equation}
K_{\ld\mu\eta}(q)=\sum_{w\in S_n}(-1)^{l(w)}K_{\Phi (\eta)}
(w(\ld +\delta )-\mu -\delta ;q). \label{3.14}
\end{equation}

\begin{rem}\label{r3.9} {\rm If $\eta =(1^n)$, then
$$\Phi(\eta)=\{(i,j)~|~1\le i<j\le n\}$$
and $K_{\Phi(1^n)}(\gamma ;q)$ coincides with the $q$--analog of
type $A$ Kostant partition function $P(\gamma ;q)$. In this case
the RHS(\ref{3.14}) takes the form
$$\sum_{w\in S_n}\epsilon (w)P(w(\ld +\delta)-\mu -\delta ;q).
$$
This sum is equal to the Kostka--Foulkes polynomial $K_{\ld\mu}(q)$, see
e.g. \cite{Ma}. Chapter~III, \S 6, Example~4.}
\end{rem}

iii) Let $\gamma\in Y_n$ and $N$ be an integer such that
$N+\gamma_n\ge 0$. Consider partitions $\ld_N=N(n,n-1,\ldots
,2,1)+\gamma$, $\mu_N=N(n,n-1,\ldots ,2,1)$, and composition
$\eta$, $|\eta|=n$. Show that
\begin{equation}
K_{\Phi(\eta)}(\gamma ;q)=K_{\ld_N,\mu_N,\eta}(q). \label{3.15}
\end{equation}

{\bf b.} Let
$$\gamma\in Y_n^+=\{\gamma\in Y_n|\gamma_i\ge 0,~1\le i\le
n-1\},~~ n\ge 3, $$
and $\eta =(1^{n-2},2)$. Show that
\begin{equation}
K_{\Phi (1^n)}(\gamma;1)=\frac{1}{3}(\gamma_1+\cdots
+\gamma_{n-3}+3\gamma_{n-2}+3)K_{\Phi(\eta)}(\gamma ;1). \label{3.16}
\end{equation}
Statement that the LHS(\ref{3.16}) is divisible by the linear
factor
$$l_n(\gamma):=\gamma_1+\gamma_2+\cdots+\gamma_{n-3}+3\gamma_{n-2}+3
$$
was first observed by J.R.~Schmidt and A.M.~Bincer, J. Math. Phys.
{\bf 26} (1984), 2367--2373, see also \cite{Kir6}. In \cite{Kir6}
we obtained also an explicit expression for the ratio
$3K_{\Phi(1^n)}(\gamma ;1)/ l_n(\gamma)$, from which one can see
that the latter can be identified with the parabolic Kostant
partition function $K_{\Phi(1^{n-2},2)}(\gamma ;1)$.

{\bf c.} Consider vector
$$\gamma =\gamma_{d,n}=(d,d+1,\ldots ,d+n-1,-n(2d+n-1)/2)\in Y_{n+1}^+,
$$
and composition $\eta_k =(1^{k},n-k,1)$. Show that
\begin{equation}K_{\Phi (\eta_k)}(\gamma_{d,n};1)=\ds\prod_{p=1}^{k-1}C_p
\!\prod_{1\le i<j\le
k}\frac{2(n-k)+i+j-1}{i+j-1}\!\prod_{(i,j)\in\Omega_{k,n}}
\frac{2d+i+j-1}{i+j-1},\label{3.17}
\end{equation}
where $\Omega_{k,n}=\Phi(1^{n+1})\setminus\Sigma_{k,n}$, and ~~
$\Sigma_{k,n}=$
$$\{ (i,j)\vert\frac{k}{2}<i<j<\frac{2n-k+4}{2}\}
\setminus\{\left(\frac{k+1+2a}{2},\frac{k+3+2a}{2}\right)\vert
0\le a\le n-k\}.
$$
In particular,
\begin{eqnarray}
K_{\Phi (1^{n+1})}(\gamma_{d,n};1)&=&\ds\prod_{p=1}^{n-1}C_p
\prod_{1\le i<j\le n}\frac{2d+i+j-1}{i+j-1},\label{3.18}\\
K_{\Phi (\eta_k)}(\gamma_{1,n};1)&=&\ds\prod_{p=1}^{k}C_p
\prod_{1\le i<j\le k+1}\frac{2(n-k)+i+j-1}{i+j-1},\label{3.19}
\end{eqnarray}
where $C_m:=\ds\frac{(2m)!}{m!(m+1)!}$ denotes the $m$-th Catalan
number.

By a theorem of Proctor~\cite{Pr}, the product $\ds\prod_{1\le i<j\le n}
\frac{2d+i+j-1}{i+j-1}$ is equal to the number $pp^{\delta_n}(d)$ of (weak)
plane partitions of staircase shape $\delta_n=(n-1,n-2,\ldots ,2,1)$ whose
parts do not exceed $d$. For example, if $n=3$, there exist 5 plane
partitions of shape $\delta_3=(2,1)$ whose parts do not exceed 1. Namely,
plane partitions
$\begin{array}{cc}0&0\\ 0\end{array}$,
$\begin{array}{cc}0&1\\ 0\end{array}$,
$\begin{array}{cc}0&0\\ 1\end{array}$,
$\begin{array}{cc}0&1\\ 1\end{array}$,
$\begin{array}{cc}1&1\\ 1\end{array}$.
More generally, it is not difficult to see that
\begin{equation}
pp^{\delta_n}(1)=\prod_{1\le i<j\le n}\frac{i+j+1}{i+j-1}=C_n,
\label{3.20}
\end{equation}
is equal to the $n$-th Catalan number, and
\begin{equation}
\prod_{(i,j)\in\Omega_{k,n}}\frac{i+j+1}{i+j-1}=C_k\prod_{j=1}^k
\frac{2n-k+j+2}{k+j}.\label{3.21}
\end{equation}

{\bf d.} Combinatorial interpretation of identities (\ref{3.18}) and
(\ref{3.19}).

$\bullet$ Let $\gamma\in Y_n$, and $\eta$ be a composition of
$n$. Denote by $SM_{\eta}(\gamma)$  the set of $n$ by $n$
skew--symmetric integer matrices $m=(m_{ij})$ such that

\vskip 0.3cm

i) $m_{ij}\ge 0$, if $(i,j)\in\Phi (\eta)$;
\vskip 0.3cm

ii) $m_{ij}=0$, if $(i,j)\in\Phi (1^n)\setminus\Phi (\eta)$;

\vskip 0.3cm
iii) $\sum_jm_{ij}=\gamma_i$, $1\le i\le n$. \vskip
0.3cm

It is clear that
$$|SM_{\eta}(\gamma)|=K_{\Phi (\eta)}(\gamma ;1).$$

$\bullet$ Denote by $Cat(n)$ the set of lattice paths from (0,0) to
$(n(n+1),0)$ with steps $(1,1)$ and $(1,-1)$ never falling below the $x$--axis
and passing through the points $P_k=(k(k+1),0)$, $k=0,1,\ldots ,n$. It is
clear that the number of elements in the set $Cat(n)$ is equal to the
product of the first $n$ Catalan numbers $C_1C_2\cdots C_n$.

$\bullet$ Denote by $PP_n(d)$ the set of (weak) plane partitions of
staircase shape $\delta_n$ whose parts do not exceed $d$.

Construct bijections
\begin{eqnarray}
SM_{(1^n)}(\gamma_{d,n-1})&\simeq &Cat(n-2)\times PP_{n-1}(d),
\label{3.22}\\
SM_{(1^k,n-k)}(\gamma_{1,n-1})&\simeq &Cat(k)\times
PP_{k+1}(n-k-1). \label{3.23*}
\end{eqnarray}

{\bf e.} Let $n\ge 3$, consider vector
$$\beta_{d,n}=(d,0,1,\ldots ,n-3,-d-(n-2)(n-3)/2).
$$

Show that
\begin{equation}
K_{\Phi (1^n)}(\beta_{d,n};1)=\prod_{j=1}^{n-3}C_j\left(
\begin{array}{c}d+(n-1)(n-2)/2\\ d\end{array}\right). \label{3.26}
\end{equation}

\begin{prb} Find a $q$--analog of identity (\ref{7.8}), i.e. to define a
statistics $L$ on the set $SM_{(1^n)}(\gamma_{d,n-1})$ with the generating
function
$$\prod_{p=1}^{n-2}C(2,p|q)\prod_{1\le i<j\le
n-1}\frac{1-q^{2d+i+j-1}}{1-q^{i+j-1}},
$$
which is agree with the decomposition (\ref{3.22}).
\end{prb}

Here
$$C(2,m|q):=\ds\frac{1-q}{1-q^{m+1}}\left[\begin{array}{c}2m\\
m\end{array}\right]_q
$$
denotes "the most obvious" $q$--analog of the Catalan number
$C_m$, see Section~\ref{kfpwm}, Exercise~{\bf 1a}.

{\bf f.} Show that

i) there exists a unique continuous piece wise polynomial function
${\cal P}_{\Phi(\eta)}(y_1,\ldots ,y_n)$ on the cone $Y_n$ such
that its restriction to the set $Y_n\cap\Z^n$ of integer points
of the cone $Y_n$ coincides with the Kostant partition function
$K_{\Phi(\eta)}$:
$${\cal P}_{\Phi(\eta)}(\bullet )\vert_{Y_n\cap\Z^n}=
K_{\Phi(\eta)}(\bullet ;1).
$$
Polynomiality domains of the continuous piecewise polynomial
function ${\cal P}_{\Phi(\eta)}(y_1,\ldots ,y_n)$ define a
subdivision of the cone $Y_n$ into convex integral polyhedron
cones. \vskip 0.3cm

\hskip -0.6cm{\bf Question.} For given $\eta$ and $n$ what is the
minimal number $r_n(\eta)$ of polynomiality domains of the
function ${\cal P}_{\Phi(\eta)}(y_1,\ldots ,y_n)$? For example,
$$r_2((1^2))=2,~~ r_3((1^3))=7,~~ r_4((1^4))=48.
$$

ii) the restriction of the function ${\cal P}_{\Phi(\eta)}$ to the
dominant chamber $Y_n^+$ is a polynomial in the variables $y_1,\ldots
,y_{n-2}$.

iii) The homogeneous degree $\left(\begin{array}{c}n-1\\
2\end{array}\right)$ component ${\cal P}_n^{\max}(y)$ of the
polynomial\\ ${\cal P}_{\Phi(1^n)}(y_1,\ldots ,y_n)|_{Y_n^+}$  has
the following form
$${\cal P}_n^{\max}(y)=\sum K_{\Phi(1^{n-2})}(l_1-n+2,l_2-n+3,\ldots
,l_{n-2}-1)\prod_{i=1}^{n-2}\frac{y_i^{l_i}}{l_i!},
$$
summed over all $(n-2)$--tuples of non--negative integers $(l_1,\ldots
,l_{n-2})$ such that

$\bullet$ $l_1+\cdots +l_k\ge kn-\ds\frac{k(k+3)}{2}$, $1\le k\le n-2$,

$\bullet$ $l_1+\cdots +l_{n-2}=\ds\frac{(n-1)(n-2)}{2}$.

\begin{con} If vector $\gamma\in\Z^n$ belongs to the cone (the so--called
dominant chamber)
$$Y_n^+=\{(y_1,\ldots ,y_n)\in Y_n\vert y_1\ge 0,\ldots ,y_{n-1}
\ge 0\},
$$
then the normalized Kostant partition function
$$\left(\begin{array}{c}n-1\\ 2\end{array}\right)!K_{\Phi (1^n)}(\gamma
;1)$$
is an inhomogeneous polynomial of degree
$\left(\begin{array}{c}n-1\\ 2\end{array}\right)$ in the variables
$\gamma_1,\ldots ,\gamma_{n-2}$ with non--negative integer
coefficients.
\end{con}

\begin{con} $K_{\Phi(\eta)}(\gamma ;q)$ is a unimodal polynomial in the
variable $q$.
\end{con}

{\bf g.} Show that
\begin{eqnarray}
{\rm RHS}(\ref{3.18})&=&
\prod_{j=d}^{n+d-2}\frac{1}{2j+1}\left(\begin{array}{c}n+d+j\\
2j\end{array}\right), \label{3.24}\\ {\rm RHS}(\ref{3.19})&=&
\prod_{j=n-k}^{n-1}\frac{1}{2j+1} \left(\begin{array}{c}n+1+j\\
2j\end{array}\right). \label{3.25}
\end{eqnarray}

{\bf h.} (Chan--Robbins polytope)\\ Let $\ld =(\ld_1,\ldots ,\ld_n)$ and
$\mu =(\mu_1,\ldots ,\mu_n)$ be compositions of the same size, without
zero components, and such that
$$\mu_1+\cdots +\mu_{i+1}\ge\ld_1+\cdots
+\ld_i~~{\rm for~all}~~ 1\le i\le n-1.
$$
Let $CR_n^{\ld\mu}$ be the convex polytope in $\R^{n^2}$ of all
points ${\mathbf x}=(x_{ij})_{1\le i,j\le n}$ satisfying the
following conditions

\vskip 0.2cm
$\bullet$ $x_{ij}=0$, if $j>i+1$;

\vskip 0.2cm
$\bullet$ $\sum_jx_{ij}=\ld_i$, $\sum_ix_{ij}=\mu_j$.

\vskip 0.2cm
This is an integral polytope (i.e. all vertices of
the polytope $CR_n^{\ld\mu}$ have integer coordinates) of
dimension $\left(\begin{array}{c}n\\ 2\end{array}\right)$. The
polytope
$$CR_n:=CR_n^{(1^n),(1^n)}$$
was introduced and studied by Chan and Robbins \cite{CR}, and
Chan, Robbins and Yuen \cite{CRY}. It is a face of dimension
$\left(\begin{array}{c}n\\ 2\end{array}\right)$ with $2^{n-1}$
vertices of the Birkhoff polytope ${\cal B}_n$, see
Example~\ref{e7.16}. We will call the polytope $CR_n$ by {\it
Chan--Robbins polytope}.

i) Show that if $k$ is a positive integer, then
\begin{eqnarray}
i(CR_n^{\ld\mu};k)&=&\#|CR_n^{k\ld,k\mu}\cap
\Z^{\scriptsize\left(\begin{array}{c}n\\
2\end{array}\right)}|\label{7.26}\\
&=&K_{\Phi(1^{n+1})}(k\mu_1,k(\mu_2-\ld_1),\ldots
,k(\mu_n-\ld_{n-1}),-k\ld_n;1).\nonumber
\end{eqnarray}
(Hint: it is enough to consider the case $k=1$. If $(x_{ij})\in
CR_n^{\ld\mu}\cap\Z^{\scriptsize\left(\begin{array}{c}n\\
2\end{array}\right)}$, then
$$\ds\sum_{1\le i<j\le n}x_{ij}e_{i,j+1}=(\mu_1,\mu_2-\mu_1,\ldots
,\mu_n-\ld_{n-1},-\ld_n),
$$
and the upper triangular matrix
$$(m_{ij}):=(x_{i,j-1}),~~ 1\le i<j\le n+1,$$
defines an element of the set
$$SM_{\Phi (1^{n+1})}(\mu_1,\mu_2-\ld_1,\ldots ,\mu_n-\ld_{n-1},-\ld_n).
$$
Conversely, if
$$(m_{ij})\in SM_{\Phi (1^{n+1})}(\mu_1,\mu_2-\ld_1,\ldots
,\mu_n-\ld_{n-1},-\ld_n),$$
then the point $x=(x_{ij})_{1\le
i,j\le n}$, where $x_{ij}=m_{j-1,i}$ if $j\le i+1$, and
$x_{ij}=0$, if $j>i+1$, belongs to the polytope $CR_n^{\ld\mu}$.)

Denote by ${\cal P}_{\Phi(1^{n+1})}(y_1,\ldots y_{n+1})$ a unique
continuous piecewise polynomial function on the cone $Y_{n+1}$
such that
$${\cal P}_{\Phi(1^{n+1})}(\bullet )|_{Y_{n+1}\cap\Z^{n+1}}=
K_{\Phi(1^{n+1})}(\bullet ;1),
$$
see \cite{Kir8}, or Exercise~{\bf f} to this Section.

ii) Show that
\begin{equation}
{\cal E}(CR_n^{\ld\mu};t)={\cal P}_{\Phi(1^{n+1})}(t\mu_1,
t(\mu_2-\ld_1),\ldots ,t(\mu_n-\ld_{n-1}),-t\ld_n), \label{7.28}
\end{equation}
where ${\cal E}(CR_n^{\ld\mu};t)$ denotes the Ehrhart polynomial of the
polytope $CR_n^{\ld\mu}$.

In particular, for the Chan--Robbins  polytope one has
\begin{equation}
{\cal E}(CR_n;t)={\cal P}_{\Phi(1^{n+1})}(t,\underbrace{0,\ldots
,0}_{n-1},-t). \label{7.29}
\end{equation}

$\bullet$ Deduce from (\ref{7.29}) and Exercise~2{\bf f}, iii) to
Section~7, that
\begin{equation}
\wt{\rm
vol}(CR_{n+1})=\prod_{j=0}^{n-2}\frac{1}{j+1}\left(\begin{array}{c} 2j\\
j\end{array}\right) \label{7.30}
\end{equation}
is the product of the first $n-1$ Catalan numbers. Formula (\ref{7.30})
for the normalized volume of the Chan--Robbins polytope $CR_n$ was
stated as a conjecture by Chan and Robbins \cite{CR}, and has been
proved by Zeilberger \cite{Z2}.

{\bf i.} Suppose that $n\ge3$ is an integer and that $2\le s\le
r\le n-1$. Follow \cite{CRY} denote by $P_n(r,s)$ the integral
convex polytope
$$P_n(r,s):=CH_n\cap\{ x_{rs}=0\}.$$
Denote by $p_n(r,s)$ the normalized volume of the polytope
$P_n(r,s)$.

\newpage
Show that

$\bullet$ $p_n(r,2)=K_{\Phi (1^{n-1})}(\gamma_r^{(n-1)};1)$, where
$$\gamma_r^{(n-1)}=\left(1,2,\ldots ,n-2,-\frac{(n-1)(n-2)}{2}\right)
-\begin{array}{c}r-1\\ (\underbrace{0,\ldots ,1,\ldots
0}_{n-1})\end{array},~~~ 2\le r\le n-1;
$$

$\bullet$ $p_n(r,s)=p_n(r,2)-p_n(s-1,2)$, $2\le s\le r\le n-1$;

$\bullet$ $p_n(r,2)+p_n(n-r,2)=p_n(n-1,2)$,

~~$p_n(r,2)=p_n(n-1,n-r+1)$, $2\le r\le n-1$,

where by definition we put $p_n(r,s)=0$, if $r<s$. \vskip 0.3cm

\hskip -0.6cm {\bf 3.} (Gelfand--Tsetlin polytope) Let $\ld
=(\ld_1\ge\ld_2\ge\cdots\ge\ld_r>\ld_{r+1}=0)$ and
$\mu=(\mu_1\ge\mu_2\ge\cdots\ge\mu_s>0)$ be partitions such that
$\ld\ge\mu$ with respect to the dominance order.

{\bf a.} Show that
$$\dim G(\ld ,\mu)=(r-1)(s-1)-\left(\begin{array}{c}r\\
2\end{array}\right)-\sum_{i=1}^{r}\left(\begin{array}{c}\ld_i'
-\ld_{i+1}'\\ 2\end{array}\right).
$$

{\bf b.} Denote by $a(\ld ,\mu)$ the lowest degree of $q$ which is
present in $K_{\ld\mu}(q)$, i.e.,
$$K_{\ld\mu}(q)=b(\ld ,\mu)q^{a(\ld ,\mu)}{\rm (1+higher ~degree~ terms)}.
$$

Show that if $n\in\Z_{\ge 1}$, then
$$a(n\ld ,n\mu)=na(\ld ,\mu).$$

\begin{con}\label{c7.12} Generating function
$$\sum_{n\ge 0}b(n\ld ,n\mu)t^n$$
is a rational function of the form
$Q_{\ld\mu}(t)/(1-t)^{r(\ld ,\mu)+1}$, where $r(\ld
,\mu)\in\Z_{\ge 0}$, and $Q_{\ld\mu}(t)$ is a polynomial of
degree $\le r(\ld,\mu)$ with non--negative integer coefficients.
More generally,
$$\sum_{k\ge 1}\sum_{n\ge 1}b((kn-1)\ld ,(kn-1)\mu
)t^{k-1}z^{n-1}=Q_{\ld\mu}(t,z)/((1-t)(1-z))^{r(\ld ,\mu)+1},
$$
where
$$Q_{\ld\mu}(z,t)=Q_{\ld\mu}(t,z)$$
is a polynomial with non--negative integer coefficients.
\end{con}

\begin{ex} {\rm Take $\ld =(7432)$ and $\mu =(4432111)$. Then $\dim
G(\ld ,\mu)=12$, $a(n\ld ,n\mu)=5n$, and
$$b(n\ld ,n\mu)=\left(\begin{array}{c}n+3\\ 4\end{array}\right)+
\left(\begin{array}{c}n+4\\ 4\end{array}\right).$$
Hence,}
\begin{eqnarray*}
&&\sum_{n\ge 0}b(n\ld ,n\mu )t^n=\frac{1+t}{(1-t)^5};\\
&&Q_{\ld\mu}(t,z)=1+z+(1+15z+6z^2)t+(6z+15z^2+z^3)t^2+(z^2+z^3)t^3.
\end{eqnarray*}
\end{ex}

{\bf c.} Show that for any $n,m\in\Z_{\ge 0}$,
$$K_{(n+m)\ld ,(n+m)\mu}(q)\le K_{n\ld ,n\mu}(q)K_{m\ld ,m\mu}(q).
$$

{\bf d.} Let ${\Pe}$ be an integral convex polytope in $\R^d$ of
dimension $d$. Show that
$$\sum_{k_1,k_2,\ldots ,k_r\ge 1}i({\Pe},k_1k_2\cdots k_r-1)
\prod_{j=1}^rt_j^{k_j-1}=Q_{\Pe}(t_1,\ldots
,t_r)/\prod_{j=1}^r(1-t_j)^{d+1},
$$
where $Q_{\Pe}(t_1,\ldots ,t_r)$ is a polynomial in the variables
$t_1,\ldots ,t_r$ with integer coefficients. It looks an interesting
combinatorial problem to describe all integral convex polytopes ${\Pe}$
such that the corresponding polynomial $Q_{\Pe}(t_1,t_2)$ has {\it
non--negative} integer coefficients.

\hskip -0.6cm{\bf Question.} Does there exist an integral convex
polytope $\Gamma (\ld ,\mu)$ such that
$$i(\Gamma (\ld ,\mu);k)=b(k\ld ,k\mu),
$$
for all integers $k\ge 1$? If so, is it true that in this case $\Gamma
(\ld ,\mu)=\Gamma (\mu',\ld')$?

{\bf e.} (MacMahon polytope and Narayana numbers again)\\
Take $\ld =(n+k,n,n-1,\ldots ,2)$ and $\mu
=\ld'=(n,n,n-1,n-2,\ldots ,2,1^k)$. Show that if $n\ge k\ge 1$,
then for any positive integer $N$
\vskip 0.2cm

$\bullet$ $a(N\ld,N\mu)=(2k-1)N$;

$\bullet$ $b(N\ld,N\mu)={\rm dim}V_{((n-k+1)^{k-1})}^{\g
l(N+k-1)}=\ds\prod_{i=1}^{k-1}\prod_{j=1}^{n-k+1}\frac{N+i+j-1}{i+j-1}.$

In other words, $b(N\ld,N\mu)$ is equal to the number of (weak)
plane partitions of rectangular shape $((n-k+1)^{k-1})$ whose
parts do not exceed $N$. According to Exercise~1, {\bf c},
Section~\ref{kfpwm}, $b(N\ld,N\mu)$ is equal also to the number
$i(\M_{k-1,n-k+1};N)$ of rational points ${\bf x}$ in the
MacMahon polytope $\M_{k-1,n-k+1}$ such that the points $N{\bf x}$
have integer coordinates. Therefore, in this example one can take
$$\Gamma (\ld,\mu)=\M_{k-1,n-k+1}.$$
It follows from (\ref{2.17}) that the generating function for
numbers $b(n\ld,n\mu)$ has the following form
$$\sum_{n\ge 0}b(n\ld,n\mu)t^n=\left(\sum_{j=0}^{(k-2)(n-k)}
N(k-1,n-k+1;j)t^j\right)/(1-t)^{(k-1)(n-k+1)+1},
$$
where $N(k,n;j)$, $0\le j\le (k-1)(n-1)$, denote rectangular
Narayana's numbers, see Section~\ref{kfpwm}, Exercise~1.

$\bullet$ Show that if $r:=k-\left(\begin{array}{c}n+2\\
2\end{array}\right)\ge 0$, then $b(\ld,\mu)=1$, and
$$a(\ld,\mu)=2\left(\begin{array}{c}n+3\\ 3\end{array}\right)
+(n+1)(2r-1)+\left(\begin{array}{c}r\\ 2\end{array}\right);
$$

$\bullet$ if $1\le k<\left(\begin{array}{c}n+2\\
2\end{array}\right)$, then there exists a unique $p$, $1\le p\le
n$, such that
$$(p-1)(2n-p+4)/2<k\le p(2n-p+3)/2.
$$
In this case
$$a(\ld,\mu)=p(2k-(p-1)n-p)+2\left(\begin{array}{c}p\\
3\end{array}\right), $$
and one can take $\Gamma(\ld,\mu)$ to be
equal to the MacMahon polytope $\M_{r(k),s(k)}$ with
$r(k):=k-1-(p-1)(2n-p+4)/2$, and $s(k):=p(2n-p+3)/2-k$.

This Exercise gives some flavor how intricate the piecewise
linear function $a(\ld,\mu)$ may be.
\vskip 0.2cm

\hskip -0.6cm{\bf 4.} (Generalized saturation conjecture)

{\bf a.} For any partition $\ld$ and a dominant sequence of
rectangular shape partitions $R=((\mu_a^{\eta_a}))^p_{a=1}$
define the numbers $a(\ld ,R)$ and $b(\ld ,R)$ via decomposition
$$K_{\ld R}(q)=b(\ld ,R)q^{a(\ld ,R)}+{\rm higher~degree~ terms}.
$$
For any non--negative integer $n$ denote by $nR$ the dominant
sequence of rectangular shape partitions
$$nR:=((n\mu_a)^{\eta_a}).$$

\begin{con}\label{c7.24} {\rm (Generalized saturation conjecture)}
$$a(n\ld ,nR)=na(\ld ,R).
$$
\end{con}

\begin{con}\label{c7.25} i) {\rm (Generalized Fulton's conjecture)}
If $b(\ld ,R)=1$, then
$$b(n\ld ,nR)=1$$
for any non--negative integer $n$; conversely, if $b(n\ld ,nR)=1$
for some positive integer $n$, then $b(\ld ,R)=1$.

ii) More precisely, there exist non--negative integers $d$ and
$a_0,a_1,\ldots, a_d$ (depending on $\ld$ and $R$) such that
$a_0=1$, and
$$b(n\ld,nR)=\ds\sum_{k=0}^da_k\left(\begin{array}{c}n+k\\
d\end{array}\right)
$$
for all $n\ge 0$. Moreover, $b(n\ld ,nR)$ is a polynomial in $n$
with non--negative rational coefficients.
\end{con}

Therefore, if $b(n\ld,nR)=1$ for some non--negative integer $n$,
then $d=0$ and $a_0=1$, and so $b(N\ld,N\mu)=1$ for all
non--negative integers $N$.

Recall that Fulton's conjecture says that a condition
$c_{\ld\mu}^{\nu}=1$ is equivalent to $c_{N\ld,N\mu}^{N\nu}$ being
one for any $N\in\N$, see e.g., A.~Buch (Enseign. Math. (2) {\bf
46} (2000), 43--60).

It is well--known that weight multiplicities $\dim
V_{\ld}(\mu)=K_{\ld\mu}$ may be considered as a limiting case of
the Littlewood--Richardson numbers, i.e. the triple tensor
product multiplicities. The Fulton conjecture is a far
generalization of the following (well--known) fact:

{\it Let $\ld$ and $\mu$ be partitions, then $K_{\ld\mu}=1$ if
and only if $K_{n\ld,n\mu}=1$ for all positive integers $n$.}

This fact is a simple corollary of the Berenstein--Zelevinsky
weight--multiplicity--one criteria, see e.g., A.~Berenstein and
A.~Zelevinsky (Functional Anal. Appl. {\bf 24} (1990), 259--269),
A.N.~Kirillov, {\it Generalization of the Gale--Ryser theorem},
European Journ. of Combinatorics, {\bf 21} (2000), 1047--1055, or
R.~Stanley \cite{St3}, Exercise~7.13, p.451.
\begin{con}\label{c7.26}
\begin{eqnarray*} i) ~~~a(\ld,\mu)&=&a(\mu',\ld');\\
ii) ~~~b(\ld,\mu)&=&b(\mu',\ld').
\end{eqnarray*}
\end{con}

Conjecture~\ref{c7.26}, i) was stated by I.~Macdonald \cite{Ma},
Chapter~III, \S 6, Example~3 on p.243. Conjecture~\ref{c7.26},
ii) was stated for the first time in \cite{Kir1}, see also
\cite{Kir6}.

Let us explain briefly why Conjecture~\ref{c7.24} may be
considered as a generalization of {\it saturation conjecture}.
First of all, let us recall a formulation of the saturation
conjecture (now theorem due to A.~Knutson and T.~Tao, J. Amer.
Math. Soc. {\bf 12} (1999), 1055--1090):

{\it Let $\ld,\mu,\nu$ be partitions, denote by
$c_{\ld,\mu}^{\nu}$ the corresponding Littlewood--Richardson
number, i.e.
$$c_{\ld,\mu}^{\nu}={\rm Mult}[V_{\nu}:V_{\ld}\otimes
V_{\mu}].
$$
If for some non--negative integer $N$ one has
$c_{N\ld,N\mu}^{N\nu}=0$, then $c_{\ld,\mu}^{\nu}=0$.}

\vskip 0.2cm

It was first shown in 1994 by A.~Klyachko (Selecta Math. {\bf 4}
(1998), 419--445) that the saturation conjecture implies the Horn
conjecture about eigenvalues of a sum of two Hermitian matrices
(A.~Horn, Pacific Journ. Math. {\bf 12} (1962), 225--241). We
refer the reader to survey articles by W.~Fulton (Bull. Amer.
Math. Soc. {\bf 37} (2000), 209--249; S\'eminaire Bourbaki, vol.
1997/98, Ast\'erisque No.252 (1998), Exp. No.845, 255--269) for
historical overview, description and solution of Horn's problem,
and explanation of connections between Horn's and saturation
conjectures.

Finally, let us show how the saturation conjecture follows from
Conjecture~\ref{c7.24}. First of all, according to Exercise~3,
Section~\ref{pkp}, if $\ld,\mu,\nu$ be partitions, then
$$K_{\wt\ld,\wt R}(q)=q^{Q_{\ld_1}(\nu)-|\ld|}\left\{
c_{\ld,\mu}^{\nu}+{\rm higher ~degree~terms}\right\}.
$$
Therefore, if $c_{N\ld,N\mu}^{N\nu}\ne 0$ for some positive
integer $N$, then
$$a(N\wt\ld,N\wt R)=Q_{N\ld_1}(N\nu )-|N\ld |= \sum_a\min
(N\ld,N\nu_a)-N|\ld |=N(Q_{\ld_1}(\nu)-|\ld |),
$$
and therefore,
$$a(\wt\ld,\wt R)=Q_{\ld_1}(\nu)-|\ld |.
$$
The latter equality means that $c_{\ld\mu}^{\nu}\ne 0$. Similarly
one can show that if $c_{\ld,\mu}^{\nu}\ne 0$ then
$c_{N\ld,N\mu}^{N\nu}\ne 0$ for any positive integer $N$.

\qed

As a byproduct, we obtain a criteria telling  when the
Littlewood--Richardson number $c_{\nu\mu}^{\nu}\ne 0$. Namely,
let $\ld,\mu,\nu$ be partitions such that $|\nu|=|\ld|+|\mu|$,
$l(\ld )\le p$, $l(\mu )\le s$. Define partition
$$\Lambda =(\ld_1+\mu_1,\ld_1+\mu_2,\ldots ,\ld_1 +\mu_s,\ld_1,\ldots
,\ld_p)
$$
and a dominant rearrangement $R$ of the sequence of
rectangular shape partitions $\{\nu,(\ld_1^s)\}$. Then
$c_{\ld\mu}^{\nu}\ne 0$ if and only if
$$a(\Lambda ,R)=Q_{\ld_1}(\nu)-|\ld|,
$$
and if the latter equality  holds, then
$$c_{\ld\mu}^{\nu}=b(\Lambda ,R).
$$
In order this criteria happened to be really effective, one has
to know a direct combinatorial interpretation of the numbers
$a(\Lambda ,R)$.

\begin{rem} {\rm More generally, for any two partitions $\ld$ and
$\mu$ of the same size, and a composition $\eta$, such that
$|\eta|\ge l(\mu)$, one can define the numbers $a(\ld,\mu;\eta)$
and $b(\ld,\mu;\eta)$ via the decomposition
$$K_{\ld\mu\eta}(q)=b(\ld,\mu;\eta)q^{a(\ld,\mu;\eta)}+ {\rm
higher~degree ~terms}.
$$
It is natural to suggest the following 

\begin{con}\label{c7.28}

i) ~~$a(n\ld,n\mu;\eta)=na(\ld,\mu;\eta)$;\\
ii) ~~$b(n\ld,n\mu;\eta)$ is a polynomial in $n$ with
non--negative
rational coefficients;\\
iii) ~$b(n\ld,n\mu;\eta)=1$ for some positive integer $n$ if and
only if $b(\ld,\mu;\eta)=1$;\\
iv) ~assume additionally that $|\eta|\ge\ld_1$, then
$a(\ld,\mu;\eta)=a(\mu',\ld';\eta)$;\\
v)~ the generating function $\sum_{n\ge 0}b(n\ld,n\mu;\eta)t^n$
is a rational function in $t$ of the form\\
$P_{\ld\mu\eta}(t)/(1-t)^{r(\ld\mu\eta)+1}$, where
$r(\ld\mu\eta)\in\Z_{\ge 0}$ and $P_{\ld\mu\eta}(t)$ is a
polynomial of degree $\le r(\ld\mu\eta)$ with non--negative
integer coefficients.
\end{con}

Let us remark that in general
$$b(\ld,\mu;\eta)\ne b(\mu',\ld';\eta).$$
For example, take
$$\ld =(7654321),~~
\mu=(4^5,2,1^6)~~{\rm and}~~ \eta=(4,1,1,2,4).
$$
Then $\mu'=(12,6,5,5)$ and $K_{\ld\mu\eta}(q)=q^8(5,17,24,17,5)$,
but
$$K_{\mu'\ld'\eta}(q)=q^8(1,1,1).$$}
\end{rem}

{\bf b.} Let $\ld,\mu$ and $\nu$ be partitions. Show that
$$a(\ld+\nu,\mu+\nu)\le a(\ld,\mu),
$$
and if $l(\nu)<l(\ld)$, then
\begin{eqnarray*}
a(\ld+\nu,\mu+\nu)&=&a(\ld,\mu),\\
b(\ld+\nu,\mu+\nu)&\ge &b(\ld,\mu).
\end{eqnarray*}

See also Exercise~1 at the end of Section~\ref{sb}.

\begin{prb} Show that $a(\ld ,R)$ and $a(\ld,\mu;\eta)$ are (continious)
piecewise linear functions of $\ld$ and $R$, and $\ld,\mu$ and
$\eta$ correspondingly, and compute these functions explicitly.
\end{prb}

{\bf c.} Let $\ld,\mu$ be partitions and $\nu$ be a composition
such that both $\ld +\nu$ and $\mu +\nu$ are partitions. Show that
$$a(\ld +\nu,\mu +\nu)\le a(\ld,\mu),$$
and if $l(\nu)<l(\ld)$, then
$$a(\ld +\nu,\mu +\nu)=a(\ld,\mu).$$

Note that in the case  under consideration the difference
$$K_{\ld +\nu,\mu +\nu}(q)-K_{\ld,\mu}(q)$$
may have negative coefficients. However, it looks plausible that
all coefficients of the latter difference have the same sign.

{\bf d.} (Geometrical interpretation of Conjecture~\ref{c7.26})
Let $\ld$ and $\mu$ be partitions of the same size, such that
$\ld\ge\mu$, and $N$ be an integer. Follow \cite{Kir2,Kir6}
consider the compact convex polytope $\Gamma^{(1)}(\ld,\mu)$ in
$\R^{N^2}$ of all points ${\bf x}=(x_{ij})_{1\le i,j\le N}$
satisfying the following conditions
\vskip 0.2cm

i) $\ds\sum_jx_{ij}=\ld_i$, ~~ $\sum_ix_{ij}=\mu_j'$;

ii) $\ds\sum_{j\le n}(x_{kj}-x_{k+1,j})\ge 0$ and $\sum_{i\ge
k+1}(x_{in}-x_{i,n+1})\ge 0$ for all integers $k,n\ge 1$.

We will call this polytope by {\it configurations} polytope. One
can show \cite{Kir6}, Section~5, that the set of integer points
$\Gamma^{(1)}_{\Z}(\ld,\mu)$ of the configuration polytope is in
one--to--one correspondence with the set $C(\ld;\mu)$ of
admissible configurations of type $(\ld;\mu)$. Therefore,
$\Gamma_{\Z}^{(1)}(\ld,\mu)\ne 0$ if and only if $\ld\ge\mu$.

For any ${\bf x}\in\R^{N^2}$ let us set
$$c({\bf x})=\sum_{i,j}x_{ij}^2,$$
and define
$$\wt a(\ld,\mu)=\min\{c({\bf x})~\vert ~{\bf x}\in\Gamma^{(1)}_{\Z}(\ld,\mu)\}.
$$
It easy to see that
$$a(\ld,\mu)=\frac{1}{2}(\wt a(\ld,\mu)-|\ld|).$$
Finally consider the convex set
$$\Gamma_0^{(1)}(\ld,\mu)=\{{\bf x}\in\Gamma^{(1)}(\ld,\mu)~|~ c({\bf x})=\wt
a(\ld,\mu)\}.
$$
$\bullet$ Show that
$$|\Gamma_0^{(1)}(\ld,\mu)\cap\Z^{N^2}|=b(\ld,\mu).
$$

On the other hand, consider the compact convex polytope
$\Gamma^{(2)}(\ld,\mu)$ in $\R^{N^2}$ of all points ${\bf
x}=(x_{ij})_{1\le i,j\le N}$ satisfying the following conditions
\vskip 0.2cm

i) $\ds\sum_jx_{ij}=\ld_i$, ~~ $\sum_ix_{ij}=\mu_j'$;

ii) $\ds\sum_{j\le n}(x_{kj}-x_{k+1,j})\le\ld_k-\ld_{k+1}$ and
$\sum_{i\ge k+1}(x_{in}-x_{i,n+1})\le\mu_n'-\mu_{n+1}'$ for all
$k,n\ge 1$.

$\bullet$ Show that the set $\Gamma_{\Z}^{(2)}(\ld,\mu)$ of
integer points of the polytope $\Gamma^{(2)}(\ld,\mu)$ is in
one--to--one correspondence with the set $C(\mu',;\ld')$ of
admissible configurations of type $(\mu';\ld')$, and construct
such correspondence explicitly.

Finally, consider the convex set
$$\Gamma^{(2)}_0(\ld,\mu)=\{{\bf x}\in\Gamma^{(2)}(\ld,\mu)~|~c({\bf
x})=\wt a(\ld,\mu)\}.
$$
Conjecture~\ref{c7.26} is equivalent to the following statements
\vskip 0.2cm

$\bullet$ $\wt a(\mu',\ld'):=\min\{c({\bf x})|{\bf
x}\in\Gamma^{(2)}_{\Z}(\ld,\mu)\}=\wt a(\ld,\mu)$;

$\bullet$ $|\Gamma_0^{(2)}(\ld,\mu)\cap\Z^{N^2}|=b(\ld,\mu)$.

\vskip 0.2cm

\hskip -0.6cm{\bf 5.} Let $\ld$ be a partition and $R$ be a
dominant sequence of rectangular shape partitions

i) Show that the generating function
$$Q_{\ld R}(q,t):=\sum_{n\ge 0}K_{n\ld,nR}(q)t^n$$
is a {\it rational} function in $q$ and $t$,
say,
$$Q_{\ld R}(q,t)=\ds\frac{N_{\ld R}(q,t)}{D_{\ld R}(q,t)},$$
where the numerator $N_{\ld R}(q,t)$ and the denominator $D_{\ld
R}(q,t)$ are suppose to be mutually prime.

ii) Let $\ld$ and $R$ be as in i). Show that there exists the
limit of {\it cocharge} parabolic Kostka polynomials
$$\lim_{n\to\infty}{\overline K}_{n\ld,nR}(q),
$$
where by definition
$${\overline K}_{\ld R}(q)=q^{n(R)}K_{\ld R}(q^{-1}).$$

\begin{con}\label{c7.30} i) The denominator of the rational
function $Q_{\ld R}(q,t)$ has a product form
$$D_{\ld R}(q,t)=\prod (1-q^at)^{k_a}$$
for some finite set of non--negative integers $\{k_{j_1},\ldots
,k_{j_p}\}$;

ii) if $b(\ld ,R)=1$, then there exists the limit of normalized
parabolic Kostka polynomials~~
$$\lim_{n\to\infty}q^{-na(\ld,R)}K_{n\ld,nR}(q);$$

iii) the generating function $\sum_{n\ge 0}K_{n\ld ,nR}(1)t^n$ is
a rational function in $t$ of the form
$$\ds\frac{P_{\ld R}(t)}{(1-t)^{d(\ld,R)+1}},$$
where $d(\ld,R)\in\Z_{\ge 0}$, $P_{\ld R}(1)\ne 0$, and $P_{\ld
R}(t)$ is a polynomial with non--negative integer coefficients.
\end{con}

\begin{prb}\label{p7.31} i) Find combinatorial inerpretation(s)
for the number $b(\ld,R)$ in terms of objects similar to either
the Berenstein--Zelevinsky triangles \cite{BZ}, or the
Knutson--Tao honeycombs and hives {\rm (A.~Knutson and T.~Tao,
[{\it ibid}]; A.~Buch, Enseign. Math. (2) {\bf 46} (2000),
43--60)}, or the Gleiser--Postnikov web functions {\rm (O.~Gleizer
and A.~Postnikov, Itern. Math. Res. Notes {\bf 14} (2000),
741--774)}, or domino tableaux {\rm (see e.g., C.~Carre and
B.~Leclerc, Journ. Alg. Comb. {\bf 4} (1995), 201--231)}.

ii) Find a $q$--analog of numbers $b(\ld ,R)$ which generalizes
the $q$--analog $c_{\ld\mu}^{\nu}(q)$ of Littlewood--Richardson's
numbers introduced by C.~Carre and B.~Leclerc [{\it ibid}].

iii) When does the number $b(\ld ,R)$ equal to 1 ?
\end{prb}

\begin{exs} {\rm 1) Take $\ld =(4422)$ and $R=(3,3,(2,2),1,1)$.
Then
\begin{eqnarray*}
\sum_{n\ge0}K_{n\ld,nR}(q)t^n&=&\frac{(1+q^4t)(1+q^5t)}{(1-q^2t)
(1-q^3t)(1-q^4t)^3(1-q^6t)},\\
\sum_{n\ge 0}K_{n\ld,nR}t^n&=&\frac{(1+t)^2}{(1-t)^6}.
\end{eqnarray*}
Note that in this case $a(\ld,R)=2$ and $b(\ld,R)=1$,
$$K_{\ld R}(q)=q^2(11411).$$

2) Take $\ld =(52)$ and $\mu =(2221)$. Then
\begin{eqnarray*}
\sum_{n\ge0}K_{n\ld,n\mu}(q)t^n&=&\frac{1+q^5t+q^6t}{(1-q^4t)(1-q^5t)(1-q^7)},\\
\sum_{n\ge 0}K_{n\ld,n\mu}t^n&=&\frac{1+2t}{(1-t)^3}.
\end{eqnarray*}
Note that in this case $a(\ld,\mu)=4$ and $b(\ld,\mu)=1$.

3) Take $\ld =(422)$ and $\mu =(22211)$. Then
\begin{eqnarray*}
\sum_{n\ge0}K_{n\ld,n\mu}(q)t^n&=&\frac{1+(q^4+q^5+q^6)t-
(q^8+q^9+q^{10})t^2-q^{14}t^3}{(1-q^3t)^2(1-q^4t)(1-q^5t)^2(1-q^7t)},\\
\sum_{n\ge 0}K_{n\ld,n\mu}t^n&=&\frac{1+4t+t^2}{(1-t)^5}.
\end{eqnarray*}
Note that in this case $a(\ld,\mu)=3$ and $b(\ld,\mu)=2$.

4) Take $\ld =(732)$ and $\mu =(432111)$. Then
$$\sum_{n\ge0}K_{n\ld,n\mu}t^n=\frac{1+34t+189t^2+228t^3+57t^4
+2t^5}{(1-t)^8}.
$$
Here $a(\ld,\mu)=5$ and $b(\ld,\mu)=3$, since
$K_{\ld\mu}(q)=q^5(3,7,9,9,7,4,2,1).$

It is not difficult to check that the Gelfand--Tsetlin polytope
$G(\ld ,\mu)$ is an integral one, and its normalized volume is
equal to $511$.

5) Take $\ld =(44)$, $\mu =(1^8)$. Then
$$\sum_{n\ge
0}K_{n\ld,n\mu}(q)=\frac{N_{\ld\mu}(q,t)}{D_{\ld\mu}(q,t)},
$$
where
\begin{eqnarray*}
N_{\ld\mu}(q,t)&=&1+q^{16}(1+q+q^2+q^3+q^4+q^5+q^6)(t-q^{53}t^4)\\
&+&q^{33}(1+q+q^2+q^3+2q^4+2q^5+2q^6+q^7+q^8\\
&+&q^9+q^{10})(t^2-q^{15}t^3)-q^{91}t^5;\\
D_{\ld\mu}(q,t)&=&(1-q^{12}t)(1-q^{14}t)(1-q^{15}t)(1-q^{16}t)
(1-q^{18}t)\\
&\cdot &(1-q^{20}t)(1-q^{24}t);\\
\sum_{n\ge 0}K_{n\ld,n\mu}t^n&=&
\frac{1+8t+22t^2+8t^3+t^4}{(1-t)^6}.
\end{eqnarray*}
In this example
$$K_{\ld\mu}(q)=q^{12}(1011212121101),$$
${\rm dim~}G(\ld,\mu)=5$, and the normalized volume $\wt{\rm
vol~}G(\ld,\mu)=40$.

6) Take $\ld =(333)$, $\mu =(1^9)$. Then
$$\sum_{n\ge
0}K_{n\ld,n\mu}t^n=\frac{1\!+31t\!+469t^2\!+2113t^3\!+3466t^4\!+2113t^5
\!+469t^6\!+31t^7\!+t^8}{(1-t)^{11}}.
$$
In this example, $K_{\ld\mu}=42$, $a(\ld,\mu)=9$, $b(\ld,\mu)=1$,
${\rm dim~}G(\ld,\mu)=10$. It is not difficult to check that the
Gelfand--Tsetlin polytope $G(\ld,\mu)$ is an integral one, and
its normalized volume is equal to 8694.

7) Take $\ld =(52^31)$, $\mu=(2^6)$. One can check that there
exists only one admissible configuration of type $(\ld;\mu)$,
i.e. $|C(\ld;\mu)|=1$. On the other hand $|C(2\ld;2\mu)|=2$. This
observation shows that in the case under consideration the
configurations polytope $\Gamma^{(1)}(\ld,\mu)$ cannot be an
integral polytope, but only a rational one.

It is not difficult to check that
\begin{eqnarray*}
K_{n\ld,n\mu}(q)&=&q^{4n}\left[\begin{array}{c}2n+4\\3\end{array}\right]_q
\left[\begin{array}{c}n+3\\3\end{array}\right]_q
\frac{1-q^{n+1}}{1-q^4},\\
\sum_{n\ge
0}K_{n\ld,n\mu}t^n&=&\frac{1+32t+128t^2+104t^3+15t^4}{(1-t)^8}.
\end{eqnarray*}
Therefore, $\dim G(\ld,\mu)=7$, ${\wt{\rm vol~}}G(\ld,\mu)=280$,
and
$$\ds\lim_{n\to\infty}q^{-na(\ld,\mu)}K_{n\ld,n\mu}(q)=
\frac{1}{(q;q)_3(q;q)_4}.
$$
}
\end{exs}

\hskip -0.6cm{\bf 6.} Example of a convex integral polytope such
that its Ehrhart's polynomial has {\it negative} integer
coefficient.

{\bf a.} Let $\tau_r\subset\R^3$ be the tetrahedron whose vertices are
the points (0,0,0), (1,0,0), (0,1,0) and $(1,1,r)$, where $r$ is a
positive integer.

Show that

$\bullet$ ${\cal E}(\tau_r;t)=(rt^3+6t^2+(12-r)t+6)/6$;

$\bullet$ $\ds\sum_{n\ge 0}i(\tau_r;n)t^n=\frac{1+(r-1)t^2}{(1-t)^4}$.\\
Hence, the tetrahedron $\tau_r$ contains no lattice points other than
its vertices, has volume $r/6$, and if $r>12$, then the Ehrhart
polynomial ${\cal E}(\tau_r;t)$ has {\it negative} coefficient. Note
that the $\delta$--vector of the tetrahedron $\tau_{r}$ is equal to
$(1,0,r-1)$, which is {\it not} unimodal.

This example is due to J.~Reeve (Proc. London Math. Soc. {\bf 7}
(1957), 378-395). Note also, see \cite{Hi}, p. 111, that there
exists an integral simplex $\Delta\subset\R^8$ of dimension 8,
which satisfies the condition $(\Delta
-\partial\Delta)\cap\Z^8\neq\emptyset$, and such that its
$\delta$--vector is equal to $\delta (\Delta
)=(1,8,15,13,15,14,14,14,5)$, which is {\it not} unimodal.

\hskip -0.6cm {\bf Question.} Is it true or not that if the
$\delta$--vector of a convex integral polytope ${\cal P}$ is
unimodal, then the Ehrhart polynomial ${\cal E}({\cal P};t)$ of
the polytope ${\cal P}$ has non--negative integer coefficients?

{\bf b.} Compute the Ehrhart polynomial and $\delta$--vector of
the simplex in the space $\R^d$ with the following vertices
$$(\underbrace{0,\ldots ,0)}_d,e_1,\ldots
,e_{d-1},(\underbrace{1,\ldots ,1}_{d-1},r),
$$
where $e_i=(e_{ik})^d_{k=1}$, and $e_{ik}=\delta_{i,k}$.

\hskip -0.6cm{\bf 7.} Let $\ld$ be a partition, $\ld
=(p_1^{m_1},\ldots ,p_k^{m_k})$, where $p_1>p_2>\cdots >p_k>0$
and each $m_i\ge 1$. Denote by $Y_k(\ld)$ the convex
(non--compact) polyhedron in $\R^{k+1}$ of all points ${\bf
y}=(y_1,\ldots ,y_{k+1})$ satisfying the following conditions

$\bullet$ $y_1+\cdots +y_i\ge 0$ for all integers $i$, $1\le i\le k$;

$\bullet$ $-m_i\le y_i\le p_{i-1}-p_i$ for all integers $i$, $2\le i\le
k$;

$\bullet$ $y_1+\cdots +y_{k+1}=0$.

i) Show that the set of integer points of the polyhedron $Y_k(\ld)$ is
in one--to--one correspondence with the set of {\it vexillary}
permutations $w\in S_{\infty}$ of shape $\ld (w)=\ld$.

ii) Denote by $Y_k(\ld ;n)$ the convex polytope in $\R^{k+1}$ of
all points ${\bf y}\in Y_k(\ld)$ that additionally satisfy the
following condition

$$p_i+y_1+\cdots +y_i+m_1+\cdots +m_i\le n
$$
for all integers $i$, $1\le i\le k$.

Show that the number of integer points of the (rational) polytope
$Y_k(\ld ;n)$ is equal to the number of {\it vexillary}
permutation $w\in S_n$ of shape $\ld (w)=\ld$. We refer the
reader to \cite{Ma1}, Chapter~I, for definition and basic
properties of {\it vexillary} permutations. Note that formula
(\ref{7.9}) with $p=3$ gives an explicit combinatorial expression
for the number of vexillary permutations in the symmetric group
$S_n$.

\hskip -0.6cm{\bf 8.} Let $P(q)$ be a polynomial with
non--negative integer coefficients.

i) Show that there exists a non--negative integer $N$, depending
on $P(q)$, such that the product $(1+q)^NP(q)$ is a {\it unimodal}
polynomial.

We denote by $ud(P)$ the least non--negative integer $N$ with
this property.

ii) Show that if $n\ge 2$, then
$$ud(1+q^n)=n^2-3.
$$

\section{Stable behavior of Kostka--Foulkes polynomials}
\label{sb}
\neweq

Given partitions $\alpha$ and $\beta$ with $l(\alpha)+l(\beta)\le n$,
define
$$V_{\alpha ,\beta}:=V_{\alpha ,\beta}^{[n]}
$$
as the Cartan piece in $V_{\alpha}\otimes V_{\beta}^*$, i.e. the
irreducible $\g l(n)$--submodule generated by the tensor product
of the highest weight vectors in each factor. It is well--known
that $V_{\alpha ,\beta}^{[n]}=V_{\gamma}^{[n]}$ for
$$\gamma:=[\alpha ,\beta]_n=(\alpha_1+\beta_1,\ldots ,\alpha_r+\beta_1,
\underbrace{\beta_1,\ldots
,\beta_1}_{n-r-s},\beta_1-\beta_s,\ldots ,\beta_1-\beta_2),
$$
where $r=l(\alpha)$ and $s=l(\beta)$. For example,
$$V_{(0),(0)}^{[n]}=\C,~~ V_{(1),(1)}^{[n]}=\g~~{\rm (the~adjoint~
representation).}
$$

In Sections~\ref{gemtr} and \ref{ipsf} we had investigated the
properties of the mixed tensor representations $V_{[\al ,\beta ]_n}$ and
their zero--weight subspaces. In this Section we are going to study the
stable behaviour of more general families of weight subspaces in
$V_{[\al ,\beta ]_n}$.

 Let $\ld ,\mu ,\nu$ be partitions, $|\ld|=|\mu|$, consider a family of
 partitions $\ld_n:=\ld\oplus\nu^n$ and $\mu_n:=\mu\oplus\nu^n$. Recall
 that if $\ld$ and $\mu$ be partitions, then $\ld\oplus\mu$ denotes the
 partition whose parts are those of $\ld$ and $\mu$, arranged in
 descending order, and $\mu\oplus\nu^n=\mu\oplus
 \underbrace{\nu\oplus\cdots\oplus\nu}_n$.

 \begin{pr}\label{p8.1}{\rm (\cite{Kir1,Kir6})}
 i) The number $\#|C(\ld_n,\mu_n)|$ of admissible
 configurations of type $(\ld_n;\mu_n)$ is finite and independent on $n$
 if $n$ is big enough.

 ii) ("Gupta conjecture") If $n\ge 0$, then
 \begin{equation}
 K_{\ld_{n+1},\mu_{n+1}}(q)\ge K_{\ld_n,\mu_n}(q),\label{8.1}
 \end{equation}
 i.e. the difference $K_{\ld_{n+1},\mu_{n+1}}(q)-K_{\ld_n,\mu_n}(q)$ is
 a polynomial with non--negative integer coefficients.

 iii) There exist the limits
 \begin{eqnarray*}
 \lim_{n\to\infty}K_{\ld_n,\mu_n}(q)&=&Z_{\ld\mu}^{\nu}(q),\\
 \lim_{n\to\infty}{\overline K}_{\ld_n,\mu_n}(q)&=&Y_{\ld\mu}^{\nu}(q),
 \end{eqnarray*}
 which are rational functions of $q$.
 \end{pr}

 \begin{rem} {\rm The statement ii) of Proposition~\ref{p8.1} was stated
 as a conjecture by Gupta \cite{Gu9}, and has been
 proved by G.-N.~Han \cite{GH}, and A.N.~Kirillov \cite{Kir1}. The proof
 given in \cite{GH} was based on the explicit combinatorial construction
 of an embedding
 $$STY(\ld ,\mu)\hookrightarrow STY(\ld\oplus (a),\mu\oplus (a))
 $$
 of the set $STY(\ld ,\mu)$ of semistandard Young tableaux of shape $\ld$
 and weight $\mu$ to that $STY(\ld\oplus (a),\mu\oplus (a))$ of semistandard
 tableaux of shape $\ld\oplus (a)$ and weight $\mu\oplus (a)$. The proof
 given in \cite{Kir1} was based on the theory of rigged configurations. It is
 still an open question whether or not a bijection constructed by
 G.-N.~Han \cite{GH} is compatible with the rigged configurations bijection.}
\end{rem}

  We are going to describe the rational functions $Z_{\ld\mu}^{\nu}(q)$
 and $Y_{\ld\mu}^{\nu}(q)$ in a particular case when the partition $\nu$
 consists of one part, namely, $\nu =(a)$.

 Without lost of generality, one can assume that partitions $\ld$ and
 $\mu$ have the following forms
 \begin{eqnarray}
 \ld &=&(\al_1+a,\ldots ,\al_r+a,\underbrace{a,\ldots ,a}_{k-r-s},
 a-\beta_s,a-\beta_{s-1},\ldots ,a-\beta_1),\label{8.2}\\
 \mu &=&(\gamma_1+a,\ldots ,\gamma_l+a,\underbrace{a,\ldots ,a}_{k-l-p},
 a-\delta_p,a-\delta_{p-1},\ldots ,a-\delta_1),\label{8.3}
 \end{eqnarray}
where $\al =(\al_1,\ldots ,\al_r)$, $\beta =(\beta_1,\ldots ,\beta_s)$,
 $\gamma =(\gamma_1,\ldots ,\gamma_l)$, and $\delta =(\delta_1,\dots
 ,\delta_p)$ are partitions such that $|\al|+|\delta|=|\beta|+|\gamma|$,
 $l(\al)+l(\beta)=l(\gamma)+l(\delta)\ge k$.

 We start with a study of the simplest case when $\gamma =\delta
 =\emptyset$, i.e. $\mu =(a^k)$. In this case one can use
 Corollaries~\ref{c4.4} and \ref{c6.7} to conclude that
 $$\lim_{n\to\infty}K_{\ld\oplus (a^n),(a^{n+k})}(q)
 =\lim_{n\to\infty}s_{\al}*s_{\beta}(q,\ldots ,q^{n+k-1})
 =q^{|\al|}s_{\al}*s_{\beta}(1,q,q^2,\ldots ),
 $$
 where $\al =((\ld_1-a)_+,(\ld_2-a)_+,\ldots )$, $\beta
 =((a-\ld_k)_+,(a-\ld_{k-1})_+,\ldots )$, and the symbol $(a)_+$
 stands for $\max (a,0)$.

 The principal specialization $s_{\al}*s_{\beta}(1,q,q^2,\ldots )$ has
 been computed by Stanley \cite{St1}. Combining with \cite{Ma},
 Chapter~VI, \S 8, Exercise~3, the Stanley result may be stated as follows
 $$s_{\al}*s_{\beta}(1,q,q^2,\ldots )=\frac{K_{\beta\al}(q,q)}{H_{\al}(q)}.
 $$
 Summarizing, we obtain the following result
 \begin{equation}
 \lim_{n\to\infty}K_{\ld\oplus (a^n),(a^{n+k})}(q)=q^{|\al|}
 \frac{K_{\beta\al}(q,q)}{H_{\al}(q)}, \label{8.4}
 \end{equation}
 where $H_{\al}(q)$ denotes the hook polynomial, see e.g. \cite{Ma},
 Chapter~I, \S 3, Example~3, and $K_{\beta\al}(q,q)$ denotes the
 specialization $q=t$ of the double Kostka--Macdonald polynomial
 $K_{\beta\al}(q,t)$, \cite{Ma}, Chapter~VI, \S 8.

 As for the limit $n\to\infty$ of the {\it cocharge} Kostka--Foulkes
 polynomials ${\overline K}_{\ld\oplus (a^n), (a^{n+k})}(q)$, it follows
 from Corollary~\ref{c4.4} that
 \begin{equation}
 \lim_{n\to\infty}{\overline K}_{\ld\oplus (a^n), (a^{n+k})}(q)=
 \frac{{\overline K}_{\beta ,(1^{|\beta|})}(q)}{H_{\al}(q)}. \label{8.5}
 \end{equation}

 The next case we are going to consider is the case when $\gamma =\emptyset$,
 i.e. $\mu_1\le a$.

 \begin{theorem}\label{t8.1}  Let $\ld$ and $\mu$ be partitions of the
 forms (\ref{8.2}) and (\ref{8.3}) correspondingly, and assume
 additionally that $\mu_1\le a$. Then
 \begin{eqnarray}
 \lim_{n\to\infty}K_{\ld\oplus (a^n),\mu\oplus (a^n)}(q)&=&
 \frac{q^{|\al|}}{H_{\al}(q)}\sum_{\eta}K_{\eta\al}(q,q)
 K_{\beta\setminus\eta ,\delta}(q);\label{8.6}\\
 \lim_{n\to\infty}\wt K_{\ld\oplus (a^n),\mu\oplus (a^n)}(q)&=&
 \frac{\wt K_{\beta ,(\delta ,1^{|\al|})}(q)}{H_{\al}(q)}.\label{8.7}
 \end{eqnarray}
 \end{theorem}

 \begin{ex}\label{e8.4} {\rm Take $\ld =(7,3,2)$, $\mu =(4,3,2,1,1,1)$
 and $a=(4)$. Then we have $\al =(3)$, $\beta =(4,4,4,2,1)$ and $\delta
 =(3,3,3,2,1)$. The Kostka--Foulkes polynomials $K_{\ld\oplus
 (a^n),\mu\oplus (a^n)}(q)$ have been computed in \cite{Kir6}, \S 7,
 Example~6, as well as the limits
 \begin{eqnarray*}
 &&\lim_{n\to\infty}K_{\ld\oplus (a^n),\mu\oplus (a^n)}(q)=
 q^5\frac{6+11q+14q^2+13q^3+9q^4+3q^5+2q^6}{(1-q)^2(1-q^3)},\\
&&\lim_{n\to\infty}\wt K_{\ld\oplus (a^n),\mu\oplus (a^n)}(q)=\\
&&\frac{1+q+3q^2+5q^3+7q^4+9q^5+11q^6+8q^7+7q^8+4q^9+q^{10}
+q^{11}}{(1-q)^2(1-q^3)}.
\end{eqnarray*}

First of all, the hook polynomial
$$H_{\al}(q)=(1-q)(1-q^2)(1-q^3),$$
and using the fermionic formula (\ref{4.1a}) for Kostka--Foulkes
polynomials, it is not difficult to check that
\begin{eqnarray*}
K_{\beta ,(\delta ,1,1,1)}(q)&=&q^3(1,2,5,11,15,19,20,16,12,8,4,2,1)\\
&=&q^3(1+q)(1,1,4,7,8,11,9,7,5,3,1,1).
\end{eqnarray*}
Hence,
$$\lim_{n\to\infty}\wt K_{\ld\oplus (a^n),\mu\oplus (a^n)}(q)=
\frac{\wt K_{\beta ,(\delta ,1,1,1)}(q)}{H_{\al}(q)},
$$
as it has to be according to Theorem~\ref{t8.1}, (\ref{8.7}).

On the other hand,
\begin{eqnarray*}&&\sum_{\eta}K_{\beta\setminus\eta
,\delta}(q)K_{\eta\al}(q,q)=K_{(44421)\setminus
(21),\beta}(q)K_{(21),\al}(q,q)\\
&+& K_{(44421)\setminus (3),\beta}(q)K_{(3),\al}(q,q)
+K_{(44421)\setminus (1^3),\beta}(q)K_{(1^3),\al}(q,q).
\end{eqnarray*}
Now one can check that
\begin{eqnarray*}
K_{(44421)\setminus (21),\beta}(q)&=&q(3,8,10,10,6,2,1)
,~~~K_{(21),\al}(q,q)=q+q^2;\\
K_{(44421)\setminus (3),\beta}(q)&=&q^2(3,5,5,3,1),~~~
K_{(3),\al}(q,q)=1;\\
K_{(44421)\setminus (1^3),\beta}(q)&=&(1,2,4,5,4,2,1),~~~
K_{(1^3),\al}(q,q)=q^3.
\end{eqnarray*}
Therefore,
\begin{eqnarray*}
\sum_{\eta}K_{\beta\setminus\eta ,\delta}(q)K_{\eta\al}(q,q)&=&
q^2(6,17,25,27,22,12,5,2)\\
&=&q^2(1+q)(6,11,14,13,9,3,2),
\end{eqnarray*}
as it should be according to Theorem~\ref{t8.1}, (\ref{8.6}).}
\end{ex}

 Finally, we consider the general case:
 \begin{theorem}\label{t8.2} Let $\ld$ and $\mu$ be partitions of the forms
 (\ref{8.2}) and (\ref{8.3}) respectively, then
 \begin{eqnarray*}
 &&\lim_{n\to\infty}K_{\ld\oplus (a^n),\mu\oplus (a^n)}(q)\bdoteq \\
 &&\frac{1}{(q;q)_{|\al|-|\gamma|}}\sum_{\eta_1,\eta_2}
 K_{\beta\setminus\eta_1,\delta}(q){\overline K}_{\al\setminus\eta_2,\gamma}
 (q)K_{\eta_1,\eta_2}(q,q)K_{\eta_2,(1^{|\al|+|\gamma|})}(q);\\
 &&\lim_{n\to\infty}{\overline K}_{\ld\oplus (a^n),\mu\oplus (a^n)}(q)
 \bdoteq \frac{1}{(q;q)_{|\al|-|\gamma|}}
 {\overline K}_{\al,(\gamma\oplus(1^{|\al|-|\gamma|}))}(q)~
 {\overline K}_{\beta,(\delta\oplus(1^{|\beta|-|\delta|}))}(q).
 \end{eqnarray*}
 \end{theorem}

\newpage
\hskip -0.6cm {\bf Exercises} \vskip 0.2cm

{\bf 1.} Let $\ld$ and $\mu$ be partitions, define $\ld +\mu$ to
be the sum of the sequences $\ld$ and $\mu$:
$$(\ld +\mu)_i=\ld_i+\mu_i.
$$
The operations + and $\oplus$ are dual to each other, i.e.
$(\ld\oplus\mu)'=\ld'+\mu'$.

Show that if $\nu$ is a partition, then
\begin{equation}
K_{\ld +\nu,\mu +\nu}(q)\ge K_{\ld ,\mu}(q). \label{8.8}
\end{equation}
\vskip 0.3cm

\begin{con}\label{con8.6} Let $\ld ,\mu$ and $\nu$ be partitions,
then
\begin{equation}
K_{\ld\oplus\nu ,\mu\oplus\nu}(q,t)\ge K_{\ld ,\mu}(q,t),
\label{8.9}
\end{equation}
i.e. the difference $K_{\ld\oplus\nu ,\mu\oplus\nu}(q,t)-K_{\ld
,\mu}(q,t)$ is a polynomial in $q$ and $t$ with non--negative
(integer) coefficients.
\end{con}

Here $K_{\ld ,\mu}(q,t)$ stands for the double Kostka polynomial
introduced by I.~Macdonald, \cite{Ma}, Chapter~VI, (8.11). Note
that Conjecture~\ref{con8.6} gives a common (conjectural)
generalization of inequalities (\ref{8.1}) and (\ref{8.8}).
\vskip 0.2cm

\hskip -0.6cm {\bf Question.} Does there exist the limit
$\ds\lim_{n\to\infty}K_{\ld\oplus\nu^n,\mu\oplus\nu^n}(q,t)$ ?

\vskip 0.2cm {\bf 2.} For given integers $n\ge 1$ and $l\ge 1$,
let $\ld$ and $\mu$ be partitions such that $|\ld|\equiv
|\mu|({\rm mod}n)$ and $\mu_1\le l$. Consider two sequences of
partitions
$$\ld_{L}:=\ld
+(\left(Ll+\ds\frac{|\mu|-|\ld|}{n}\right)^n)~~{\rm and}~~
\mu_L:=((l)^{nL},\mu).
$$

Show that there exist the limits
$\ds\lim_{L\to\infty}K_{\ld_L,\mu_L}(q)$ and
$\ds\lim_{L\to\infty}{\overline K}_{\ld_L,\mu_L}(q)$, which are
formal power series in $q$, but both do not equal to any rational
function.

\begin{prb} Find combinatorial and representation theoretical
interpretations of the latter power series.
\end{prb}
For partial answer on this problem see \cite{HKKOTY,Kir4,NY}.

\vskip -0.5cm

\section*{Acknowledgments}
This notes grew out of the series of lectures given at the
Hokkaido, Kyushu, and Nagoya Universities, and at the RIMS and
IIAS (Kyoto). I would like to thank those who attended for helpful
comments, suggestions and support. In particular M.~Kashiwara,
T.~Miwa, H.~Umemura, H.-F.~Yamada and M.~Yoshida. My special
thanks go to Alexander Postnikov for very helpful discussions. I
wish thank with much gratitude the colleagues at the Graduate
School of Mathematics, Nagoya University, for their hospitality
which made it possible to finish this work.

I would like to acknowledge my special indebtedness to
Dr.~N.A.~Liskova for the inestimable help and support on all
stages of the paper writing.


\end{document}